\documentclass[11pt]{article}
\usepackage{fullpage, amsmath, amssymb, amsthm, epsfig, color, graphicx}

\def\allow{\mathop{\rm Allow}\nolimits}
\def\qallow{\mathop{\rm q\!-\!Allow}\nolimits}
\def\Int{\mathop{\rm Int}\nolimits}

\def\eqbd{\mathop{{:}{=}}}
\def\bdeq{\mathop{{=}{:}}}
\def\R{\mathbb{R}}
\def\S{\mathcal{S}}
\def\G{\mathcal{G}}
\def\C{\mathbb{C}}
\def\Z{\mathbb{Z}}
\def\k{{\cal K}}
\def\i{{\cal I}}

\newtheorem{theorem}{Theorem}[section]

\newtheorem{result}{Result}[theorem]
\newtheorem{proposition}[theorem]{Proposition}
\newtheorem{remark}[theorem]{Remark}
\newtheorem{lemma}[theorem]{Lemma}
\newtheorem{corollary}[theorem]{Corollary}

\newtheorem{definition}[theorem]{Definition}
\newtheorem{example}[theorem]{Example}
\newcommand{\ignore}[1]{}

\title{Toward accurate polynomial evaluation \\ 
in rounded arithmetic}
\author{James Demmel\thanks{Mathematics Department and CS Division,
University of California, Berkeley, CA 94720. 
The author acknowledges the support of NSF under grants CCF-0444486, ACI-00090127, 
CNS-0325873 and of DOE under grant DE-FC02-01ER25478.},
Ioana Dumitriu\thanks{Mathematics Department, University of California, Berkeley, CA
94720. The author acknowledges the support of the Miller Institute for Basic Research 
in Science.},
and Olga Holtz\thanks{Mathematics Department, University of California, Berkeley, CA 
94720.}}

\begin{document}
\date{November 1, 2005}
\maketitle

\begin{abstract}
Given a multivariate real (or complex) polynomial $p$ and a domain $\cal D$,
we would like to decide whether an algorithm exists to evaluate $p(x)$ accurately 
for all $x \in {\cal D}$ using rounded  real (or complex) arithmetic. 
Here ``accurately'' means with relative error less than 1, i.e., with some correct 
leading digits. The answer depends on the model of rounded arithmetic:
We assume that for any arithmetic operator  $op(a,b)$, for example $a+b$ or 
$a \cdot b$,  its computed value is $op(a,b) \cdot (1 + \delta)$, 
where $| \delta |$ is bounded by some constant $\epsilon$ where $0 < \epsilon \ll 1$, 
but $\delta$ is otherwise arbitrary. This model is the traditional one used to analyze 
the accuracy of floating point algorithms.


Our ultimate goal is to establish a decision procedure that, for any $p$ and $\cal D$, 
either exhibits an accurate algorithm or proves that none exists. In contrast to the 
case where numbers are stored and manipulated as finite bit strings (e.g., as floating 
point numbers or rational  numbers)  we show that some polynomials $p$ are impossible to 
evaluate accurately.  The existence of an accurate algorithm will depend not just
on $p$ and $\cal D$, but on which arithmetic operators are available (perhaps beyond 
$+$, $-$, and $\cdot$), which constants are available to the algorithm (integers, 
algebraic numbers, ...), and whether branching is permitted in the algorithm. For 
floating point computation,  our model can be used to identify which accurate operators
beyond $+$, $-$ and $\cdot$ (e.g., dot products, 3x3 determinants, ...) are necessary to 
evaluate a particular $p(x)$.

Toward this goal, we present necessary conditions on $p$ for it to be 
accurately evaluable on open real or complex domains ${\cal D}$.
We also give sufficient conditions, and describe progress toward
a complete decision procedure. We do present a complete 
decision procedure for homogeneous polynomials $p$ with integer coefficients,
${\cal D} = \C^n$, and using only the arithmetic operations
$+$, $-$ and $\cdot$.


\end{abstract} 

\newpage

\tableofcontents

\section{Introduction}
\label{sec_Intro}

Let $x = (x_1,...,x_n)$ be a vector of real (or complex) numbers,
let $p(x)$ denote a multivariate polynomial, and let
${\cal D}$ be a subset of $\R^n$ (or $\C^n$).
We would ideally like to evaluate $p(x)$ 
{\em accurately} for all $x \in {\cal D}$,
despite any rounding errors in arithmetic operations.
The nature of the problem depends on 
how we measure accuracy, 
what kinds of rounding errors we consider,
the class of polynomials $p(x)$,
the domain $\cal D$, 
and what operations and constants our evaluation algorithms may use.
Depending on these choices, an accurate algorithm for evaluating 
$p(x)$ may or may not exist. 
Our ultimate goal is a decision procedure that will either exhibit 
an algorithm that evaluates $p(x)$ accurately for all $x \in \cal D$, 
or else exhibit a proof that no such algorithm exists.

By {\em accuracy}, we mean that we compute an approximation
$p_{comp}(x)$ to $p(x)$ that has small {\em relative} error:
$|p(x) - p_{comp}(x)| \leq \eta |p(x)|$ for some desired $0 < \eta < 1$.
In particular, $\eta < 1$ implies that $p(x)=0$ if and only if
$p_{comp}(x) = 0$. This requirement that $p$ and $p_{comp}$ define
the same variety will be crucial in our development.
We justify this definition of accuracy in more detail in
Section~2.

Our motivation for this work is two-fold. First, it is common for
numerical analysts to seek accurate formulas for particularly
important or common expressions.
For example, in computational geometry and mesh
generation, certain geometric predicates like ``Is point $x$
inside, on or outside circle $C$?'' are expressed as multivariate
polynomials $p( \cdot )$  whose signs determine the answer; 
the correctness
of the algorithms depends critically on the correctness of the
sign, which is in turn guaranteed by having a relative error
less than 1 in the value of $p( \cdot )$ \cite{shewchuk}.
We would like to automate the process of finding such formulas.

The second motivation is based on recent work of Koev and one of the
authors \cite{demmelkoevICIAM03} which identified several 
classes of structured 
matrices (e.g., Vandermonde, Cauchy, totally positive, certain discretized
elliptic partial differential equations, ...) for which 
algorithms exist to accurately perform some (or all) computations from linear
algebra: determinants, inversion, Gaussian elimination, computing singular 
values,
computing eigenvalues, and so on. The proliferation of these classes of
structured matrices led us to ask what common algebraic structure these 
matrix classes possess that made these accurate algorithms possible.
This paper gives a partial answer to this question; 
see section~\ref{sec_LinearAlgebra}.

Now we consider our model of rounded arithmetic. Let $op(\cdot)$ denote
a {\em basic arithmetic operation}, for example $op(x,y) = x+y$ or
$op(x,y,z) = x+y \cdot z$. Then we assume that the rounded
value of $op(\cdot)$, which we denote $rnd(op(\cdot))$, satisfies
\begin{equation}\label{def_roundingmodel}
rnd(op(\cdot)) = op(\cdot)(1 + \delta)
\end{equation}
where we call $\delta$ the {\em rounding error}.
We assume only that $| \delta |$ is tiny, $| \delta | \leq \epsilon$,
where $0 < \epsilon < 1$ and typically $\epsilon \ll 1$;
otherwise $\delta$ is an arbitrary real (or complex) number.
The constant $\epsilon$ is called the {\em machine precision},
by analogy to floating point computation, since this model
is the traditional one used to analyze the accuracy of
floating point algorithms \cite{higham96,wilkinsonroundingerror}.

To illustrate the obstacles to accurate evaluation that this
model poses, consider evaluating $p(x) = x_1 + x_2 + x_3$
in the most straightforward way: add (and round) $x_1$ and $x_2$, 
and then add (and round) $x_3$. If we let $\delta_1$ be the
first rounding error and $\delta_2$ be the second rounding error,
we get the computed value 
$p_{comp}(x) = ((x_1+x_2)(1+\delta_1) + x_3)(1 + \delta_2)$.
To see that this algorithm is not accurate, simply choose
$x_1 = x_2 = 1$ and $x_3=-2$ (so $p(x)=0)$ and
$\delta_1 \neq 0$. Then
$p_{comp}(1,1,-2) = 2\delta_1 (1 + \delta_2) \neq 0$, so the relative error is infinite.
Indeed, it can be shown that there is an open set of $x$
and of $(\delta_1, \delta_2)$ where the relative error is
large, so that this loss of accuracy occurs on a ``large'' set.
We will see that unless $x_1 + x_2 + x_3$ is itself a 
basic arithmetic operation, or unless the variety
$\{ x_1 + x_2 + x_3 = 0 \}$ is otherwise constructible from
varieties derived from basic operations as described 
in Theorem~\ref{gen_result}, then {\em no\/} algorithm exists to 
evaluate $x_1 + x_2 + x_3$ accurately for all arguments.

\ignore{
Here is what can be said about the errors in straightforward polynomial
evaluation using  Horner's rule \cite[Chap. 5]{higham96}: The
computed value of $p(x)$ is the exact value of $\hat{p}(x)$ where
the polynomial $\hat{p}$ differs from $p$ by making a small relative
change of $O( \epsilon )$ in each of $p$'s coefficients. This
in turn lets us bound $|p(x) - \hat{p}(x)|  = O( \epsilon ) |p|(|x|)$,
where $|p|(|x|)$ means the value of the polynomial where absolute values
are taken of each coefficient of $p$ and each component of $x$. This does not
guarantee high relative accuracy when $|p(x)| \ll |p|(|x|)$, 
i.e., at or near a near a zero of $p(x)$,
though it is adequate
in many instances, especially where one may consider the coefficients
of $p$ to be known only approximately; see Section~\ref{sec_OtherModels}
for further discussion.
} 

In contrast, if we were to assume that the $x_i$ and coefficients of $p$ were given 
as exact rational numbers  (e.g., as floating point numbers), then by performing 
integer arithmetic with sufficiently large integers  it would clearly be a 
straightforward matter to evaluate any $p(x)$  as an exact rational number.  (One 
could also use floating point arithmetic to accomplish this; see 
Sections~\ref{sec_BasicArithmeticOperations} and~\ref{sec_OtherModels}.)  In other 
words, accurate evaluation is always possible, and the only question is cost. Our 
model addresses this by identifying which composite operations have to be provided 
with high precision in order to evaluate $p(x)$ accurately.
For further discussion of the challenge of evaluating a simple polynomial
like $x_1 + x_2 + x_3$ accurately, see section~\ref{sec_OtherModels}.


\ignore{ 
The domain of evaluation $\cal D$ may in principle be any semialgebraic subset
of ${\R}^n$ or ${\C}^n$. 
Some choices, for example ${\cal D} = \{x: p(x) = 0\}$, make
evaluating $p(x)$ easy but beg the question of how one would
decide whether $x \in \cal D$.
For simplicity we 
will mostly consider ${\cal D} = {\R}^n$ or ${\C}^n$, and point out 
where our theory works for smaller $\cal D$
(but with ``large'' interiors).
Section 2 discusses this point further.
} 

We give some examples to illustrate our results.
Consider the family of homogeneous polynomials
\[
M_{jk}(x) = j \cdot x_3 ^6 + x_1^2 \cdot x_2^2 \cdot 
(j \cdot x_1^2 + j \cdot x_2^2 - k \cdot x_3^2)
\]
where $j$ and $k$ are positive integers, 
${\cal D} = \R^n$, and
we allow only addition, subtraction and multiplication
of two arguments as basic arithmetic operations, along
with comparisons and branching.
\begin{itemize}
\item
When $k/j < 3$, $M_{jk}(x)$ is {\em positive definite,\/} i.e.,
zero only at the origin and positive elsewhere.  This will mean that 
$M_{jk}(x)$ is easy to evaluate accurately using a simple method
discussed in Section~\ref{sec_PositivePolys}.
\item
When $k/j > 3$, then we will show that $M_{jk}(x)$
cannot be evaluated accurately by {\em any\/} algorithm
using only addition, subtraction and multiplication of
two arguments. This will follow from a simple necessary
condition on the real variety $V_{\R}(M_{jk})$, 
the set of real $x$ where $M_{jk}(x)=0$, see Theorem~\ref{conj}.
This theorem requires the real variety $V_{\R}(p)$ 
(or the complex variety $V_{\C}(p)$) 
to lie in a certain explicitly given finite set of varieties 
called {\em allowable varieties}
in order to be able to 
evaluate $p(x)$ accurately in real arithmetic (or in complex
arithmetic, resp.).
\item
When $k/j=3$, i.e., on the boundary between the above two cases,
$M_{jk}(x)$ is a multiple of the Motzkin polynomial
\cite{reznick2000}. 
Its real variety $V_{\R}(M_{jk}) = \{x: |x_1| = |x_2| = |x_3| \}$
satisfies the necessary condition of Theorem~\ref{conj}, and the simplest
accurate algorithm to evaluate it that we know is shown (in part) below:
\begin{eqnarray*}
{\rm if} &   & |x_1-x_3| \leq |x_1+x_3| \wedge |x_2-x_3| \leq |x_2+x_3| \; \; \; {\rm then}  \\
       & p = & x_3^4 \cdot [4((x_1-x_3)^2 + (x_2-x_3)^2 + (x_1-x_3)(x_2-x_3))] \\
         &   & +x_3^3 \cdot [2(2(x_1-x_3)^3 + 5(x_2-x_3)(x_1-x_3)^2 + 5(x_2-x_3)^2(x_1-x_3) + \\
         &   &   \; \; \; \hspace*{.5in} 2(x_2-x_3)^3)] \\
         &   & +x_3^2 \cdot [(x_1-x_3)^4 + 8(x_2-x_3)(x_1-x_3)^3 + 9(x_2-x_3)^2(x_1-x_3)^2 + \\
         &   &  \; \; \; \hspace*{.5in} 8(x_2-x_3)^3 (x_1-x_3) + (x_2-x_3)^4] \\
         &   & +x_3 \cdot [2(x_2-x_3)(x_1-x_3)((x_1-x_3)^3 + 2(x_2-x_3)(x_1-x_3)^2 + \\
         &   &  \; \; \; \hspace*{.5in} 2(x_2-x_3)^2(x_1-x_3) + (x_2-x_3)^3)] \\
         &   & + (x_2-x_3)^2(x_1-x_3)^2((x_1-x_3)^2 + (x_2-x_3)^2) \\
       & p = & j \cdot p \\
{\rm else} &  &  ... \; \; 7 \; \;
         {\rm more\ analogous\ cases.}
\end{eqnarray*}
In general, for a Motzkin polynomial in $n$ real variables ($n=3$ above), 
the algorithm has $2^n$ separate cases. Just $n$ tests and branches
are needed to choose the correct case for any input $x$,
so that the cost of running the algorithm
is still just a polynomial function of $n$ for any particular $x$.
\end{itemize}

In contrast to the real case, when ${\cal D} = \C^n$ then
Theorem~\ref{conj} will show that $M_{jk}(x)$ is not 
accurately evaluable using only addition, subtraction and
multiplication. 

If we still want to evaluate $M_{jk}(x)$ accurately in
one of the cases where addition, subtraction and multiplication
alone do not suffice, it is natural to ask which
composite or ``black-box'' operations we
would need to implement accurately to do so.
Section~\ref{bb} addresses this question.

The necessary condition for accurate evaluability of $p(x)$ 
in Theorem~\ref{conj} depends only on the variety of $p(x)$.
The next example shows that the variety alone is not enough
to determine accurate evaluability, at least in the real case.
Consider the two irreducible, homogeneous, degree $2d$, real polynomials
\begin{equation}
\label{eqn_NeedInduction}
p_i(x) = (x_1^{2d} + x_2^{2d}) + (x_1^2 + x_2^2)(q_i(x_3,...,x_n))^2
\; \; \; {\rm for} \; i=1,\;2
\end{equation}
where $q_i(\cdot)$ is a homogeneous polynomial of degree $d-1$.
Both $p_1(x)$ and $p_2(x)$ have the same real variety
$V_{\R}(p_1) = V_{\R}(p_2) = \{x: x_1 = x_2 = 0 \}$, which
is allowable, i.e.,  satisfies the necessary condition for 
accurate evaluability in Theorem~\ref{conj}. However, near
$x_1=x_2=0$, $p_i(x)$ is ``dominated'' by 
$(x_1^2 + x_2^2)(q_i(x_3,...,x_n))^2$, so accurate
evaluability of $p_i(x)$ in turn depends on
accurate evaluability of $q_i(x_3,...,x_n)$.
Since $q_1( \cdot )$ may be accurately evaluable while $q_2 ( \cdot )$
is not, we see that $V_{\R}(p_i)$ alone cannot determine whether
$p_i(x)$ is accurately evaluable. Applying the same principle to
$q_i( \cdot )$, we see that any decision procedure must be recursive,
expanding $p_i(x)$ near the components of its variety and so on.
We show current progress toward a decision procedure
in Section~\ref{suf_r}. In particular, Theorem~\ref{p=>pdom} shows that, at least
for algorithms without branching, being able to compute dominant terms
of $p$ (suitably defined) accurately on $\R^n$ is a necessary condition for
computing $p$ accurately on $\R^n$.   Furthermore, Theorem~\ref{pdom=>p}  shows that
accurate evaluability of the dominant terms, along with branching,  is
sufficient to evaluate $p$ accurately.

In contrast to the real case,  Theorem~\ref{sufficiency_c} shows that for 
the complex case, with ${\cal D} = \C^n$, and using only addition,
subtraction and multiplication of two arguments, a
homogeneous polynomial $p(x)$ with integer coefficients is
accurately evaluable if and only if it satisfies the necessary
condition of Theorem~\ref{conj}. More concretely, $p(x)$ is accurately
evaluable for all $x \in \C^n$ if and only if $p(x)$ can
be completely factored into a product of factors of the form
$x_i$, $x_i + x_j$ and $x_i - x_j$.

The results described so far from Section~\ref{class} consider
only addition, subtraction, multiplication and (exact) negation
(which we call {\em classical arithmetic}).
Section~\ref{bb} considers the same questions
when accurate {\em black-box} operations beyond addition, 
subtraction and multiplication are permitted, 
such as fused-multiply-add \cite{FMA}, 
or indeed any collection of polynomials at all
(e.g., dot products, 3x3 determinants, ...).
The necessary condition on the variety of $p$ from 
Theorem~\ref{conj} is generalized to black-boxes in
Theorem~\ref{gen_result}, and the sufficient conditions 
Theorem~\ref{sufficiency_c} 
in the complex case are generalized in
Theorems~\ref{q-suff-c1} and \ref{q-suff-c2}.

The rest of this paper is organized as follows.
Section~\ref{sec_Models} discusses further details of our 
algorithmic model, explains why it is a useful model of
floating point computation, and otherwise justifies the choices 
we have made in this paper.  
Section~\ref{sec_PositivePolys} discusses the evaluation
of positive polynomials.
Section~\ref{class} discusses 
necessary conditions (for real and complex data) and
sufficient conditions (for complex data)
for accurate evaluability, when using only 
classical arithmetic.
Section~\ref{suf_r} describes progress toward
devising a decision procedure for accurate evaluability
in the real case using classical arithmetic.
Section~\ref{bb} extends Section~\ref{class}'s 
necessary conditions to
arbitrary black-box arithmetic operations,
and gives sufficient conditions in the complex case.
Section~\ref{sec_LinearAlgebra} describes implications for 
accurate linear algebra
on structured matrices.

\newpage

\section{Models of Algorithms and Related Work}
\label{sec_Models}

Now we state more formally our decision question.
We write the output of our algorithm as $p_{comp}(x, \delta )$, where
$\delta = ( \delta_1  , \delta_2, ... \delta_k )$ is the
vector of rounding errors made during the algorithm.

\begin{definition}
We say that $p_{comp}(x, \delta)$ is an {\em accurate algorithm}
for the evaluation of $p(x)$ for $x \in {\cal D}$ if
\\
\hspace*{0.25in}    $\forall \; 0 < \eta < 1 \; \; \; \; \;$ ...
                 for any $\eta$ = desired relative error \\
\hspace*{0.50in} $\exists \; 0 < \epsilon < 1 \; \; \; \;$ ...
                 there is an $\epsilon$ = machine precision \\
\hspace*{0.75in} $\forall \; x \in {\cal D} \; \; \; \; \;$ ...
                 so that for all $x$ in the domain \\
\hspace*{1.00in} $\forall \; |\delta_i| \leq \epsilon \; \; \;$ ...
                 and for all rounding errors bounded by $\epsilon$ \\
\hspace*{1.25in} $|p_{comp}(x,\delta) - p(x)| \leq \eta \cdot |p(x)|$ ...
                 the relative error is at most $\eta$.
\\
\end{definition}

Our ultimate goal is a decision procedure (a ``compiler'') 
that takes $p(\cdot)$ and $\cal D$ as input, and
either produces an accurate algorithm $p_{comp}$
(including how to choose the machine precision $\epsilon$ given 
the desired relative error $\eta$)
or exhibits a proof that none exists.

To be more precise, we must say what our set of possible
algorithms includes.
The above decision question is apparently not Tarski-decidable
\cite{renegar92,tarski_book}
despite its appearance, because we see no way to
express ``there exists an algorithm'' in that format.

The basic decisions about algorithms that we make
are as follows, with details given in the 
indicated sections:

\begin{description}

\item[Sec.~\ref{sec_ExactRoundedInputs}:]
We insist that the inputs $x$ are given exactly, rather than
approximately.

\item[Sec.~\ref{sec_FiniteConvergence}:] 
We insist that the algorithm compute the exact value of $p(x)$ in finitely 
many steps when all rounding errors $\delta = 0$. 
In particular, we exclude iterative algorithms which 
might produce an approximate value of $p(x)$ even when
$\delta = 0$.

\item[Sec.~\ref{sec_BasicArithmeticOperations}:]
We describe the basic arithmetic operations we consider, 
beyond addition, subtraction and multiplication.
We also describe the constants available to our algorithms.

\item[Sec.~\ref{sec_ComparisonsBranching}:]
We consider algorithms both with and without comparisons and
branching, since this choice may change the set of polynomials
that we can accurately evaluate.


\item[Sec.~\ref{sec_Nondeterminism}:]
If the computed value of an operation depends only the
values of its operands, i.e., if the same operands $x$ and $y$
of $op(x,y)$ always yield the same $\delta$ in
$rnd(op(x,y)) = op(x,y) \cdot (1+ \delta)$, then we
call our model {\em deterministic,\/} else it is
{\em nondeterministic.\/} We show that comparisons and branching
let a nondeterministic machine simulate a deterministic one,
and subsequently restrict our investigation to the easier 
nondeterministic model. 

\item[Sec.~\ref{sec_Domain}:]
What domains of evaluation $\cal D$ do we consider?
In principle, any semialgebraic set $\cal D$ is a
possibility, but for simplicity we mostly consider
open ${\cal D}$, especially ${\cal D} = \R^n$ or ${\cal D} = \C^n$.
We point out issues in extending results to other $\cal D$.
\end{description}

Finally, Section~\ref{Sec_Axioms} summarizes the axioms our model satisfies, 
and Section~\ref{sec_OtherModels} compares our model to other models of arithmetic, 
and explains the advantages of our model.

\subsection{Exact or Rounded Inputs}
\label{sec_ExactRoundedInputs}
We must decide whether we assume that the arguments
are given exactly 
\cite{BCSS96,CuckerGrigoriev99,renegar92,tarski_book}
or are known only approximately
\cite{cuckersmale99,EdelatSunderhauf98,KoKI,pourelrichards89}.
Not knowing the input $x$ exactly means that at best 
(i.e., in the absence of any further error) we could 
only hope to compute the exact value of $p(\hat{x})$
for some $\hat{x} \approx x$, an algorithmic property 
known as {\em backward stability} \cite{demmelMA221,higham96}.
Since we insist that zero outputs be computed exactly
in order to have bounded relative error, this means there 
is no way to guarantee that $p(\hat{x}) = 0$
when $p(x)=0$, for nonconstant $p$. This is true even for simple
addition $x_1+x_2$.  So we insist on exact inputs in our model.

\subsection{Finite Convergence}
\label{sec_FiniteConvergence}
Do we consider algorithms that take a bounded amount of
time for all inputs $x \in {\cal D}$, and return 
$p_{comp}(x,0) = p(x)$, i.e., the exact answer 
when all rounding errors are zero?
Or do we consider possibly iterative algorithms that might
take arbitrarily long on some inputs to produce an adequately
accurate answer? We consider only the former,
because (1) it seems natural to use a finite algorithm to evaluate
a finite object like a polynomial, (2) we have seen
no situations where an iterative algorithm offers any
advantage to obtaining guaranteed {\em relative} accuracy
and (3) this lets us write any algorithm as a 
piecewise polynomial function and so use
tools from algebraic geometry.

\subsection{Basic Arithmetic Operations and Constants}
\label{sec_BasicArithmeticOperations}

What are the basic arithmetic operations? For most of
the paper we consider addition, subtraction and multiplication
of two arguments, since this is necessary and sufficient for
polynomial evaluation in the absence of rounding error.
Furthermore, we consider negation as a basic operation that
is always exact (since this mimics all implementations
of rounded arithmetic).
Sometimes we will also use (rounded) multiplication by a constant 
$op(x) = c \cdot x$.
We also show how to extend our results 
to include additional basic arithmetic operations like 
$op(x,y,z) = x + y \cdot z$.  
The motivations for considering such additional
``black-box'' operations are as follows:

\begin{enumerate}

\item By considering operations like $x+c$, $x-c$ and $c \cdot x$ 
for any $c$ in a set $C$ of constants, we may investigate how
the the choice of $C$ affects accurate evaluability.
For example, if $C$ includes the roots of a polynomial like 
$p(x) = x^2-2$, then we can accurately evaluate 
$p(x)$ with the algorithm 
$(x-\sqrt{2}) \cdot (x+\sqrt{2})$,
but otherwise it may be impossible.
We note that having $C$ include all algebraic numbers
would in principle let us evaluate any
univariate polynomial $p(x)$ accurately by using its
factored form $p(x) = c \prod_{i=1}^d (x-r_i)$.

In the complex case, it is natural to consider
multiplication by $\sqrt{-1}$ as an exact operation,
since it only involves ``swapping'' and possibly
negating the real and imaginary parts. We can
accommodate this by introducing operations like 
$x+\sqrt{-1} \cdot y$ and $x-\sqrt{-1} \cdot y$.

The necessary conditions in Theorem~\ref{gen_result} and
sufficient conditions in Theorems~\ref{q-suff-c1} 
and \ref{q-suff-c2} do not depend
on how one chooses an operation $x-r_i$ from a
possibly infinite set, just whether that operation
exists in the set. On the other hand, a decision
procedure must effectively choose that operation,
so our decision procedures will restrict 
themselves to enumerable (typically finite!) sets
of possible operations.\footnote{We could in principle
deal with the set of all instructions 
$x - r$ for $r$ an arbitrary algebraic number, 
because the algebraic numbers are enumerable.}

\item Many computers now supply operations like 
$x+ y \cdot z$ 
in hardware with the accuracy we demand
(the {\em fused-multiply-add} instruction \cite{FMA}).
It is natural to ask how this operation extends the class of
polynomials that we can accurately evaluate.

\item It is natural to build a library 
(in software or perhaps even hardware) 
containing several such accurate operations, 
and ask how much this extends the class of polynomials 
that can be evaluated accurately.
This approach is taken in computational geometry,
where the library of accurate operations is chosen 
to implement certain geometric predicates precisely 
(e.g., ``is point $x$ inside, outside or on circle $C$?''
written as a polynomial whose sign determines the answer).
These precise geometric predicates are critical
to performing reliable mesh generation \cite{shewchuk}.

\item A common technique for extending floating point
precision is to simulate and manipulate extra precision 
numbers by representing a high precision number $y$ as a sum 
$y = \sum_{i=1}^k y_i$ of numbers satisfying
$|y_i| \gg |y_{i+1}|$, the idea being that each $y_i$ represents
(nearly) disjoint parts of the binary expansion of $y$
(see \cite{bailey1a,bailey1b,dekker,demmelhida,KhachiyanICM84,moller,pichat,priest} 
and the references therein; similar techniques were used by Gill as early as 1951).
This technique can be modeled by the correct choice of
black-box operations as we now illustrate. 
Suppose we include the enumerable set of black-box 
operations $\sum_{i=1}^n p_i$, where $n$ is any finite
number, and each $p_i$ is the product of 1 to $d$ arguments.
In other words, we include the accurate evaluation of arbitrary
multivariate polynomials in $\Z[x]$ of degree at most $d$
among our black-box operations.
Then the following sequence of operations produces 
as accurate an approximation of any such polynomial
\[
\sum_{i=1}^n p_i = y_1 + y_2 + \cdots + y_k
\]
as desired:
\begin{eqnarray*}
y_1 & = & rnd ( \sum_{i=1}^n p_i ) = (1+ \delta_1)( \sum_{i=1}^n p_i ) \\
y_2 & = & rnd ( \sum_{i=1}^n p_i - y_1 ) = (1+\delta_2)( \sum_{i=1}^n p_i - y_1 ) \\
    & \cdots & \\
y_k & = & rnd ( \sum_{i=1}^n p_i - \sum_{j=1}^{k-1}y_j )
          = (1 + \delta_k)( \sum_{i=1}^n p_i - \sum_{j=1}^{k-1}y_j )
\end{eqnarray*}
Induction shows that 
\[
\sum_{j=1}^{k}y_j = 
\left( 1 - (-1)^k ( \prod_{j=1}^k \delta_j ) \right) \cdot (\sum_{i=1}^n p_i)
\]
so that $y = \sum_{j=1}^{k}y_j$ approximates the desired quantity
with relative error at most $\epsilon^k$.
Despite this apparent power, our necessary conditions in Theorem~\ref{gen_result} 
and Section~\ref{vander} will still show limits on what can be evaluated accurately.
For example, no irreducible polynomial of degree $\ge 3$ can be
accurately evaluable over $\C^n$ if only dot products (degree $d=2$) are available.

\item Another standard technique for extending floating point precision
is to split a floating point number $x$ with $b$ bits in its fraction
into the exact sum $x = x_{hi} + x_{lo}$, where $x_{hi}$ and $x_{lo}$
each have only $b/2$ bits in their fractions. Then products like
\[
x \cdot y = (x_{hi}+x_{lo}) \cdot (y_{hi}+y_{lo}) = 
x_{hi} \cdot y_{hi} + x_{hi} \cdot y_{lo} +
x_{lo} \cdot y_{hi} + x_{lo} \cdot y_{lo}
\]
can be represented exactly as a sum of 4 floating point numbers,
since each product like $x_{hi} \cdot y_{hi}$ has at most $b$ bits in
its fraction and so can be computed without error. Arbitrary products may be
computed accurately by applying this technique repeatedly, which is the basis 
of some extra-precise software floating point libraries like \cite{shewchuk}.
The proof of accuracy of the algorithm for splitting $x = x_{hi} + x_{lo}$ 
is intrinsically discrete, and depends on a sequence of classical operations 
some of which can be proven to be free of error \cite[Thm 17]{shewchuk}, 
and similarly for the exactness of $x_{hi} \cdot y_{hi}$. 
Therefore this exact multiplication operation
cannot be built from simpler, inexact operations in our classical model.
But we may still model this approach as follows:
We imagine an exact multiplication operation $x \cdot y$, and note that
all we can do with it is feed it into the inputs of other operations.
This means that from the operation $rnd(z+w)$ we also get $rnd(x \cdot y + w)$,
$rnd(x \cdot y + r \cdot s)$, $rnd(x \cdot y \cdot z + w)$, and so on.
In other words, we take the other operations in our model and from
each create an enumerable set of other black-boxes, to which we can apply our
necessary and sufficient conditions. 

\end{enumerate}

\subsection{Comparisons and Branching}
\label{sec_ComparisonsBranching}
Are we permitted to do comparisons and then branch based
on their results? Are comparisons {\em exact}, i.e.,
are the computed values of $x > y$, $x=y$ and $x<y$ (true or false) always correct
for real $x$ and $y$?
(For complex $x$ and $y$ we consider only the comparison $x=y$.)
We consider algorithms both without comparisons (in which case
$p_{comp} (x,\delta)$ is simply a polynomial), 
and with exact comparisons and branching
(in which case $p_{comp}(x, \delta)$ is a piecewise polynomial,
on semialgebraic sets determined by inequalities among
other polynomials in $x$ and $\delta$).
We conjecture that using comparisons and
branching strictly enlarges the set of polynomials that we can 
evaluate accurately.

We note that by comparing $x-r_i$ to zero for selected constants
$r_i$, we could extract part of the bit representation of $x$. 
Since we are limiting ourselves to a finite number of operations, 
we could at most approximate $x$ this way, and as stated in
Section~\ref{sec_ExactRoundedInputs}, this means we could not exploit
this to get high relative accuracy near $p(x)=0$.
We note that the model of arithmetic in 
\cite{cuckersmale99} 
excludes real$\rightarrow$integer conversion
instructions.




\subsection{Nondeterminism}
\label{sec_Nondeterminism}

As currently described, our model is {\em nondeterministic},
e.g., the rounded result of $1+1$ is not necessarily identical 
if it is performed more than once. This is certainly
different behavior than the deterministic computers whose
behavior we are modeling. However, it turns out that this
is not a limitation, because we can always simulate a
deterministic machine with a nondeterministic one using
comparisons and branching. The idea is simple: The first
addition instruction (say) records its input arguments and 
computed sum in a list.
Every subsequent addition instruction compares its arguments
to the ones in the list (which it can do exactly), 
and either just uses the precomputed sum if it finds them,
or else does the addition and appends the results to the list.
In other words, the existence (or nonexistence) of an accurate
algorithm in a model with comparisons and branching does not
depend on whether the machine is deterministic. So for
simplicity, we will henceforth assume that our machines are
nondeterministic.


\subsection{Choice of Domain $\cal D$}
\label{sec_Domain}
As mentioned in the introduction, it seems natural to consider
any semialgebraic set as a possible domain of evaluation for
$p(x)$. While some choices, like ${\cal D} = \{ x: p(x) = 0\}$
make evaluating $p(x)$ trivial, they beg the question of how
one would know whether $x \in \cal D$. Similarly, if $\cal D$
includes a discrete set of points, then $p(x)$ can be evaluated
at these points by looking up the answers in a table. 
To avoid these pathologies, it may seem adequate restrict $\cal D$ to be
a sufficiently ``fat'' set,  say open. But this still leads to
interesting complications; for example the algorithm
\[
p_{comp}(x,\delta) = ((x_1+x_2)(1+ \delta_1)+x_3)(1+\delta_2)
\]
for $p(x) = x_1+x_2+x_3$, which is inaccurate on $\R^n$, 
is accurate on the open set
\linebreak
$\{ |x_1+x_2|>2|x_3| \}$, whose closure intersects
the variety $V_{\R}(p)$ on $\{x_1+x_2=0 \wedge x_3=0\}$.

In this paper we will mostly deal with open $\cal D$, 
especially ${\cal D} = \R^n$ or ${\cal D} = \C^n$,
and comment on when our results apply to smaller $\cal D$.

\subsection{Summary of Arithmetic and Algorithmic Models}
\label{Sec_Axioms}

We summarize the axioms our arithmetic and algorithms must satisfy.
We start with the axioms all arithmetic operations and algorithms
satisfy:

\begin{description}
\item[Exact Inputs.] Our algorithm will be given the input exactly.

\item[Finite Convergence.] An accurate algorithm must, when all roundoff
errors $\delta = 0$, compute the exact value of the polynomial in a finite
number of steps.
\item[Roundoff Model.] Except for negation, which is always exact,
the rounded value of any arithmetic operation $op(x_1,....,x_k)$ satisfies
\[
rnd(op(x_1,...,x_k)) = op(x_1,...,x_k)\cdot (1+ \delta)
\]
where
$\delta$ is arbitrary number satisfying $|\delta| \leq \epsilon$,
where $\epsilon$ is a nonzero value called the {\em machine precision}.
If the data $x$ is real (or complex), then $\delta$ is also real (resp. complex).
\item[Nondeterminism.] Every arithmetic operation produces an independent
roundoff error $\delta$, even if the arguments to different operations are
identical. 
\item[Domain $\cal D$.] Unless otherwise specified, the domain of evaluation
$\cal D$ is assumed to be all of $\R^n$ (or all of $\C^n$).
\end{description}

We now list the alternative axioms our algorithms may satisfy.
In each category, an algorithm must satisfy one set of axioms
or the other.

\begin{description}
\item[Branching or Not.] Some of our algorithms will permit exact comparisons of
intermediate quantities ($<$, $=$ and $>$ for real data), and subsequent branching 
based on the result of the comparison. Other algorithms will not permit branching.
In the complex case, we will see that branching does not matter (see Sections~\ref{suf_c}
and~\ref{suf_c_bb}).

\item[Classical or ``black-box'' operations.]
Some of our algorithms will use only ``classical'' arithmetic operations,
namely addition, subtraction and multiplication.
Others will use a set of arbitrary polynomial 
``black-box'' operations,
like $op(x,y,z) = x + y \cdot z$ or $op(x) = x - \sqrt{2}$, of our choice.
In particular, we omit division.

\end{description}








\subsection{Other Models of Error and Arithmetic}
\label{sec_OtherModels}

Our goal in this paper is to model rounded, finite precision computation,
i.e., arithmetic with numbers represented in scientific notation,
and rounded to their leading $k$ digits, for some fixed $k$.
It is natural to ask about models related to ours.

First, we point out some positive attributes of our model:

\begin{enumerate}

\item The model $rnd(op(a,b)) = op(a,b)(1+ \delta)$ has been the most widely
used model for floating point error analysis \cite{higham96}
since the early papers of 
von Neumann \cite{VonNeumannGoldstine47}, 
Turing \cite{turing} and 
Wilkinson \cite{wilkinsonroundingerror}.

\item The extension to include black-boxes includes widely used 
floating point techniques for extending the precision.

\item Though the model is for real (or complex) arithmetic, it can
be efficiently simulated on a conventional Turing machine by
using a simple variation of floating point numbers $m \cdot 2^e$,
stored as the pair of integers $(m,e)$,
where $m$ is of fixed length, and $|e|$ grows as necessary.  
In particular, any sequence of $n$ addition, subtraction, 
multiplication or division (by nonzero) operations can increase the 
largest exponent $e$ by at most $O(n)$ bits, and so can be done
in time polynomial in the input size. 
See \cite{demmelkoev99} for further discussion.
This is in contrast to 
repeated squaring in the BSS model \cite{Blum04} which
can lead to exponential time simulations.

\end{enumerate}

Models of arithmetic may be categorized according to 
several criteria (the references below are not exhaustive,
but illustrative):

\begin{itemize}
\item Are numbers (and any errors) represented discretely 
(e.g., as bit strings such as floating point numbers) 
\cite{demmelkoev99,higham96,wilkinsonroundingerror},
or as a (real or complex) continuum \cite{BCSS,cuckerdedieu01}? 

\item
Is arithmetic exact \cite{BCSS,BCSS96a}
or rounded
\cite{CuckerGrigoriev99,cuckersmale99,higham96,wilkinsonroundingerror}?
If it is rounded, is the error bounded in a 
relative sense \cite{higham96},
absolute sense \cite{BCSS}, 
or something else \cite{OlverLozier,demmel87d,demmel84}
\cite[Sec. 2.9]{higham96}?

\item
In which of these metrics is the final error assessed?

\item
Is the input data exact \cite{BCSS}
or considered ``rounded'' 
from its true value 
\cite{chatelinfraysse,EdelatSunderhauf98,KoKI,pourelrichards89,renegar94}
(and if rounded, again how is the error bounded)?

\item
Do we want a ``worst case'' error analysis 
\cite{higham96,wilkinsonroundingerror},
or by modeling rounding errors as random variables, a
statistical analysis \cite{vignes,KahanImprob,SpielmanTengICM02}
\cite[Sec. 2.8]{higham96}?
Does a condition number appear explicitly in the complexity
of the problem \cite{cuckersmale99}?

\end{itemize}

First we consider floating point arithmetic itself, i.e., where
real numbers are represented by a pair of integers $(m,n)$
representing the real number $m \cdot r^n$, 
where $r$ is a fixed number called the {\em radix} 
(typically $r=2$ or $r=10$).
Either by using one of many techniques in the literature
for using an array $(x_1,...,x_s)$ of floating point numbers 
to represent $x = \sum_{i=1}^s x_i$ to very high 
accuracy and to perform arithmetic on such high precision
numbers (e.g., \cite{bailey1a,bailey1b,priest}),
or by converting $m \cdot r^n$ to an exact rational
number and performing exact rational arithmetic, 
one can clearly evaluate {\em any} polynomial $p(x)$ without error,
and the only question is cost.
In light of this, our results on classical vs black-box arithmetic
can be interpreted as saying when  such high precision techniques are necessary,
and which black-box operations must be implemented this way,
in order to evaluate $p$ accurately.

Let us revisit the accurate evaluation of the simple polynomial $y_1 + y_2 + y_3$. 
The obvious algorithm is to carry enough digits so that the sum is computed exactly,
and then rounded at the end. But then to compute $(2^e+1)-2^e$ accurately
would require carrying at least $e$ bits, which is {\em exponential} in the size
of the input ($\log_2 e$ bits to represent $e$). Instead, most practical
algorithms rely on the technique in the above paragraph, repeatedly replacing
partial sums like $y_1 + y_2$ by $x_1 + x_2$ where $|x_1| \gg |x_2|$ and
in fact the bits of $x_1$ and $x_2$ do not ``overlap.'' These techniques
depend intrinsically on the discreteness of the number representation
to prove that certain intermediate additions and subtractions are in
fact exact. Our model treats this by modeling the entire operation
as a black-box (see Section~\ref{sec_BasicArithmeticOperations}).

Second, consider our goal of guaranteed high {\em relative} accuracy.
One might propose that {\em absolute} accuracy is a more
tractable goal, i.e., guaranteeing
$|p_{comp}(x,\delta) - p(x)| \leq \eta$ instead of
$|p_{comp}(x,\delta) - p(x)| \leq \eta |p(x)|$.
However, we claim that as long as our basic arithmetic
operations are defined to have bounded relative error
$\epsilon$, then trying to attain relative error
in $p_{comp}$ is the most natural goal.

Indeed, we claim that tiny absolute accuracy is
impossible to attain for {\em any} nonconstant
polynomial $p(x)$
when ${\cal D} = \R^n$ or ${\cal D} = \C^n$.
For example, consider $p(x) = x_1 + x_2$, for which
the obvious algorithm is
$p_{comp}(x,\delta) = (x_1+x_2)(1+ \delta)$.
Thus the absolute error 
$|p_{comp}(x,\delta) - p(x)| = |x_1+x_2|\delta \leq |x_1+x_2| \epsilon$.
This absolute error is at most $\eta$ precisely when
$|x_1+x_2| \leq \eta / \epsilon$, i.e., for $x$ in a diagonal strip in the 
$(x_1,x_2)$ plane. For $p(x) = x_1 \cdot x_2$ we analogously get
accuracy only for $x$ in a region bounded by hyperbolas.
In other words, even for the simplest possible polynomials that take
one operation to evaluate, they cannot be evaluated to high absolute
accuracy on most of ${\cal D} = \R^n$ or $\C^n$.
The natural error model to consider when trying to attain
low absolute error in $p(x)$ is to have low absolute error
in the basic arithmetic operations, and this is indeed the
approach taken in \cite{cuckersmale99}
(though as stated before, repeated squaring can lead to
an exponential growth in the number of bits a real number
represents \cite{Blum04}).

One could also consider more complicated error models,
for example {\em mixed absolute/relative error}:
$|p_{comp} (x, \delta) - p(x)| \leq \eta \cdot \max( |p(x)|, 1 )$.
Similar models have been used to model underflow error in floating 
point arithmetic~\cite{demmel84}.
A small mixed error implies that either the relative error or
the absolute error must be small, and so may be easier to attain than
either small absolute error or small relative error alone. 
But we argue that, at least
for the class of homogeneous polynomials evaluated on homogeneous
$\cal D$, the question of whether $p(x)$ is accurately evaluable
yields the same answer whether we mean
accuracy in the relative sense or mixed sense.
To see why, note that $x \in \cal D$ if and only if
$\alpha x \in \cal D$ for any scalar $\alpha$, since
$\cal D$ is homogeneous, and 
that $p(\alpha x) = \alpha^d p(x)$, where $d = degree(p)$.
Thus for any nonzero $p(x)$, scaling $x$ to $\alpha x$
will make $\eta \cdot \max( |p(\alpha x)|, 1 ) = \eta | p(\alpha x) |$
once $\alpha$ is large enough, i.e., relative error $\eta$ must be
attained. By results in Section~\ref{suf_r}, this will mean that
$p_{comp}(x, \delta)$ must also be homogeneous in $x$ of the
same degree, i.e., 
$p_{comp}(\alpha x, \delta) = \alpha^d p_{comp}(x,\delta)$.
Thus for any $x \in {\cal D}$ at which we can evaluate $p(x)$
with high mixed accuracy, we can choose $\alpha$ large enough so that
\[
\alpha^d |p_{comp}(x,\delta) - p(x)| = 
|p_{comp}(\alpha x,\delta) - p(\alpha x)| \leq
\eta \cdot \max( |p( \alpha x)|, 1 ) = 
\eta \cdot |p( \alpha x)| = 
\alpha^d \cdot \eta \cdot |p(x)| 
\]
implying that $p(\alpha x)$ can be evaluated with high
relative accuracy for all $\alpha$.
In summary, changing our goal from relative accuracy to mixed
relative/absolute accuracy will not change any of our results,
for the case of homogeneous $p$ and homogeneous $\cal D$.

Yet another model is to assume that the input $x$ is given
only approximately, instead of exactly as we assume.
This corresponds to the approach taken in 
\cite{EdelatSunderhauf98,KoKI,pourelrichards89}, 
in which one can imagine reading as many leading bits as 
desired of each input $x_i$ from an infinite tape,
after which one tries to compute the answer using a conventional
Turing machine model. This gives yet different results, since,
for example, the difference $x_1 - x_2$ cannot be computed
with small relative error in a bounded amount of time, since
$x_1$ and $x_2$ may agree in arbitrarily many leading digits.
Absolute error is more appropriate for this model.

It is worth commenting on why high accuracy of the sort we want
is desirable in light of inevitable uncertainties in the inputs.
Indeed, many numerical algorithms are successfully analyzed using
{\em backward error analysis} \cite{higham96,demmelMA221},
where the computed results are shown to be the exact result for
a slightly perturbed value of the input. This is the case, for
example, for polynomial evaluation using Horner's rule
where one shows that one gets the exact value of a polynomial at
$x$ but with slightly perturbed coefficients.
Why is this not always accurate enough?

We already mentioned mesh generation \cite{shewchuk}, where the
inputs are approximately known physical coordinates of some physical
object to be triangulated, but where geometric predicates about the
vertices defining the triangulation must be answered consistently;
this means the signs of certain polynomials must be computed exactly,
which is in turn guaranteed by guaranteeing any relative accuracy
$\eta < 1$.

More generally, in many physical simulations, the parameters describing
the physical system to be simulated are often known to only a few 
digits, if that many. Nonetheless, intermediate computations must
be performed to much higher accuracy than the input data is known,
for example to make sure the computed system conserves energy
(which it should to high accuracy for the results to be meaningful,
even if the initial conditions are uncertain). 

Another example where high accuracy is important
are the trigonometric functions: 
When $x$ is very large and slightly 
uncertain, the value of $\sin x$ may be completely uncertain. 
Still, we want the computed trigonometric functions to
(nearly) satisfy identities like $\sin^2 x + \cos^2 x = 1$
and $\sin 2x = 2 \sin x \cos x$ so that we can reason about
program correctness. Many other examples of this sort can be
found in articles posted at \cite{KahanWebPage}.

In the spirit of backward error analysis, one could consider
the polynomial $p$ fixed, but settle for accurately computing
$p(\hat{x})$ where $\hat{x}$ differs from $x$ by only a 
small relative change in each component $x_i$. This is not
guaranteed by Horner's rule, which is equivalent to changing 
the polynomial $p$ slightly but not $x$. Would it be easier to 
compute $p(\hat{x})$ accurately than $p(x)$ itself?
This is the case for some polynomials, like $x_1 + x_2 + x_3$
or $c_1 x_2^2  x_3^3 + c_2 x_1^2 x_3^3 + c_3 x_1 x_2^4$, where there is a 
unique $x_i$ that we can associate with each monomial to 
``absorb'' the rounding error from Horner's rule. 
In particular, with Horner's rule, 
the number of monomials in $p(x)$ may at most be equal to the
number of $x_i$. In analogy to this paper, one could ask for 
a decision procedure to identify polynomials that permit
accurate evaluation of $p(\hat{x})$ using any algorithm.
This is a possible topic for future work.

Another possibility is to consider error probabilistically
\cite[Sec. 2.8]{higham96}.
This has been implemented in a practical system \cite{vignes},
where a program is automatically executed repeatedly with 
slightly different rounding errors made at each step in order 
to assess the distribution of the final error. This approach
is criticized in \cite{KahanImprob} for improperly modeling
the discrete, non-random behavior of roundoff, and for 
possibly invalidating (near) identities like 
$\sin 2x = 2 \sin x \cos x$ upon which correctness may depend.

In {\em smoothed analysis} \cite{SpielmanTengICM02},
one considers complexity (or for us, relative error) 
by averaging over a Gaussian distribution around each input. 
For us, input could mean either the argument $x$ of a fixed polynomial 
$p$, or the polynomial itself, or both.
First consider the case of a fixed polynomial $p$
with a randomly perturbed $x$.
This case 
is analogous to the previous paragraph, because the inputs
can be thought of as slightly perturbed before starting the
algorithm. Indeed, one could imagine rounding the inputs
slightly to nearby rational or floating point numbers, and
then computing exactly. But in this case, it is easy to
see that, at least for codimension 1 varieties of $p$,
the ``smoothed'' relative error is finite or infinite
precisely when the worst case relative error is finite or
infinite. So smoothing does not change our basic analysis.
\footnote{The logarithm of the relative error, like the logarithm
of many condition numbers, does however have a finite average.}
Now suppose one smooths over the polynomial $p$, i.e., over its
coefficients. If we smooth using a Gaussian distribution,
then as we will see, 
the genericity of 
``bad'' $p$ will make the smoothed relative error infinite
for all polynomials. Changing the distribution from Gaussian
to one with a small support would only distinguish between
positive definite polynomials, the easy case discussed
in section~\ref{sec_PositivePolys}, and polynomials that
are not positive definite.

In {\em interval arithmetic} \cite{moore,neumaier-book,alefeldherzberger}
one represents each number by a floating point interval guaranteed to
contain it. To do this one rounds interval endpoints ``outward'' to 
ensure that, for example, the sum $c=a+b$ of two intervals yields an 
interval $c$ guaranteed to contain the sum of any two numbers in $a$ 
and $b$.
It is intuitive that if an interval algorithm existed to evaluate
$p(x)$ for $x \in {\cal D}$ that always computed an interval whose 
width was small compared
to the number of smallest magnitude in the interval,
and if the algorithm obeyed the rules in Section~\ref{Sec_Axioms}, 
then it would satisfy our accuracy requirements.
Conversely, one might conjecture that an algorithm accurate by our
criteria would straightforwardly provide an accurate interval algorithm,
where one would simply replace all arithmetic operation by 
interval operations. The issue of interpreting comparisons and
branches using possibly overlapping intervals makes this question
interesting, and a possible subject for future work.

Finally, many authors use condition numbers in their analysis
of the complexity of solving certain problems. This is classical
in numerical analysis \cite{higham96}; 
more recent references are 
\cite{chatelinfraysse,cuckersmale99,cuckerdedieu01}. 
In this approach, one is willing to do more and more work to get 
an adequate answer as the condition number grows, perhaps without bound.
Such a conditioning question appears in our approach, if we ask
how small the machine precision $\epsilon$ must be as a function of
the desired relative error $\eta$, as well as $p$, $\cal D$, and allowed
operations. Computing this condition number
(outside the easy case described in Section~\ref{sec_PositivePolys}) 
is an open question.


\section{Evaluating positive polynomials accurately}
\label{sec_PositivePolys}

Here we address the simpler case where the polynomial $p(x)$ 
to be evaluated has no zeros in the domain of 
evaluation $\cal D$. It turns out that we need
more than this to guarantee accurate evaluability: we
will require that $|p(x)|$ be bounded both above and below 
in an appropriate manner on $\cal D$.

We let $\bar{\cal D}$ denote the closure of $\cal D$.

\begin{theorem}
\label{thm_positive_compact}
Let $p_{comp} (x,\delta)$ be {\em any} algorithm for
$p(x)$ satisfying $p_{comp}(x,0) = p(x)$, 
i.e. it computes the right value in the absence of rounding error.
Let $p_{min} \eqbd \inf_{x \in \bar{\cal D}} |p(x)|$.
Suppose $\bar{\cal D}$ is compact and $p_{min} > 0$.
Then $p_{comp}(x,\delta)$ is an accurate algorithm 
for $p(x)$ on $\cal D$.
\end{theorem}

\begin{proof}
Since the relative error on $\cal D$ is 
$|p_{comp}(x,\delta) - p(x)|/|p(x)| \leq |p_{comp}(x,\delta) - p(x)|/p_{min}$,
it suffices to show that the numerator approproaches 0 uniformly
as $\delta \rightarrow 0$. This follows by writing the value of
$p_{comp}(x,\delta)$ along any branch of the algorithm as
$p_{comp}(x,\delta) = p(x) + \sum_{\alpha > 0} p_{\alpha}(x) \delta^{\alpha}$,
where $\alpha > 0$ is a multiindex with at least one component exceeding 0.
By compactness of $\bar{\cal D}$,
$|\sum_{\alpha > 0} p_{\alpha}(x) \delta^{\alpha} | \leq C 
\sum_{\alpha > 0} |\delta|^{\alpha}$ for some constant $C$,
which goes to 0 uniformly as the upper bound $\epsilon$ on
each $|\delta_i|$ goes to zero.
\end{proof}

Next we consider domains $\cal D$ whose closure is not compact.
To see that merely requiring $p_{min}>0$ is not enough,
consider evaluating $p(x) = 1 + (x_1 + x_2 + x_3)^2$ on
$\R^3$. Intuitively, $p(x)$ can only be accurate if
its ``dominant term'' $(x_1+x_2+x_3)^2$ is accurate, 
once it is large enough, and this is not possible
using only addition, subtraction and multiplication.
(These observations will be formalized in 
Sections~\ref{suf_r} and \ref{class}, respectively.)

Instead, we consider a homogeneous polynomial $p(x)$ 
evaluated on a homogeneous $\cal D$, i.e. one where
$x \in \cal D$ implies $\gamma x \in \cal D$ for any scalar
$\gamma$.  Even though such $\cal D$ are
unbounded, homogeneity of $p$ will let us consider just
the behavior of $p(x)$ on $\cal D$ intersected with the
unit ball $S^{n-1}$ in $\R^n$ (or $S^{2n-1}$ in $\C^n$). 
On this intersection
we can use the same compactness argument as above:

\begin{theorem}
\label{thm_positive_homo}
Let $p(x)$ be a homogeneous polynomial, let $\cal D$
be a homogeneous domain, and let $S$ denote the unit
ball in $\R^n$ (or $\C^n$). 
Let 
\[
p_{min,homo} \eqbd \inf_{x \in \bar{\cal D} \cap S} |p(x)| \;
\]
Then $p(x)$ can be evaluated accurately if $p_{min,homo} > 0$.
\end{theorem}

\begin{proof}
We describe an algorithm $p_{comp} (x, \delta )$ for 
evaluating $p(x)$. There are many such algorithms,
but we only describe a simple one.
(Indeed, we will see that the set of all accurate algorithms 
for this situation
can be characterized completely by Definition~\ref{homalg} 
and Lemma~\ref{1}.)
Write $p(x) = \sum_{\alpha} c_\alpha x^\alpha$, where
$\alpha$ is a multiindex $(\alpha_1,...,\alpha_n)$,
$x^{\alpha} \eqbd x_1^{\alpha_1} \cdots x_n^{\alpha_n}$,
and $c_{\alpha} \neq 0$ is a scalar. 
Homogeneity implies $|\alpha| = \sum_i \alpha_i$ is constant.
Then the algorithm simply
\begin{enumerate}
\item computes each $x^{\alpha}$ term by repeated 
multiplication by $x_i$s,
\item computes each $c_{\alpha} x^{\alpha}$ either 
by multiplication by $c_{\alpha}$ or by repeated 
addition if $c_{\alpha}$ is an integer, and
\item sums the $c_{\alpha} x^{\alpha}$ terms.
\end{enumerate}
Since each multiplication, addition and subtraction
contributes a $(1 + \delta_i)$ term, it is easy to see that 
\[
p_{comp}(x,\delta) = \sum_{\alpha} c_{\alpha} x^{\alpha} 
\Delta_{\alpha}
\]
where each $\Delta_{\alpha}$ is the product of at most 
some number $f$ of factors of the form $1+ \delta_i$.

Now let $\|x\|_2 = (\sum_i |x_i|^2)^{1/2}$,
so $\hat{x} = x/\|x\|_2$ is in the unit ball $S$. 
Then the relative error may be bounded by
\begin{eqnarray*}
\left| \frac{p_{comp}(x,\delta) - p(x)}{p(x)} \right|
& = &
\left| \frac{\sum_{\alpha} c_{\alpha} x^{\alpha} \Delta_{\alpha}
- \sum_{\alpha} c_{\alpha} x^{\alpha}}
{\sum_{\alpha} c_{\alpha} x^{\alpha}} \right| \\
& = &
\left| \frac{\sum_{\alpha} c_{\alpha} \hat{x}^{\alpha} 
(\Delta_{\alpha} - 1)} 
{\sum_{\alpha} c_{\alpha} \hat{x}^{\alpha}} \right| \\
& \leq &
\frac{\sum_{\alpha} | c_{\alpha} | \cdot  |\Delta_{\alpha} - 1|}
{p_{min}} \\
& \leq &
\frac{\sum_{\alpha} | c_{\alpha} | \cdot  ((1+\epsilon)^f-1)}
{p_{min}}
\end{eqnarray*}
which goes to zero uniformly in $\epsilon$.
\end{proof}

\section{Classical arithmetic} \label{class}

In this section we consider the simple or classical arithmetic over the real or complex fields, with the three basic operations $\{+, -, \cdot\}$, to which we add negation. The model of arithmetic is governed by the laws in Section~\ref{Sec_Axioms}. We remind the reader that this arithmetic model \emph{does not allow} the use of constants.

In Section~\ref{class_nec} we find a necesary condition for accurate evaluability over either field, and in Section~\ref{suf_c} we prove that this condition is also sufficient for the complex case.

Throughout this section, we will make use of the following definition of allowability.

\begin{definition} \label{allwa} Let $p$ be a polynomial over $\mathbb{R}^n$ or $\mathbb{C}^n$, with variety $V(p)\eqbd \{ x~:~p(x)=0\}$. 
We call $V(p)$ \emph{allowable} if it can be represented as a union of intersections of sets of the form
\begin{eqnarray} \label{unu}
& 1.&  Z_i =\{x~:~ x_i ~=~ 0\}~, \\
\label{doi}
& 2.&  S_{ij} = \{x~:~ x_i+x_j ~=~ 0\}~, \\
\label{trei}
& 3.&  D_{ij} = \{x~:~ x_i - x_j ~=~ 0\}~.
\end{eqnarray}
If $V(p)$ is not allowable, we call it unallowable.
\end{definition}

\begin{remark} For a polynomial $p$, having an allowable variety $V(p)$ is obviously a Tarski-decidable property (following \cite{tarski_book}), since the number of unions of intersections of hyperplanes \eqref{unu}-\eqref{trei} is finite. \end{remark}

\subsection{Necessity: real and complex} \label{class_nec}

All the statements and proofs in this section work equally well for both the real and the complex case, and thus we may treat them together. At the end of the section we use the necessity condition to obtain a partial result relating to domains.

\begin{definition} From now we will refer to the space of variables as $\S \in \{\R^n, \C^n\}$.
\end{definition}

To state and prove the main result of this section, we need to introduce some additional notions and notation.

\begin{definition} \label{general_p} Given a polynomial $p$ over $\S$ with unallowable variety $V(p)$, consider 
all sets $W$ that are finite intersections of allowable hyperplanes defined by \eqref{unu}, 
\eqref{doi}, \eqref{trei}, and subtract from $V(p)$ those $W$ for which $W \subset V(p)$. 
We call the remaining subset of the variety {\em points in general position\/} and denote it by $G(p)$.
\end{definition}

\begin{remark} If $V(p)$ is not allowable, then from definition~\ref{general_p} it follows that $G(p)\neq \emptyset$. One may also think of points in $G(p)$ as ``unallowable'' or ``problematic'', because, as we will see, we necessarily get large relative errors in their vicinity.
\end{remark}

\begin{definition} \label{allowance} 
Given $x \in \S$, define the set $\allow(x)$ as the intersection of
all allowable hyperplanes going through $x$:
$$ \allow(x)\eqbd \left( \cap_{x\in Z_i} Z_i\right)  \cap \left( \cap_{x\in S_{ij}} S_{ij}\right) 
 \cap \left( \cap_{x\in D_{ij}} D_{ij} \right) , $$
with the understanding that 
$$ \allow(x)\eqbd \mathcal{S} \qquad {\rm whenever} \qquad x\notin Z_i, \; S_{ij}, \; D_{ij}
\quad \hbox{\rm for all} \quad i,j.$$
Note that $\allow(x)$ is a linear subspace of $S$.
\end{definition}

We will be interested in the sets $\allow(x)$ primarily when $x\in G(p)$. 
For such cases we make the following observation.

\begin{remark} For each $x\in G(p)$, the set $\allow(x)$ is not a subset of $V(p)$:
$$ \allow(x)\not\subseteq V(p),$$ 
which follows directly from the definition of $G(p)$.
\end{remark}  

We can now state the main result of this section, which is a necessity condition for the evaluability of polynomials over domains. 

\begin{theorem} \label{conj}
Let $p$ be a polynomial over a domain $\mathcal{D} \in \S$. Let $G(p)$ be the set of points in general position on the variety $V(p)$. If there exists $x \in \mathcal{D} \cap G(p)$ such that $\allow(x) \cap \Int(\mathcal{D}) \neq \emptyset$, then $p$ is not accurately evaluable on $\mathcal{D}$.
\end{theorem}

To prove Theorem~\ref{conj}, we need to recall the notion of Zariski topology (see, e.g.,~\cite{KH}). 

\begin{definition} \label{zariski} A subset $Y\subseteq \R^n$ (or $\mathbb{C}^n$) is called a {\sl Zariski closed\/} 
set if there a subset $T$ of the polynomial 
ring $\R[x_1,\ldots,x_n]$ (or $\C[x_1, \ldots, x_n]$) such that 
$Y$ is the variety of $T$: $Y=V(T):=\cap_{p \in T} V(p)$. A complement of a Zariski 
closed set is said to be {\sl Zariski open\/}. The class of Zariski open
sets defines the {\sl Zariski topology\/} on $\S$.
\end{definition}

In this paper, we consider the Zariski topology not on $\S$, but on a hypercube centered at the origin in $\delta$-space (the space in which the vector of error variables $\delta$ lies). This topology is defined in exactly the same fashion. 

Note that a Zariski closed set has measure zero unless it is defined by 
the zero polynomial only; then the set is the whole space. In the coming proof 
we will deal with nonempty Zariski open sets, which are all of full measure. Finally, it is worth noting that the Zariski sets we will work with are algorithm-dependent. 

Finally, we represent any algorithm as in~\cite{smale2, AHU} by a directed acyclic graph (DAG)
with input nodes, branching nodes, and output nodes. 
For simplicity in dealing with negation (given that negation is an 
\emph{exact\/} operation), we define a special type of edge which indicates 
that the value carried along the edge is negated. We call these special edges 
\emph{dotted,\/} to distinguish them from the regular \emph{solid} ones. 

 Every computational node has two inputs (which may both come from a
single other computational node); depending on the source of these inputs
we have computational nodes with inputs from two distinct nodes and
computational nodes with inputs from the same node. The latter type
correspond either to \begin{enumerate} \item doubling ($(x, x)
\stackrel{+}{\mapsto} 2x$), \item doubling and negating ($(-x, -x)
\stackrel{+}{\mapsto} -2x$), \item computing zero exactly ($(-x,
x)\stackrel{+}{\mapsto} 0$, $(-x, -x)\stackrel{-}{\mapsto} 0 $, or $(x,
x)\stackrel{-}{\mapsto} 0$), \item squaring ($(x, x)
\stackrel{\cdot}{\mapsto} x^2$ or $(-x, -x) \stackrel{\cdot}{\mapsto}
x^2$), \item squaring and negating ($(-x, x) \stackrel{\cdot}{\mapsto}
-x^2$). \end{enumerate}

All nodes are labeled by $(op(\cdot), \delta_i)$ with $op(\cdot)$ representing the operation 
that takes place at that node. It means that at each node, the algorithm takes in two inputs, executes 
the operation, and multiplies the result by $(1+\delta_i)$. 

Finally, for each branch, there is a single destination node, with one input and no output, whose input value is the result of the algorithm.

{Throughout the rest of this section, unless specified, we consider only non-branching algorithms.} 

\begin{definition} \label{non-trivial} For a given $x\in \S$, we say that a computational node $N$ is {\em of non-trivial type\/} if its output is a nonzero polynomial in the variables $\delta$ when 
the algorithm is run on the given $x$ and with symbolic $\delta$s. \end{definition}




\begin{definition} \label{left_and_right}
For a fixed $x$, let $N$ be any non-trivial computational node in an algorithm. 
We denote by $L(N)$ (resp., $R(N)$) the set of computational nodes in the left 
(resp., right) subgraphs of $N$. If both inputs come from the same node, i.e. $L(N)$ and $R(N)$ overlap, we will only talk about $L(N)$.
\end{definition}

\begin{definition} \label{eps-hypercube} For a given $\epsilon>0$, we denote by $H_{\epsilon}$ 
the hypercube of edge length $2\epsilon$ centered at the origin, in $\delta$-space. 
\end{definition}

We will need the following Proposition.

\begin{proposition} \label{unu_unu} Given any algorithm, any $\epsilon>0$,  
and a point 
$x\in G(p)$, there exists a  Zariski open set $\Delta$ in 
$H_{\epsilon}$ such that no non-trivial computational node has a zero 
output on the input $x$ for all $\delta \in \Delta$.
\end{proposition}

\begin{proof}
The proof follows from the definition of the non-trivial computational node. 

Since every non-trivial computational node outputs a non-trivial polynomial in $\delta$, it follows that each non-trivial computational node is nonzero on a Zariski open set (corresponding to the output polynomial in $\delta$) in $H_{\epsilon}$. Intersecting this finite number of Zariski open sets we obtain a Zariski open set which we denote by $\Delta$; for any $\delta \in \Delta$ the output of any non-trivial computational node is nonzero.
\end{proof}

We can now state and prove the following crucial lemma.

\begin{lemma} \label{prima}
For a given algorithm, any $x\in G(p)$, and $\epsilon>0$, exactly one of the following holds:
\begin{enumerate} \item there exists a  Zariski open set $\Delta \subseteq H_{\epsilon}$ such that the value $p_{comp}(x)$ computed by the algorithm is not zero when the algorithm 
is run with source input $x$ and $\delta \in \Delta$;
\item $p_{comp}(y, \delta)=0$ for all $y\in \allow(x)$ and all $\delta$ in $H_{\epsilon}$.
\end{enumerate}
\end{lemma}

\noindent \textit{Proof of Lemma~\ref{prima}.} We recall that the algorithm can be represented as a DAG, as described in the paragraphs preceding Definition~\ref{non-trivial}. 

Fix a point $x\in G(p)$. Once $x$ is fixed, the result of each computation is a polynomial expression in the $\delta$s. Consider the Zariski open set $\Delta$ whose existence is guaranteed by Proposition~\ref{unu_unu}. There are now two possibilities: either the output node is of non-trivial type, in which case $p_{comp}(x, \delta) \neq 0$ for all $\delta \in \Delta$, or the output node is not of non-trivial type, in which case $p_{comp}(x, \delta_0) = 0$ for some $\delta_0 \in \Delta$.

 In the latter case the output of the computation is zero; we trace back this zero to its origin, by marking in descending order all computational nodes that produced a zero (and thus we get a set of paths in the DAG, all of whose nodes produced exact zeros). Note that we are not interested in \emph{all\/} nodes that produced a $0$; only those which are on paths of zeros to the output node. 

We will examine the last occurrences of zeros on paths of marked vertices, i.e. the zeros that are farthest from the output on such paths. 

\begin{lemma} \label{enumer} 
The \emph{last\/} zero on such a path must be either 
\begin{enumerate}
\item a source;
\item the output of a node where $(-x, x) \stackrel{+}{\mapsto} 0$, $(-x, -x) \stackrel{-}{\mapsto} 0$, or $(x, x) \stackrel{-}{\mapsto} 0$ are performed;
\item the output of an addition or subtraction node with two nonzero source inputs.
\end{enumerate}
\end{lemma}

\noindent \textit{Proof of Lemma~\ref{enumer}.}  Note that a nonzero non-source output will be a non-constant polynomial in the $\delta$ specific to that node.

Clearly the last zero output cannot happen at a multiplication node; we have thus to show that the last occurrence of a zero output cannot happen at an addition or subtraction node which has two nonzero inputs from different nodes, at least one of which is a non-source. We prove the last statement by reductio ad absurdum. 

Assume we could have a zero output at a node $N$ with two nonzero inputs, at least one of which is not a source. Let $R(N)$ and $L(N)$ be as in Definition~\ref{left_and_right}. Let $\delta(L(N))$ and $\delta(R(N))$ be the sets of errors $\delta_i$ corresponding to the left, respectively the right subtrees of $N$. 

By assumption, $\delta(R(N)) \cup \delta(L(N)) \neq \emptyset$ (since at least one of the two input nodes is a non-source). Let $\delta_l$ ($\delta_r$) denote the $\delta$ associated to the 
left (right) input node of $N$. Then we claim that either $\delta_l \notin \delta(R(N))$ or $\delta_r\notin \delta(L(N))$.  (There is also the possibility that one of the two input nodes \emph{is\/} a source and does not have a $\delta$, but in that case the argument in the next paragraph becomes trivial.)

Indeed, since each $\delta$ is specific to a node, if $\delta_l$ were in $\delta(R(N))$, there would be a path from the left input node to the right input node. Similarly, if  $\delta_r$ were in $\delta(L(N))$, then there would be 
a path from the right input node of $N$ to the left input node of $N$. So if both events were to happen at the same time, there would be a cycle in the DAG. This cannot happen, hence either $\delta_l \notin \delta(R(N))$ or $\delta_r\notin \delta(L(N))$.

Assume w.l.o.g. $\delta_l \notin \delta(R(N))$. Then the left input of $N$ is a non-trivial polynomial in $\delta_l$, 
while the right input does not depend on $\delta_l$ at all. Hence their sum or difference is still a non-trivial 
polynomials in $\delta_l$. Contradiction. \qed 

Now that Lemma~\ref{enumer} has been proven, we can state the crucial fact of the proof of Lemma~\ref{prima}: \emph{all last occurrences of a zero appear at nodes which either correspond to allowable constraints (i.e., zero sources, or sums and differences of sources), or are addition/subtraction nodes with both inputs from the same node, which always, on any source inputs, produce a zero}. 

Take now any point $y\in \allow(x)$; then $y$ produces the same chains of consecutive zeros constructed (marked) in Lemma~\ref{prima} as $x$ does, with errors given by $\delta_0 \in \Delta$. Indeed, any node on such a chain that has a zero output at $x$ when the error variables are $\delta_0$ can trace this zero back to an allowable constraint (which is satisfied by both $x$ and $y$) or to an addition/subtraction node with both inputs from the same node; hence the node will also have a zero output at $y$ with errors $\delta_0$. In particular, if $p_{comp}(x, \delta_0)=0$ for $\delta_0\in \Delta$, then $p_{comp}(y, \delta_0)=0$. Moreover, changing $\delta_0$ can only introduce additional zeros, but cannot 
eliminate zeros on the zero paths that we traced for $x$ (by the choice of $\Delta$). Therefore, $p_{comp}(y, \delta)=0$ for all $y\in \allow(x)$
and $\delta \in H_{\epsilon}$. This completes the proof of Lemma~\ref{prima}. \qed

From Lemma~\ref{prima} we obtain the following corollary.

\begin{corollary} \label{partsial} For any $\epsilon>0$ and any $x\in G(p)$, exactly one of the following holds: the relative error of computation, 
$|p_{comp}-p|/|p|$, is either infinity at $x$ for all $\delta$ in a  Zariski open set or $1$ at all
points $y\in (\allow(x)\setminus V(p))$ and all $\delta \in H_{\epsilon}$. 
\end{corollary}

We now consider algorithms with or without branches.

\begin{theorem} \label{finala} 
Given a (branching or non-branching) algorithm with output function $p_{comp}(\cdot)$, $x\in G(p)$, and $\epsilon>0$, then one of the following is true:
\begin{enumerate}
\item there exists a set $\Delta_1$ of positive measure in $H_{\epsilon}$  such that $p_{comp}(x, \delta)$ is nonzero whenever the algorithm is run with errors $\delta \in \Delta_1$, or
\item there exists a set $\Delta_2$ of positive measure in $H_{\epsilon}$ such that for every $\delta \in \Delta_2$, there exists a neighborhood $N_{\delta}(x)$ of $x$ such that 
for every $y \in N_{\delta}(x) \cap \left ( \allow(x) \setminus V(p) \right )$, $p_{comp}(y, \delta) = 0$ when the algorithm is run with errors $\delta$.
\end{enumerate}
\end{theorem}

\begin{remark}
This implies that, on a set of positive measure in $H_{\epsilon}$, the relative accuracy 
of any given algorithm is either $\infty$ or $1$.  
\end{remark}

\begin{proof}
With $p_{comp}(\cdot)$ the output function and $x$ a fixed point in general position, we keep the $\delta$s symbolic. 
Depending on the results of the comparisons, the algorithm splits into a finite number of non-branching algorithms,
which all start in the same way (with the input nodes) and then differ in accordance with a finite set of polynomial 
constraints on the $\delta$s and $x$s.  

Some of these branches will be chosen by sets of $\delta$s of measure zero; at least one of the branches will have 
to be chosen by a set of $\delta$s of positive measure whose interior is nonempty (all constraints being polynomials).
Call that branch $B$, and let the set of $\delta$s that choose it be called $\Delta_B$. 

By Proposition~\ref{unu_unu}, there exists a  Zariski open set $\Delta \in H_{\epsilon}$ such that, for all $\delta \in \Delta$, no non-trivial node in the subgraph representing our branch $B$ 
has a zero output. In particular, this includes all quantities computed for comparisons that define $B$. 
Let $\Delta_2\eqbd \Int(\Delta_B \cap \Delta)$, where
$\Int$ denotes the interior of a set. By the choice of $\Delta_B$ and $\Delta$, the obtained set $\Delta_2$
is non-empty.  

Suppose the algorithm is run with errors $\delta_0\in \Delta_2$ and $p_{comp}(x, \delta_0)\neq 0$. Then, by continuity, there 
must be a neighborhood $\Delta_1$ in the set $\Delta_2$ on which the computation will still be directed to 
branch $B$ and $p_{comp}(x, \cdot)$ will still be nonzero, so we are in Case 1. 

Assume now that we are not in Case 1, i.e. there is no $\delta \in \Delta_2$ such that $p_{comp}(x, \delta) \neq 0$. In this case we show by contradiction that $p_{comp}(y, \delta)=0$ for all 
$y\in \allow(x)$ if $y$ is sufficiently close to $x$ (since $\allow(x)$ is a linear subspace containing $x$, there exist points in $\allow(x)$ which are arbitrarily close to $x$), thus, that Case 2 must be fulfilled.

If this claim is not true, then there is no neighborhood 
$N_{\delta}(x)$ of $x$ such that when $y \in N_{\delta}(x) \cap \allow(x)$, the algorithm is directed to branch $B$ 
on $\delta$. In that case, there must be a sequence $\{y_n\} \in \left ( \allow(x) \setminus V(p) \right)$ such that $y_n \rightarrow x$ and 
$y_n$ is always directed elsewhere for this choice of $\delta$. The reason for this is that $\allow(x)$ is a linear subspace which is \emph{not\/} contained in $V(p)$; hence no neighborhood of $x$ in $\allow(x)$ can be contained in $V(p)$, and then such a sequence $y_n$ must exist. 

Since there is a finite number of branches, we might as well assume that all $y_n$ will be directed to the same branch $B'$ for this $\delta$ and that they
split off at the same branching node (pigeonhole  principle). 

Now consider the branching node where the splitting occurs, and let $r(z, \delta)$ be the quantity to be compared to $0$ at that node. 
Since we always go to $B'$ with $y_n$ but to $B$ with $x$, it follows that we \emph{necessarily must have\/} $r(y_n,\delta)\neq 0$ whereas
$r(x,\delta)=0$. On the other hand, until that splitting point the algorithm followed the same path with
$y_n$ and with $x$, computing with the same errors $\delta$. Applying then case 2 of Lemma \ref{prima} (which can be read to state that any algorithm computing $r$, and obtaining $r(x, \delta) = 0$, will also obtain $r(y_n, \delta) = 0$), we get a contradiction.

This completes the proof of Theorem~\ref{finala}. \end{proof}

\begin{corollary} \label{in_sfirsit} Let $p$ be a polynomial over $\S$ with unallowable variety $V(p)$. Choose any algorithm with output function $p_{comp}(\cdot)$, any point $x \in G(p)$, $\epsilon>0$, and $\eta<1$. Then there exists a set $\Delta_x$ of positive measure {\sl arbitrarily close} to $x$ and a set $\Delta$ of positive measure in $H_{\epsilon}$, such that $|p_{comp} - p|/|p|$ is strictly larger than $\eta$ when computed at a point $y \in \Delta_x$ using any vector of relative errors $\delta \in \Delta$.
\end{corollary}

\begin{proof}
On symbolic input $x$ and with symbolic $\delta$, the algorithm will have $m$ branches $B_1, \ldots, B_m$ that correspond to constraints yielding (semi-algebraic) sets of positive measure $S_1, \ldots, S_m$ in $(x, \delta)$-space. 
Choose $x \in G(p)$, and let $(x, 0)$ be a point in $(x, \delta)$-space. 
\begin{enumerate} 
\item \label{primo} If $(x, 0)$ is in $\Int(S_i)$ (the interior of some region $S_i$), then by Lemma \ref{prima} and Corollary~\ref{partsial} there exists either 
\begin{enumerate} \item a $\delta_0$ in $\delta$-space sufficiently small such that $(x, \delta_0)$ is in $\Int(S_i)$ \emph{and\/} $p_{comp}(x, \delta_0) \neq 0$. The relative error at $(x, \delta_0)$ is in this case $\infty$, and (by continuity) there must be a small ball around $(x, \delta_0)$ which is still in $\Int(S_i)$,  on which the minimum relative error is arbitrarily large, certainly  larger than $1$;
\item a $\delta_0$ in $\delta$-space sufficiently small and a $y \in \allow(x) \setminus V(p)$ sufficiently close to $x$ such that $(y, \delta_0)$ is in $\Int(S_i)$ and $p_{comp}(y, \delta_0) = 0$. In this case the relative error at $(y, \delta_0)$ is $1$, and (by continuity) there must be a small ball around $(y, \delta_0)$ which is still in $\Int(S_i)$, on which the relative error is strictly larger than our $\eta<1$.
\end{enumerate} 
\item Otherwise, $(x, 0)$ must be on the boundary of some of the regions $S_i$; assume w.l.o.g. that it is on the boundary of the regions $S_1, \ldots, S_l$. In this case, we choose a small hyperdisk $B_{\tilde{\epsilon}}((x, 0))$ in the linear subspace $(x, \cdot)$ such that $B_{\tilde{\epsilon}}((x, 0))$ intersects the closures of $S_1, \ldots, S_l$ (and no other $S_i$s). We can do this because the sets $S_i$ are all semi-algebraic.

\begin{enumerate} \item If there exists a $\delta_0$ in $\delta$-space such that $(x, \delta_0) \in B_{\tilde{\epsilon}}((x, 0))$ and $(x, \delta_0) \in \Int(S_i)$ for some $i \in \{1, \ldots, l\}$, then by the same argument as in case~\ref{primo} we obtain a small ball included in $\Int(S_i)$ on which the relative error is greater than $\eta$;
\item Otherwise, if there exists a $\delta_0$ such that $(x, \delta_0) \in B_{\tilde{\epsilon}}((x, 0))$ is on the boundary of some region $S_i$ for which the local algorithm corresponding to it would yield $p_{comp}(x, \delta_0) \neq 0$, then (by continuity) there exists a small ball around $(x, \delta_0)$ such that the intersection of that small ball with $S_i$ is of positive measure, and the relative error on that small ball as computed by the algorithm corresponding to $S_i$ is greater than $1$;
\item Finally, otherwise, choose some point $(x, \delta_1) \in B_{\tilde{\epsilon}}((x, 0))$, so that $(x, \delta_1)$ is on the boundary of a subset of regions $S \subset \{S_1, \ldots, S_l\}$. We must have that 
$p_{comp}(x, \delta_1) = 0$ when computed using any of the algorithms that correspond to any $S_i \in S$. 

Let now $B(x)$ be a small ball around $x$ in $x$-space, and consider $\tilde{B}(x) \equiv B(x) \cap (\allow(x) \setminus V(p))$. 

There exists some $y \in \tilde{B}(x) $, close enough to $x$, such that $(y, \delta_1)$ is either in the interior or on the boundary of some $S_k \in S$.

By Lemma~\ref{prima}, since we must necessarily have $p_{comp}(y, \delta_1) = 0$ as computed by the algorithm corresponding to $S_k$, if follows (by continuity) that there is a small ball around $(y, \delta_1)$ on which the relative error, when computed using the algorithm corresponding to $S_k$, is greater than $\eta$. The intersection of that small ball with $S_k$ must have positive measure.
\end{enumerate}

\end{enumerate}

From the above analysis, it follows that there is always a set of positive measure, arbitrarily close to $(x, 0)$, on which the algorithm will produce a relative error larger than $\eta$.
\end{proof}

\vspace{.3cm}

\noindent \textit{Proof of Theorem \sl{\ref{conj}}.}  Follows immediately from Theorem~\ref{finala} and  Corollary~\ref{in_sfirsit}. \qed

\vspace{.3cm}

\begin{remark} Consider the polynomial $p(x,y) = (1-xy)^2+x^2$, whose variety is at infinity. We believe that Theorem~\ref{conj} can be extended to show that polynomials like $p(x,y)$ cannot be evaluated accurately on $\R$; this is future work.
\end{remark}



\subsection{Sufficiency: the complex case} \label{suf_c}

Suppose we now restrict input values to be complex numbers and use the same 
algorithm types and the notion of accurate evaluability from the previous 
sections. By Theorem~\ref{conj}, for a polynomial $p$ of $n$ complex 
variables to be accurately evaluable over $\C^n$ it is necessary that 
its variety $V(p)\eqbd 
\{z\in \C^n : p(z)=0 \}$ be allowable.

The goal of this section is that this condition is also sufficient, as stated in the following theorem.

\begin{theorem} \label{sufficiency_c}
Let $p: \C^n \to \C$ be a polynomial with integer
coefficients and zero constant term. Then $p$ is accurately
evaluable on $\mathcal{D} = \C^n$ if and only if the variety $V(p)$ is allowable.
\end{theorem}

To prove this we first investigate what allowable complex varieties can look like.
We start by recalling a basic fact about complex polynomial varieties, which can for example be deduced from Theorem 3.7.4 in~\cite[page 53]{T}. 
Let $V$ denote any complex variety. To say that $\dim_\C(V)=k$
means that, for each $z\in V$ and each $\delta>0$, there exists 
$w\in V\cap B(z,\delta)$ such that $w$ has a $V$-neighborhood that
is homeomorphic to a real $2k$-dimensional ball.

\begin{theorem} \label{dimensions} Let $p$ be a non-constant polynomial over $\C^n$. 
Then $$\dim_\C (V(p))=n-1.$$
\end{theorem}

\begin{corollary} \label{no_intrs} 
Let $p: \C^n \to \C$ be a nonconstant polynomial 
whose variety $V(p)$ is allowable. Then $V(p)$ is a union of allowable 
hyperplanes.  
\end{corollary}

\begin{proof} Suppose $V(p)=\cup_{j} S_j$, where each $S_j$ is an intersection of the sets in Definition~\ref{allwa} and, for some $j_0$, $S_{j_0}$
is not a hyperplane but an irredundant intersection of hyperplanes.
Let $z\in {S_{j_0}\setminus \cup_{j\neq j_0} S_j}$. Then, for some
$\delta>0$, $B(z,\delta)\cap V(p)\subset S_{j_0}$. Since $\dim_C(S_{j_0})<n-1$,
no point in $B(z,\delta)\cap V(p)$ has a $V(p)$-neighborhood that is 
homeomorphic to a real $2(n-1)$-dimensional ball. Contradiction.
\end{proof}
 
\begin{corollary} \label{factors} If $p: \C^n \to \C$ is a 
polynomial whose variety $V(p)$ is allowable, then it is a product 
$p=c \prod_j p_j$, where each $p_j$ is a power of $x_i$, $(x_i-x_j)$, or $(x_i +x_j)$.
\end{corollary}

\begin{proof} By Corollary~\ref{no_intrs}, the variety $V(p)$ is a union of allowable hyperplanes. Choose a hyperplane $H$ in that union. If $H = Z_{j_0}$ for some $J_0$, expand $p$ into a Taylor series in $x_{j_0}$. If $H =D_{i_0 j_0}$ (or $H = S_{i_0 j_0}$)  for some $i_0$, $j_0$, expand $p$ into a Taylor series in $(x_{i_0} - x_{j_0})$ (or  $(x_{i_0} + x_{j_0})$). In either case, in this expansion, the zeroth coefficient of $p$ must be the zero polynomial in $x_j$, $j\neq j_0$ (or $ j \notin \{i_0, j_0\}$). 
Hence there is a $k$ such that $p(x)= x_{j_0}^{k} ~\widetilde{p}(x)$ in the first case, or $p(x) = (x_{i_0} \pm x_{j_0})^k ~ \widetilde{p}(x)$ in the second (third) one. In any case,
we choose $k$ maximal, so that the variety $V(\widetilde{p})$ is the closure of the set
$V(p)\setminus Z_{j_0}$ in the first case, or $V(p)\setminus D_{i_0 j_0}$ ($V(p) \setminus S_{i_0, j_0}$) in the second (third) case. Then proceed by factoring $\widetilde{p}$ in the same fashion. 
\end{proof}

\noindent \textit{Proof of Theorem \sl{\ref{sufficiency_c}}.}   
By Corollary~\ref{factors},  $ p=c\prod_j p_j $, with each $p_j$ a power of $x_k$ or $(x_k \pm x_l)$. It also follows that $c$ must be an integer since all coefficients of $p$ are integers.

Since each
of the factors is accurately evaluable, and we can get any integer constant $c$ in front of $p$ by repeated addition (followed, if need be, by negation), which are again accurate operations, the algorithm that forms their product and then adds/negates to obtain $c$ evaluates $p$ accurately. 
\qed

\begin{remark} From Theorem {\sl \ref{sufficiency_c}}, it follows that only homogeneous polynomials are accurately evaluable over $\C^n$. \end{remark}



\subsection{Toward a necessary and sufficient condition in the real case} \label{suf_r}

In this section we show that accurate evaluability of a polynomial over
$\R^n$ is ultimately related to accurate evaluability of its ``dominant
terms''. This latter notion is formally defined later in this section. 
Informally, it describes the terms of the polynomial that dominate
the remaining terms  in a particular semialgebraic set close to a particular 
component of its variety; thus it depends on how we ``approach'' the
variety of a polynomial. 

For reasons outlined in Section \ref{sec_PositivePolys}, we consider here only homogeneous polynomials. Futhermore, most of this section is devoted to non-branching algorithms, but we do need branching for our statements at the end of the section. The reader will be alerted
to any change in our basic assumptions.  
 
Here is a short walk through this section:
\begin{itemize}
\item In Section \textbf{\ref{hom}. Homogeneity}, we discuss an expansion of the relative error $|p_{comp}(x, \delta) - p(x)|/|p(x)|$ as a function of $x$ and $\delta$, and prove a result about accurate evaluability of homogeneous polynomials that will be used in Section~\textbf{\ref{sec_prune}. Pruning}.
\item In Section \textbf{\ref{dom}. Dominance}, we introduce the notion of dominance and present different ways of looking at an irreducible component of the variety $V(p)$ using various simple linear changes of variables. These changes of variables allow us to identify all the dominant terms of the polynomial, together with the ``slices'' of space where they dominate.
\item In Section \textbf{\ref{sec_prune}. Pruning}, we explain how to ``prune'' an algorithm to manufacture an algorithm that evaluates one of its dominant terms, and prove a necessary condition for the accurate evaluation of a homogeneous polynomial by a non-branching algorithm. Roughly speaking, this condition says that accurate evaluation of the dominant terms we identified in Section~\textbf{\ref{dom}. Dominance}, is necessary.
\item In Section \textbf{\ref{sec_suff}. Sufficiency of evaluating dominant terms}, we identify a special collection of dominant terms, together with the slices of space where they dominate. If accurately evaluable by (branching or non-branching) algorithms, these dominant terms allow us to construct a branching  algorithm for the evaluation of the polynomial over the entire space (Theorem~\ref{pdom=>p}). These are just some of the terms present in the statement of Theorem~\ref{p=>pdom}. 
\end{itemize}

\subsubsection{Homogeneity} \label{hom}

We begin by establishing some basic facts about non-branching algorithms that 
evaluate homogeneous polynomials.

\begin{definition} \label{homalg}
We call an algorithm $p_{comp}(x,\delta)$ with error set $\delta$ for computing $p(x)$
{\em homogeneous of degree $d$\/} if 
\begin{enumerate} 
\item \label{unu1} the final output is of degree $d$ in $x$;
\item \label{doi2} no output of a computational node exceeds degree $d$ in $x$;
\item \label{trei3} the output of every computational node is homogeneous in $x$.
\end{enumerate}
\end{definition}

\begin{lemma} \label{1} If $p(x)$ is a homogeneous polynomial of degree $d$ and 
if a non-branching algorithm  evaluates $p(x)$ accurately by computing 
$p_{comp}(x,\delta)$, the algorithm must itself be homogeneous of degree $d$.
\end{lemma}

\begin{proof} First note that the output of the algorithm must be of degree
at least $d$ in $x$,  
since $p_{comp}(x,\delta)=p(x)$ when $\delta = 0$. Let us now write the overall
relative error as \[ rel_{err}(x,\delta) = { p_{comp}(x,\delta) - p(x) \over p(x)} =
\sum_{\alpha} \frac{p_{\alpha}(x)}{p(x)}~\delta^\alpha \] where $\alpha$
is a multi-index. If $p_{comp}(x,\delta)$ is accurate then $p_{\alpha}(x)/p(x)$
must be a bounded rational function on the domain ($\R^n$, or in
the homogeneous case the sphere $S^{(n-1)}$). This implies, in particular, 
that the output cannot be of degree higher than $d$ in $x$. So, Condition~1
of Definition~\ref{homalg} must be satisfied.

Now suppose Condition~\ref{trei3} of Definition~\ref{homalg} is
violated. We would like to show that the final output is also inhomogeneous.
We can assume without loss of generality that the algorithm does not
contain nodes that do operations like $x-x$ or $0 \cdot x$ (these can be
``pre-pruned'' and replaced with a $0$ source). There exists a highest node
$(op(\cdot), \delta_i)$ whose output is not homogeneous. If it is the
output node, we are done. Otherwise, look at the next node  
$(op(\cdot), \delta_j)$ on the path toward the output node. 
The output of $op(\cdot, \delta_j)$ is homogeneous. On the other hand,
the output of $(op(\cdot), \delta_i)$ (which is one of the two inputs to
$(op(\cdot), \delta_j)$) must be inhomogeneous in $x$ and must contain a
term $\delta_i r(x)$ with $r(x)$ an inhomogeneous polynomial in $x$.
 
If $op(\cdot,\delta_i)$ is the only input to $op(\cdot,\delta_j)$, then inhomogeneity 
will be present in both outputs, since neither doubling nor squaring can cancel it; 
contradiction. Otherwise there is another input to $op(\cdot,\delta_j)$
 (call it $op(\cdot,\delta_k)$). The output of $op(\cdot,\delta_k)$
must therefore also be inhomogeneous to cancel the inhomogeneous $r(x)$. 
Since the DAG is acyclic, $\delta_i$ is not 
present in the output of $op(\cdot,\delta_k)$ or $\delta_k$ is not present 
in the output $op(\cdot,\delta_i)$. Without loss of generality, assume the former case.
Then the term  $\delta_i r(x)$ will create inhomogeneity in the output of
$(op(\cdot), \delta_j)$, and hence $(op(\cdot), \delta_i)$ is
\emph{not} a highest node with inhomogeneous output, contradiction.
 Hence $p_{comp}(x,\delta)$ is not homogeneous in $x$, thus one
of the $p_{\alpha}(x)$'s has to contain terms in $x$ of higher or
smaller degree than $d$. 

Similarly, if Condition~\ref{doi2} of Definition~\ref{homalg} were
violated, then for some $\delta$s the final output would be a polynomial of
higher degree in $x$, and that would also mean some $p_{\alpha}(x)$ would
be of higher degree in $x$.

In either of these cases, if some $p_{\alpha}(x)$ contained terms of smaller 
degree than $d$, by scaling the variables appropriately and letting some of 
them go to $0$, we would deduce that $p_{\alpha}(x)/p(x)$ could not be bounded. 
If some $p_{\alpha}(x)$ contained terms of higher degree than $d$, by scaling 
the variables appropriately and letting some of them go to $\infty$, we would 
once again obtain that $p_{\alpha}(x)/p(x)$ could not be bounded.
\end{proof}

This proof shows that an algorithm evaluates a homogeneous polynomial $p$ accurately 
on $\R^n$ if and only if each fraction $p_\alpha/p$ is bounded on $\R^n$. It also
shows each $p_\alpha$ has to be homogeneous of the same degree as $p$. Therefore,
each fraction $p_\alpha/p$ is bounded on $\R^n$ if and only if it is bounded on
the unit sphere $S^{(n-1)}$. We record this as a corollary.

\begin{corollary} A non-branching homogeneous algorithm is accurate on $\R^n$ if and only
if it is accurate on $\S^{(n-1)}$. 
\end{corollary}

\subsubsection{Dominance} \label{dom}

Now we begin our description of ``dominant terms'' of a polynomial.
Given a polynomial $p$ with an allowable variety $V(p)$, let us fix 
an irreducible component of $V(p)$. Any such component is described by linear allowable
constraints, which, after reordering variables, can be grouped into 
$l$ groups as
$$ x_1=\cdots=x_{k_1}=0, \quad x_{k_1+1}=\cdots=\pm 
x_{k_2}, \;\; \ldots, \quad x_{k_{l-1}+1}=\cdots=
\pm x_{k_l}.    $$
To consider terms of $p$ that ``dominate'' in a neighborhood 
of that component, we will change variables to map any component 
of  a variety to a set of the form
\begin{equation} \tilde{x}_1=\cdots=\tilde{x}_{k_1}=0, \quad \tilde{x}_{k_1+2}=\cdots=
\tilde{x}_{k_2}=0, \quad \ldots, \quad \tilde{x}_{k_{l-1}+2}=\cdots=
\tilde{x}_{k_l}=0.  \label{new_var} 
\end{equation}
The changes of variables we will use are defined inductively as follows.

\begin{definition}  \label{chng_var}
We call a change of variables associated with a set of the
form  $$  \sigma_1 x_1=\sigma_2 x_2 =\cdots =\sigma_k x_k, \qquad \sigma_l=\pm 1, \;\; \;\; l=1, \ldots, k, $$
{\em basic\/} if it leaves one of the variables unchanged, which we will refer
 to as the {\em representative\/} of the group, and replaces the
remaining  variables by their sums (or differences) with the representative
of the group. In other words, 
$$  \widetilde{x}_j\eqbd x_j, \quad \widetilde{x}_l\eqbd x_l -\sigma_j \sigma_l x_j
\qquad {\rm for} \;\;\; l\neq j,  $$
where $x_j$ is the representative of the group $x_1, \ldots, x_k$.
A change of variables associated with a set of all $x$ satisfying conditions
\begin{equation} \begin{array}{l}
 x_1=\cdots=x_{k_1}=0, \\ \sigma_{k_1+1}x_{k_1+1}=\sigma_{k_1+2}x_{k_1+2}=\cdots=
\sigma_{k_2} x_{k_2}, \\
 \ldots \ldots \ldots \ldots \ldots \ldots \ldots \ldots \ldots \ldots \ldots \ldots \ldots \\ 
 \sigma_{k_{l-1}+1}x_{k_{l-1}+1}=\sigma_{k_{l-1}+2}x_{k_{l-1}+2}=\cdots=\sigma_{k_l} x_{k_l},  \\ \\
\sigma_j=\pm 1 \;\; \hbox{\rm for all pertinent} \; j  
\end{array}  \label{branch}
\end{equation}
is {\em basic\/} if it is a composition of the identity map on the first
$k_1$ variables and $(l-1)$ basic changes of variables associated 
with each set $\sigma_{k_1+1}x_{k_1+1}=\cdots=\sigma_{k_2}x_{k_2}$ through 
$\sigma_{k_{l-1}+1}x_{k_{l-1}+1}= \cdots=\sigma_{k_l} x_{k_l}$. 

Finally, a change of variables associated with a set $S$ of type~(\ref{branch})
is {\em standard\/} if it is a basic change of variables associated with
some allowable irreducible superset $\widetilde{S}\supseteq S$ and it maps 
$S$ to~(\ref{new_var}).
\end{definition}

Thus, a standard change of variables amounts to splitting the group $x_1, \ldots,
 x_{k_1}$ into smaller groups and either keeping the conditions $x_r=\cdots =x_q=0$ 
or assigning arbitrary signs to members of each group so as 
to obtain a set $\sigma_r x_r=\cdots = \sigma_q x_q$. It may also involve splitting the 
chains of  conditions $\sigma_{k_m+1}x_{k_m+1}=\cdots=\sigma_{k_{m+1}}x_{k_{m+1}}$ into 
several  subchains. The standard change of variables is then just one of the basic
changes of variables associated with the obtained set.

\begin{example} There are  $5\times 3$ standard changes of variables associated 
with  the set  $$ x_1=x_2=0, \qquad x_3=-x_4=x_5:$$
$$ \begin{array}{llllll}
 \tilde{x}_1=x_1, &  \tilde{x}_2=x_2, & \tilde{x}_3=x_3, 
& \tilde{x}_4=x_4+x_3, & \tilde{x}_5=x_5-x_3, & {\rm or} \\  
 \tilde{x}_1=x_1, &  \tilde{x}_2=x_2, & \tilde{x}_3=x_3+x_4, & 
\tilde{x}_4=x_4, & \tilde{x}_5=x_5+x_4, & {\rm or} \\  
 \tilde{x}_1=x_1, & \tilde{x}_2=x_2, & \tilde{x}_3=x_3-x_5, &
\tilde{x}_4=x_4+x_5, & \tilde{x}_5=x_5, & {\rm or} \\
 \tilde{x}_1=x_1, &  \tilde{x}_2=x_2-x_1, & \tilde{x}_3=x_3, 
& \tilde{x}_4=x_4+x_3, & \tilde{x}_5=x_5-x_3, & {\rm or} \\  
 \tilde{x}_1=x_1, &  \tilde{x}_2=x_2-x_1, & \tilde{x}_3=x_3+x_4, & 
\tilde{x}_4=x_4, & \tilde{x}_5=x_5+x_4, & {\rm or} \\  
 \tilde{x}_1=x_1, & \tilde{x}_2=x_2-x_1, & \tilde{x}_3=x_3-x_5, &
\tilde{x}_4=x_4+x_5, & \tilde{x}_5=x_5, & {\rm or} \\
 \tilde{x}_1=x_1, &  \tilde{x}_2=x_2+x_1, & \tilde{x}_3=x_3, 
& \tilde{x}_4=x_4+x_3, & \tilde{x}_5=x_5-x_3, & {\rm or} \\  
 \tilde{x}_1=x_1, &  \tilde{x}_2=x_2+x_1, & \tilde{x}_3=x_3+x_4, & 
\tilde{x}_4=x_4, & \tilde{x}_5=x_5+x_4, & {\rm or} \\  
 \tilde{x}_1=x_1, & \tilde{x}_2=x_2+x_1, & \tilde{x}_3=x_3-x_5, &
\tilde{x}_4=x_4+x_5, & \tilde{x}_5=x_5, & {\rm or} \\
 \tilde{x}_1=x_1-x_2, &  \tilde{x}_2=x_2, & \tilde{x}_3=x_3, 
& \tilde{x}_4=x_4+x_3, & \tilde{x}_5=x_5-x_3, & {\rm or} \\  
 \tilde{x}_1=x_1-x_2, &  \tilde{x}_2=x_2, & \tilde{x}_3=x_3+x_4, & 
\tilde{x}_4=x_4, & \tilde{x}_5=x_5+x_4, & {\rm or} \\  
 \tilde{x}_1=x_1-x_2, & \tilde{x}_2=x_2, & \tilde{x}_3=x_3-x_5, &
\tilde{x}_4=x_4+x_5, & \tilde{x}_5=x_5, & {\rm or} \\
 \tilde{x}_1=x_1+x_2, &  \tilde{x}_2=x_2, & \tilde{x}_3=x_3, 
& \tilde{x}_4=x_4+x_3, & \tilde{x}_5=x_5-x_3, & {\rm or} \\  
 \tilde{x}_1=x_1+x_2, &  \tilde{x}_2=x_2, & \tilde{x}_3=x_3+x_4, & 
\tilde{x}_4=x_4, & \tilde{x}_5=x_5+x_4, & {\rm or} \\  
 \tilde{x}_1=x_1+x_2, & \tilde{x}_2=x_2, & \tilde{x}_3=x_3-x_5, &
\tilde{x}_4=x_4+x_5, & \tilde{x}_5=x_5, & 
\end{array} $$
The supresets $\widetilde{S}$ for this example are the set $S$ itself
together with the set $\{x: x_1+x_2=0, x_3=-x_4=x_5 \}$ and the set 
$\{x : x_1-x_2=0, x_3=-x_4=x_5 \}$.
\end{example}

Note that we can write the vector of new variables $\tilde{x}$
as $Cx$ where $C$ is a matrix, so can label the change of variables 
by the matrix $C$. 

Now let us consider components of the variety $V(p)$. We have seen that
any given component of $V(p)$ can be put into the form
$x_1 = x_2 =...= x_k =0$ using a standard change of variables,
provided $V(p)$  is allowable. (To avoid cumbersome notation, 
we renumber all the variables set to zero as $x_1$ through $x_k$
for our discussion that follows. We will return to the original
description  to introduce the notion of pruning.)

Write the polynomial $p(x)$ in the form 
\begin{eqnarray} \label{expression}
p(x) = \sum_{\lambda \in \Lambda} c_{\lambda} x_{[1{:}k]}^{\lambda} q_{\lambda}(x_{[k+1:n]})~,
\end{eqnarray}
where, almost following MATLAB notation, we write $x_{[1{:}k]} \eqbd (x_1, \ldots, x_k)$, 
$x_{[k{+}1{:}n]}\eqbd (x_{k+1}, \ldots, x_n)$. Also, we let $\Lambda$ be the set of 
all multi-indices $\lambda \eqbd (\lambda_1, \ldots, \lambda_k)$ occuring in the 
monomials  of $p(x)$.

To determine all dominant terms associated with the component $x_1 = x_2 =...= x_k =0$,
consider the Newton polytope $P$ of the polynomial $p$ with respect to the 
variables $x_1$ through $x_k$ only,  i.e., the convex hull of the exponent vectors $\lambda \in \Lambda$ 
(see, e.g.,~\cite[p.~71]{MS}). Next, consider the normal fan $N(P)$ of $P$ (see~\cite[pp.~192--193]{Z}) 
consisting of the cones of all row vectors $\eta$ from the dual space $(\R^k)^*$ whose dot products with 
$x\in P$ are maximal for $x$ on a fixed face of $P$. That means that for every nonempty face $F$ of $P$  we take
$$ N_F\eqbd \{ \eta=(n_1, \ldots, n_k)\in (\R^k)^* : F\subseteq \{ x\in P : \eta x(\eqbd
 \sum_{j=1}^k n_j x_j )= \max_{y\in P}  \eta y  \}   \}  $$   
and $$ N(P)\eqbd \{ N_F : \;  F \;\hbox{\rm is a face of} \; P  \}.   $$

%

Finally, consider the intersection of the negative of the normal fan $-N(P)$ and the nonnegative  quadrant 
$(\R^k)^*_+$. This splits the first quadrant $(\R^k)^*_+$ into several regions
$S_{\Lambda_j}$ according to which subsets $\Lambda_j$ of exponents $\lambda$ 
``dominate'' close to the considered component of the variety $V(p)$, in the
following sense:

\begin{definition}
Let $\Lambda_j$ be a subset of $\Lambda$ that determines a face of the Newton
polytope $P$ of $p$ such that the negative of its normal cone $-N(P)$ 
intersects $(\R^k)^*_+$ nontrivially
(not only at the origin). Define $S_{\Lambda_j} \in (\R^k)^*_+$ to be the 
set of all nonnegative row vectors $\eta$ such that 
\[ \eta{\lambda_1} = \eta{\lambda_2} < \eta {\lambda}, ~~\forall \lambda_1, \lambda_2 \in \Lambda_j, ~~\mbox{and}~ \lambda \in \Lambda \setminus \Lambda_j. \]
\end{definition}


Note that if $x_1$ through $x_k$ are small, then the exponential change of variables 
$x_j \mapsto -\log |x_j|$ gives rise to a correspondence between the nonnegative part 
of $-N(P)$ and the space of original variables $x_{[1{:}k]}$.  
We map back the sets $S_{\Lambda_j}$ into a neighborhood of $0$ in
 $\R^k$ by lifting.\footnote{This is reminiscent of the concept 
of an amoeba introduced in~\cite{GKZ}.}

\begin{definition}
Let $F_{\Lambda_j} \subseteq [-1,1]^k$ be the set of all points $x_{[1:k]}\in \R^k$ such that 
\[\eta\eqbd (-\log |x_1|, \ldots, - \log|x_k|) \in S_{\Lambda_j}. \]
\end{definition}

\begin{remark} 
For any $j$, the closure of $F_{\Lambda_j}$ contains the origin in $\R^k$.
\end{remark}

\begin{remark} \label{rem_F} Given a point $x_{[1:k]} \in F_{\Lambda_j}$, and given $\eta=(n_1, n_2, \ldots, n_k) \in S_{\Lambda_j}$, for any $t\in (0,1)$, the vector $(x_1 t^{n_1}, \ldots, x_k t^{n_k})$ is in $F_{\Lambda_j} $.  Indeed, if $(-\log|x_1|, \ldots, -\log |x_k|)\in S_{\Lambda_j}$,
then so is  $(-\log|x_1|,\ldots,-\log|x_k|)-\log|t|\eta$, since all
equalities and inequalities that define $S_{\Lambda_j}$ will be preserved,
the latter because $\log|t|<0$. 
\end{remark}

\begin{example} \label{example} Consider the following polynomial
$$p(x_1,x_2,x_3)=x_2^8 x_3^{12}+x_1^2x_2^2x_3^{16}+x_1^8 x_3^{12}+
x_1^6x_2^{14}+x_1^{10}x_2^6 x_3^4.$$
We show below the Newton polytope $P$ of $p$ with respect to 
the variables $x_1$, $x_2$, its normal fan $N(P)$, the intersection
$-N(P)\cap R^2_+$, the regions $S_{\Lambda_j}$, and the regions $F_{\Lambda_j}$. 
\end{example}
\vskip 1cm

\begin{tabular}{lllll}
\multicolumn{2}{l}{\epsfxsize=6.5cm \epsfbox{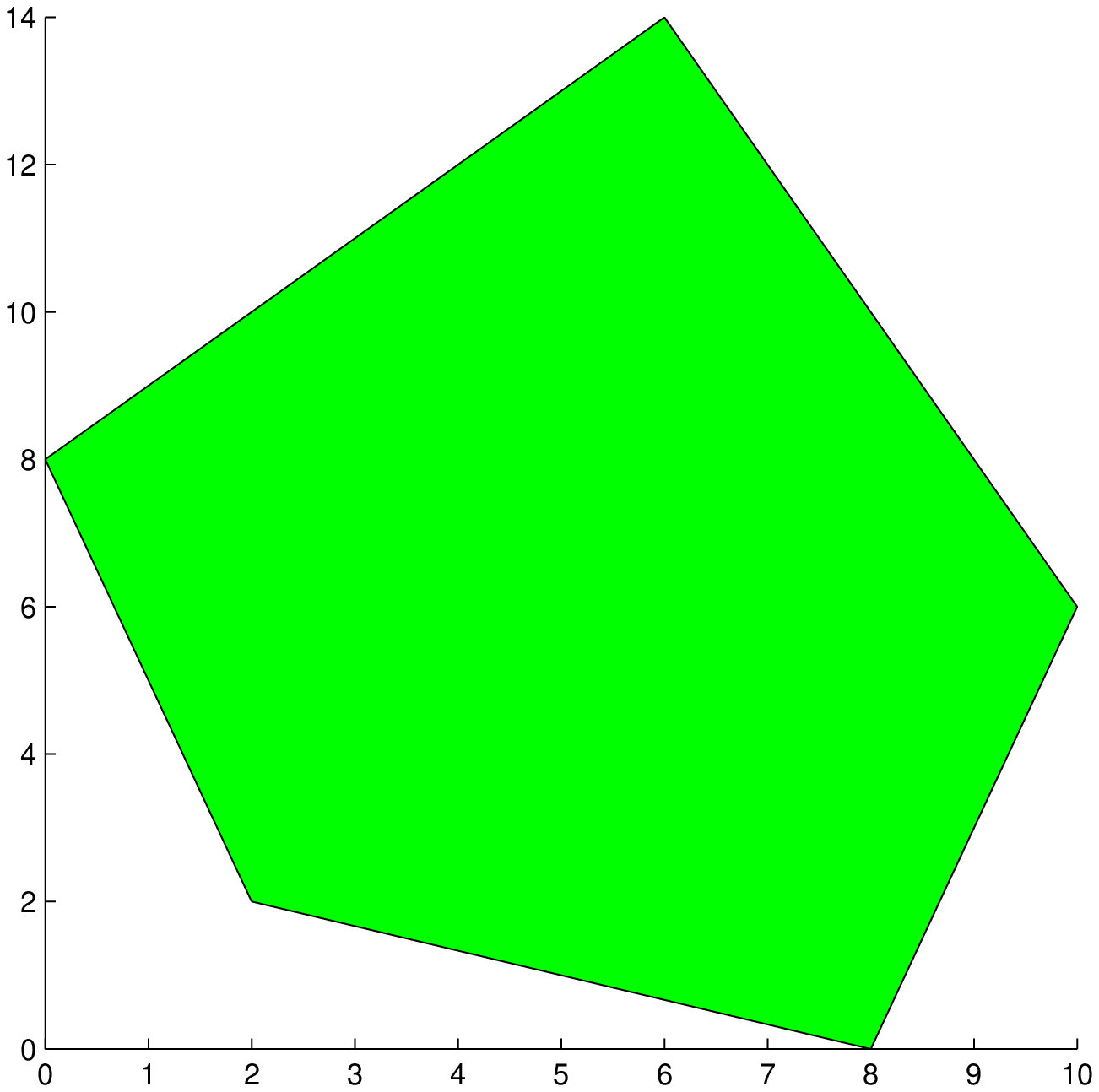}} & &
\multicolumn{2}{l}{\epsfxsize=7.0cm \epsfbox{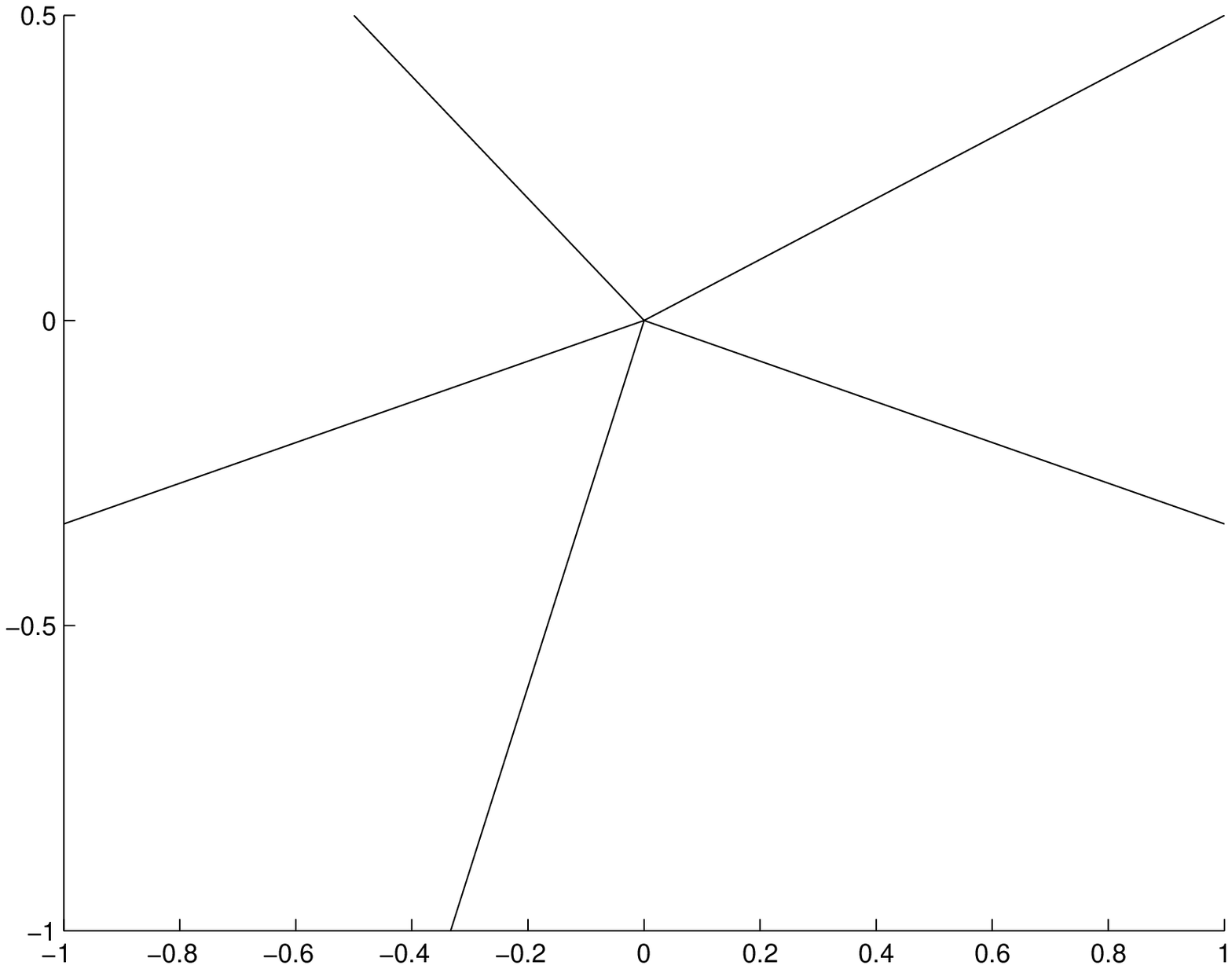}} \\
Figure 1. & A Newton polytope $P$. & & Figure 2. & Its normal fan $N(P)$. \\
\end{tabular} 
\vskip 1cm
 
\begin{tabular}{lllll}
\multicolumn{2}{l}{\epsfxsize=7.0cm \epsfbox{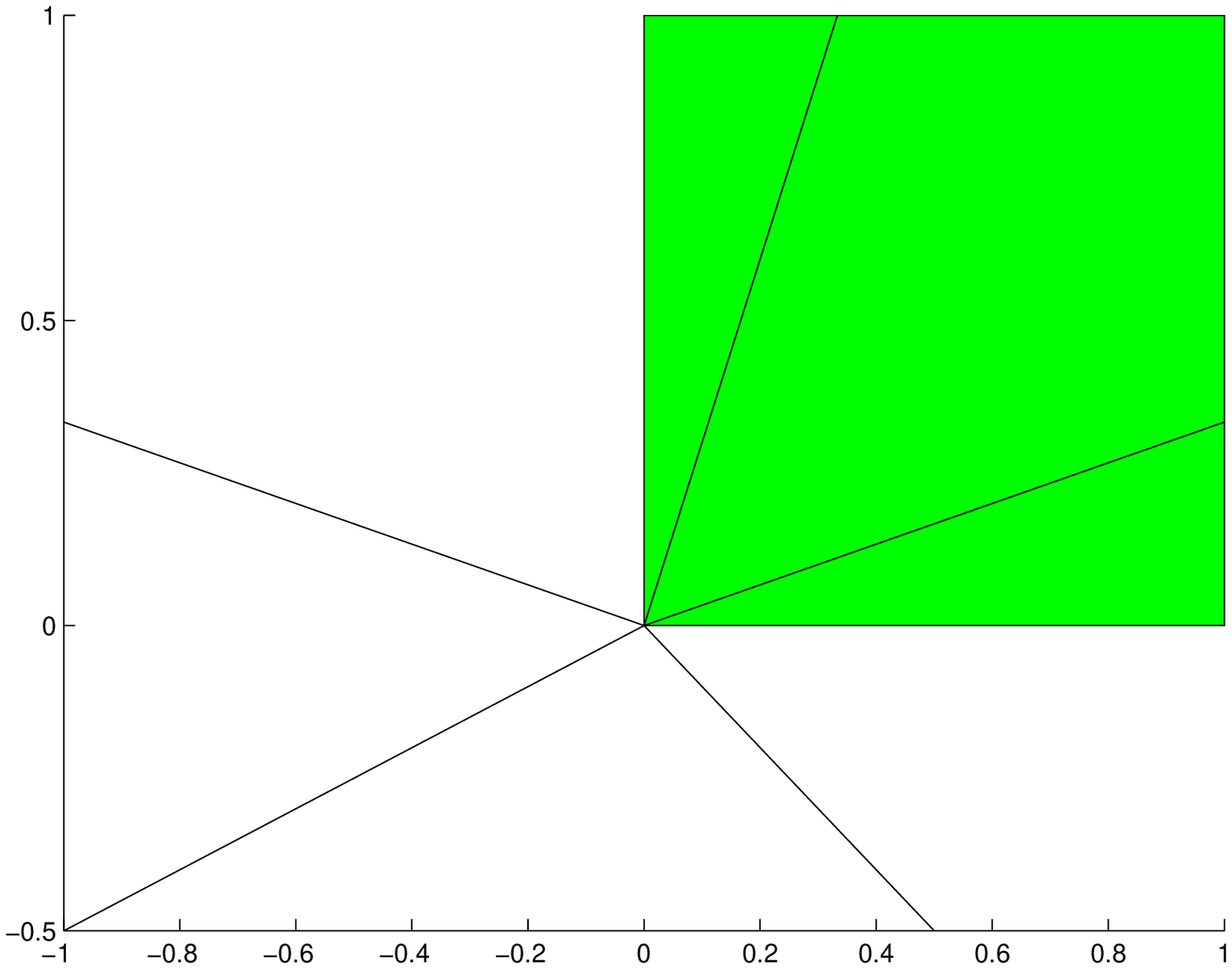}} & &
\multicolumn{2}{l}{\epsfxsize=7.0cm \epsfbox{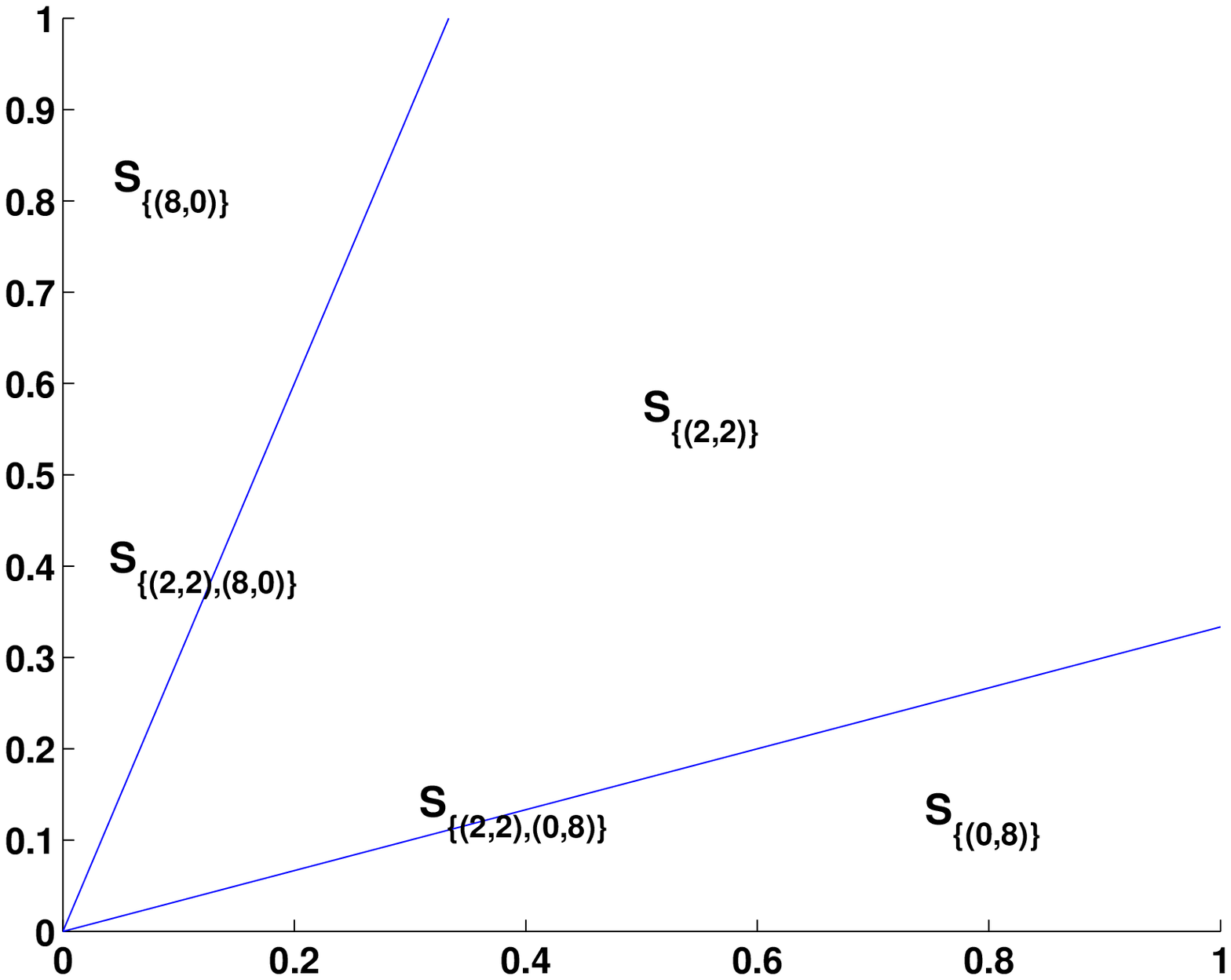}}  \\
Figure 3. & The intersection $-N(P)\cap \R^k_+$. & &
Figure 4. & The regions $S_{\Lambda_j}$.  \\
\end{tabular}
\vskip 1cm

\begin{tabular}{ll}
\multicolumn{2}{l}{\epsfxsize=13.0cm \epsfbox{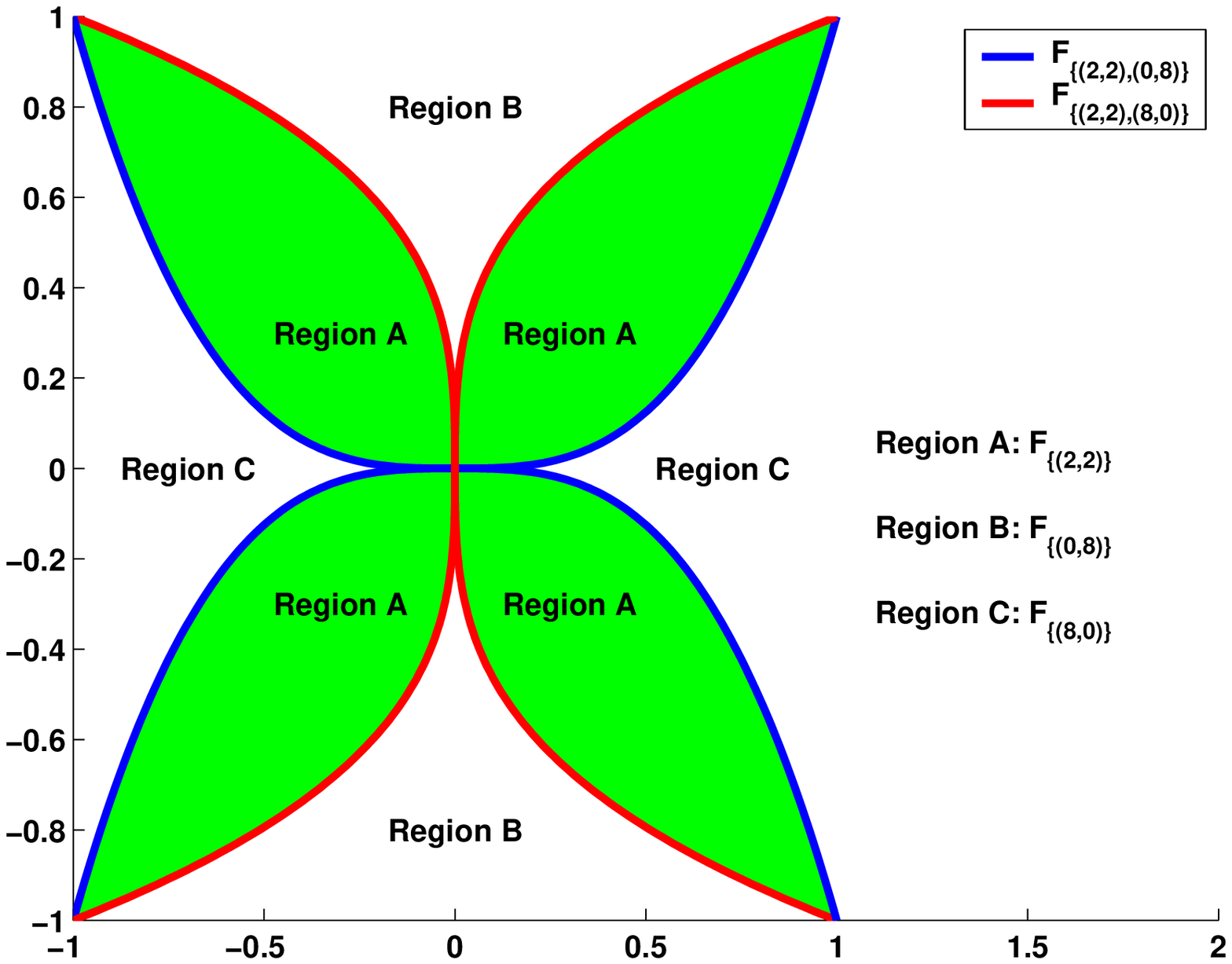}} \\
Figure 5. & The regions $F_{\Lambda_j}$.  \\ \\
\end{tabular}
\vspace{0.2cm}

\begin{definition}
%
We define the {\em dominant term\/} of $p(x)$ corresponding  to the component 
$x_1=\cdots=x_k=0$ and the region $F_{\Lambda_j}$ by
\[
p_{dom_j}(x) \eqbd \sum_{\lambda \in \Lambda_j} c_{\lambda} x_{[1:k]}^{\lambda} q_{\lambda}(x_{[k+1:n]})~.
\]  \end{definition}

The following observations about dominant terms are immediate. 

\begin{lemma} \label{leading} Let $\eta=(n_1, \ldots, n_k)\in S_{\Lambda_j}$ and let
$d_j\eqbd \sum_{\lambda_i\in \Lambda_j} \lambda_i n_i$.  Let $x^0$ be fixed and let 
$$ x(t)\eqbd (x_1(t), \ldots, x_n(t)), \qquad  x_j(t) \eqbd  \left\{ \begin{array}{ll} t^{n_j} x^0_j &  
j=1,\ldots, k,  \\ x^0_j, & j=k+1, \ldots, n. \end{array} \right. $$ 
Then $p_{dom_j}(x(t))$ has degree $d_j$ in $t$ and is the lowest degree 
term of $p(x(t))$ in $t$, that is
$$p(x(t))=p_{dom_j}(x(t))+o(t^{d_j})  \quad {\rm as}\;\; t\to 0, \qquad \deg_t p_{dom_j}(x(t))=d_j.  $$
\end{lemma}

\begin{proof} Follows directly from the definition of a dominant term. 
\end{proof}

\begin{corollary} Under the assumptions of Lemma~\ref{leading} suppose that
$p_{dom_j}(x^0)\neq 0$. Then 
$$\lim_{t\to 0} { p_{dom_j}(x(t)) \over  p(x(t))}=1.   $$  
\end{corollary}

Thus $p_{dom_j}$ is the leading term along each curve traced by $x(t)$ as $t$ tends to 
zero from above. An important question now is whether the dominant term $p_{dom_j}$
indeed dominates the remaining terms of $p$ in the region $F_{\Lambda_j}$ 
in the sense that $p_{dom_j}(x)/p(x)$ is close to $1$ sufficiently 
close to the component  $x_1=\cdots= x_k=0$ of the variety $V(p)$. This requires, at a 
minimum, that the variety $V(p_{dom_j})$ does not have a component strictly larger than
the set $x_1=\cdots= x_k=0$. Note that most dominant terms of a polynomial actually
fail this requirement. Indeed, most dominant terms of $p$ are monomials,
which correspond to regions $F_{\Lambda_j}$ indexed by singletons $\Lambda_j$, hence
to the vertices of the Newton polytope of $p$. The dominant terms corresponding to larger
sets $\Lambda_j$ are more useful, since they pick up terms relevant not only in the region
 $F_{\Lambda_j}$ but also in its neighborhood. In Example~\ref{example} above the dominant 
terms for $F_{\{(2,2),(8,0)\}}$ and $F_{\{(2,2),(0,8)\}}$, corresponding to the edges of 
the Newton polygon, are the useful ones. This points to the fact that we should be 
ultimately interested only in dominant terms corresponding to the facets, i.e., the
highest-dimensional faces, of the Newton polytope of $p$. Note that the convex hull
of $\Lambda_j$ is a facet of the Newton polytope $N$ if and only if the set
$S_{\Lambda_j}$ is a one-dimensional ray. 
 
The next lemma will be instrumental for our results in Section~\ref{sec_suff}.  It shows 
that each 
dominant term $p_{dom_j}$ such that the convex hull of $\Lambda_j$ is a facet of the 
Newton polytope of $p$  and whose variety $V(p_{dom_j})$ does not have a  component 
strictly larger than the  set $x_1=\cdots= x_k=0$ indeed dominates the remaining terms 
in $p$ in a certain ``slice'' $\widetilde{F}_{\Lambda_j}$  around $F_{\Lambda_j}$.

\begin{lemma} \label{true_dom} Let $p_{dom_j}$ be the dominant term of a homogeneous 
polynomial $p$ corresponding to the component $x_1=\cdots =x_k=0$ of the variety $V(p)$ 
and to the set $\Lambda_j$ whose convex hull is a facet of the Newton polytope $N$.
 
Let $\widetilde{S}_{\Lambda_j}$ be any closed pointed cone in $(\R^k)^*_+$ with vertex at
$0$ that does not intersect other one-dimensional rays $S_{\Lambda_l}$, $l\neq j$, and 
contains $S_{\Lambda_j} \setminus\{0\}$ in its interior. Let  $\widetilde{F}_{\Lambda_j}$
be the closure of the set
\begin{equation}
 \{x_{[1:k]} \in [-1,1]^k :  
 (-\log |x_1|, \ldots, - \log|x_k|) \in \widetilde{S}_{\Lambda_j} \}. \label{slices}
\end{equation} 
Suppose the variety $V(p_{dom_j})$ of $p_{dom_j}$ is allowable and intersects 
$\widetilde{F}_{\Lambda_j}$ 
only at $0$. Let $\| \cdot\|$ be any norm. Then, for any $\delta=\delta(j)>0$, there 
exists  $\varepsilon=\varepsilon(j)>0$  such that 
\begin{equation}
 \left|{p_{dom_j}(x_{[1:k]}, x_{[k+1:n]}) \over  p(x_{[1:k]}, x_{[k+1:n]})} -1  \right|<
\delta  \quad {\rm whenever} \;\; {\|x_{[1:k]}\|\over\|x_{[k+1:n]}\|}\leq \varepsilon \;\; {\rm and} \;\;
x_{[1:k]}\in \widetilde{F}_{\Lambda_j}.   \label{dom_inq}
\end{equation}
\end{lemma}

\begin{proof} We prove the lemma in the case $\widetilde{F}_{\Lambda_j}$ does not
intersect nontrivially any of the coordinate planes (the proof extends to the other 
case via limiting arguments).  
Let $\|x_{[k+1:n]}\|=1$ and let $x_1$ through $x_k$ be $\pm 1$.
If $\eta=(n_1, \ldots, n_k) \in \widetilde{S}_{\Lambda_j}$, then, directly from 
the definition of the set $\widetilde{F}_{\Lambda_j}$, the curve $(t^{n_1}x_1, \ldots, 
t^{n_k}x_k)$, $t\in (0,1]$, lies in $\widetilde{F}_{\Lambda_j}$ (and every point in
$\widetilde{F}_{\Lambda_j}$ lies on such a curve).
Denote $(t^{n_1}x_1, \ldots, t^{n_k}x_k, x_{[k+1:n]})$ by $x(t)$ and let $t$ decrease
from $1$ to $0$, keeping the $x_m$, $m=1, \ldots, n$, fixed. 
By the assumption of the Lemma,  $p_{dom_j}(x(t))$  does not vanish for sufficiently 
small $t>0$. Moreover, by Lemma~\ref{leading}, $p_{dom_j}$ is the leading term of $p$
in $F_{\Lambda_j}$. Since the cone $\widetilde{S}_{\Lambda_j}$ around $S_{\Lambda_j}$
does not intersect any other one-dimensional rays $S_{\Lambda_l}$, $l\neq j$, 
all the monomials present in any term that dominates in $\widetilde{F}_{\Lambda_j}
\setminus F_{\Lambda_j}$ are already present in $p_{dom_j}$. Thus $p_{dom_j}$ contains 
all terms that dominate in $\widetilde{F}_{\Lambda_j}$.  Therefore, there exists 
$\varepsilon(x)>0$ such that  $|p_{dom_j}(x(t))/p(x(t))-1|<\delta$ whenever 
$t<\varepsilon(x)$. The function $f: x\to \varepsilon(x)$ is lower semicontinuous.
Since the set $S\eqbd \{ x : x_m=\pm 1, \; m=1, \ldots, k, \; \|x_{[k+1:n]}\|=1 \}$
is compact, the minimum $\varepsilon\eqbd \min f(S)$ is necessarily positive and
satisfies~(\ref{dom_inq}).   \end{proof}

The above discussion of dominance was based on the transformation of a given irreducible
component of the variety to the form $x_1=\cdots =x_k=0$. We must reiterate that 
the identification of dominant terms becomes possible only after a suitable change 
of variables $C$ is used to put a given irreducible component into the standard 
form $x_1=\cdots =x_k=0$ and then the sets $\Lambda_j$ are determined. Note however 
that the polynomial $p_{dom_j}$ is given in terms of the original variables,
i.e., as a sum of monomials in the original variables $x_q$ and sums/differences 
$x_q\pm x_r$. We will therefore use the more precise notation $p_{dom_j,C}$ in the sequel. 

Without loss of generality we can assume that any standard change of variables 
has the form 
\begin{equation}
\begin{array}{l}  x=(x_{[1:k_1]}, x_{[k_1+1:k_2]}, \ldots, x_{[k_{l-1}+1:k_l]}) \mapsto 
\widetilde{x}=(\widetilde{x}_{[1:k_1]}, \widetilde{x}_{[k_1+1:k_2]}, \ldots, 
\widetilde{x}_{[k_{l-1}+1:k_l]}), \;\; {\rm where} \\ 
\widetilde{x}_{k_m+1}\eqbd x_{k_m+1}, \;\; \widetilde{x}_{k_m+2}\eqbd 
x_{k_m+2} -\sigma_{k_m+2} x_{k_m+1}, \;\; \ldots, \;\; \widetilde{x}_{k_{m+1}}
\eqbd x_{k_{m+1}} -\sigma_{k_{m+1}} x_{k_{m+1}}, \\
k_0\eqbd 0, \quad \sigma_r=\pm 1 \;\;\; \hbox{\rm for all pertinent}\;\; r 
 \end{array} 
\label{change*}
\end{equation}
Note also that we can think of the vectors $\eta\in S_{\Lambda_j}$ as being indexed by 
integers $1$ through $k_l$, i.e., $\eta=(n_1, \ldots, n_{k_l})$. Moreover, to define 
pruning in the next subsection we will assume that
\begin{equation}
n_{k_m+1}\leq n_r \quad \hbox{\rm for all}\;\; \;r=k_m+2, \ldots, k_{m+1}
\quad \hbox{\rm and for all} \;\; m=0, \ldots, l-1.
\label{exp_cond}
\end{equation} 

\begin{remark}  \label{WLOG}
 This condition is trivially satisfied if 
$n_{k_m+1}=0$, as is the case for any group $x_{k_m+1}=\sigma_{k_m+2} 
x_{k_m+2}=\cdots= \sigma_{k_{m+1}} x_{k_{m+1}}$ of original conditions
that define the given irreducible component of $V(p)$, since $x_{k_m+1}$
does not have to be close to $0$ in the neighborhood of that component of
$V(p)$. If, however, the same group of equalities was created from the original 
conditions  $x_{k_m+1}= x_{k_m+2}=\cdots= x_{k_{m+1}}=0$ due to the particular
change of variables $C$, the condition~(\ref{exp_cond}) is no longer 
forced upon us.  Yet~(\ref{exp_cond}) can be assumed without loss of
generality. Indeed, if, say, $n_{k_m+2}<n_{k_m+1}$, then we can always switch
to another standard change of variables by taking $x_{k_m+2}$ to be the 
representative of the group $x_{k_m+1}, \ldots, x_{k_{m+1}}$ and taking the 
sums/differences with $x_{k_m+2}$ as the other new variables.
Also note that~(\ref{exp_cond}) is satisfied either by all or by
no vectors in  $S_{\Lambda_j}$. In other words,~(\ref{exp_cond}) is a 
property of the entire set $S_{\Lambda_j}$. So, with a slight abuse
of terminology we will say that a set $S_{\Lambda_j}$ satisfies 
or fails~(\ref{exp_cond}). 
\end{remark}

Finally note that the curves $(x(t))$
corresponding to the change of variables~(\ref{change*}) are described
as follows:
\begin{equation}
\begin{array}{l}
x(t)\eqbd (x_{[1:k_1]}(t), x_{[k_1+1:k_2]}(t), \ldots, x_{[k_{l-1}+1:k_l]}(t),
x_{[k_l+1:n]}), \;\; {\rm where} \\

x_{[k_m+1:k_{m+1}]}(t)\eqbd \\ 
\;\; (t^{n_{k_m+1}}x_{k_m+1},  t^{n_{k_m+2}}x_{k_m+2}+
\sigma_{k_m+2}t^{n_{k_m+1}} x_{k_m+1}, \ldots,  t^{n_{k_{m+1}}}x_{k_{m+1}}+ \sigma_{k_{m+1}} 
t^{n_{k_m+1}}x_{k_m+1} ) \\
{\rm where} \;\;\; k_0\eqbd 1,\;\;\;  m=0, \ldots, l. 
\end{array}
\label{input}
\end{equation}
This description will be instrumental in our discussion of pruning, which follows immediately. 
 
\subsubsection{Pruning}  \label{sec_prune}

Now we discuss how to convert an accurate algorithm that evaluates a polynomial $p$
into an accurate  algorithm that evaluates a selected dominant term $p_{dom_j,C}$. 
This process, which we will refer to as {\em pruning,\/}  will consist of deleting some 
vertices and edges and redirecting certain other edges in the DAG that represents the 
algorithm.

\begin{definition}[Pruning] \label{def_prune} 
Given a non-branching algorithm represented by a DAG for computing 
$p_{comp}(x,\delta)$, a standard change of variables $C$ of the form~(\ref{change*})
and a subset $\Lambda_j \in \Lambda$ satisfying~(\ref{exp_cond}), we choose any 
$\eta\in S_{\Lambda_j}$, we input (formally) the expression~(\ref{input}), 
and then perform the following process. 

We can perform one of two actions: {\em redirection\/} or\/
{\em deletion.\/} By {\em redirection\/} (of sources) we mean replacing an edge from a source node corresponding
to a variable $x_j$ to  a computational node $i$ by an edge from the representative $x_{rep}$
 of $x_j$ followed by exact negation if $\sigma_j=-1$. This corresponds to replacing
$x_j$ by the product $\sigma_j x_{rep}$. To define {\em deletion,\/} consider a node $i$ 
with distinct input nodes $j$ and $k$. Then deletion of node $i$ from node $j$ means  
deleting the out-edge to node $i$ from node $j$, changing the origin  of all out-edges 
from node $i$ to input node $k$, and deleting node $i$. 
 
Starting at the sources, we process each node as follows, provided that both its 
inputs have already been processed (this can be done because of acyclicity). 
Let the node being processed be node $i$, i.e.,  $(op(\cdot),\delta_i)$, and assume 
it has input nodes $k$ and $l$. Both inputs being polynomials in $t$, we determine 
the lowest degree terms in $t$ present in either of them and denote these degrees
by $\deg(k)$ and $\deg(l)$.

\begin{description}

\item[if $op(\cdot)=\cdot$ and one or both inputs are sources,] then
\begin{description}
\item  redirect each source input.
\end{description}

\item[if $op(\cdot)=\pm$,] then 

\begin{description}

\item[if $\deg(k)\neq \deg(l)$, say $\deg(k)>\deg(l)$,] delete input node $i$ from node $k$.

\item[else] If nodes $k$ and $l$ are sources and the operation $op(\cdot)$ leads to 
cancellation of their lowest degree terms in $t$, examine their second-lowest degree terms.
If those degrees coincide or if one second-lowest term is missing, we change nothing. If 
one is bigger than the other, we do not 
change the source containing the lower degree term in $t$, but redirect the other source.

If only one of nodes $k$ and $l$ is a source or if both inputs are sources but there is no 
cancellation of lowest degree terms, redirect each source.


\end{description}
\end{description}

We then delete inductively all nodes which no longer are on any path to the output. 

We call this process {\em pruning,\/} and denote the output of the pruned
algorithm by  $p_{dom_j,C, comp}(x,\delta)$. \end{definition}

\begin{remark} Note that the outcome of pruning does not depend on the choice 
of $\eta\in S_{\Lambda_j}$. Since each region $S_{\Lambda_j}$ is determined
by linear homogeneous  equalities and inequalities with integer coefficients, 
the vector $\eta$ can  always be chosen to have all integer entries.
\end{remark}

\begin{example} Figure~6 shows an example of pruning an algorithm that evaluates the 
polynomial  $$   x_1^2 x_2^2 + (x_2-x_3)^4+(x_3-x_4)^2x_5^2  $$
using the substitution $$  (tx_1,x_2,tx_3+x_2, tx_4+x_2,x_5) $$
near the component $$ x_1=0, \;\; x_2=x_3=x_4.$$

\vskip 1cm

\begin{tabular}{ll}
\multicolumn{2}{l}{\epsfxsize=19.0cm \epsfbox{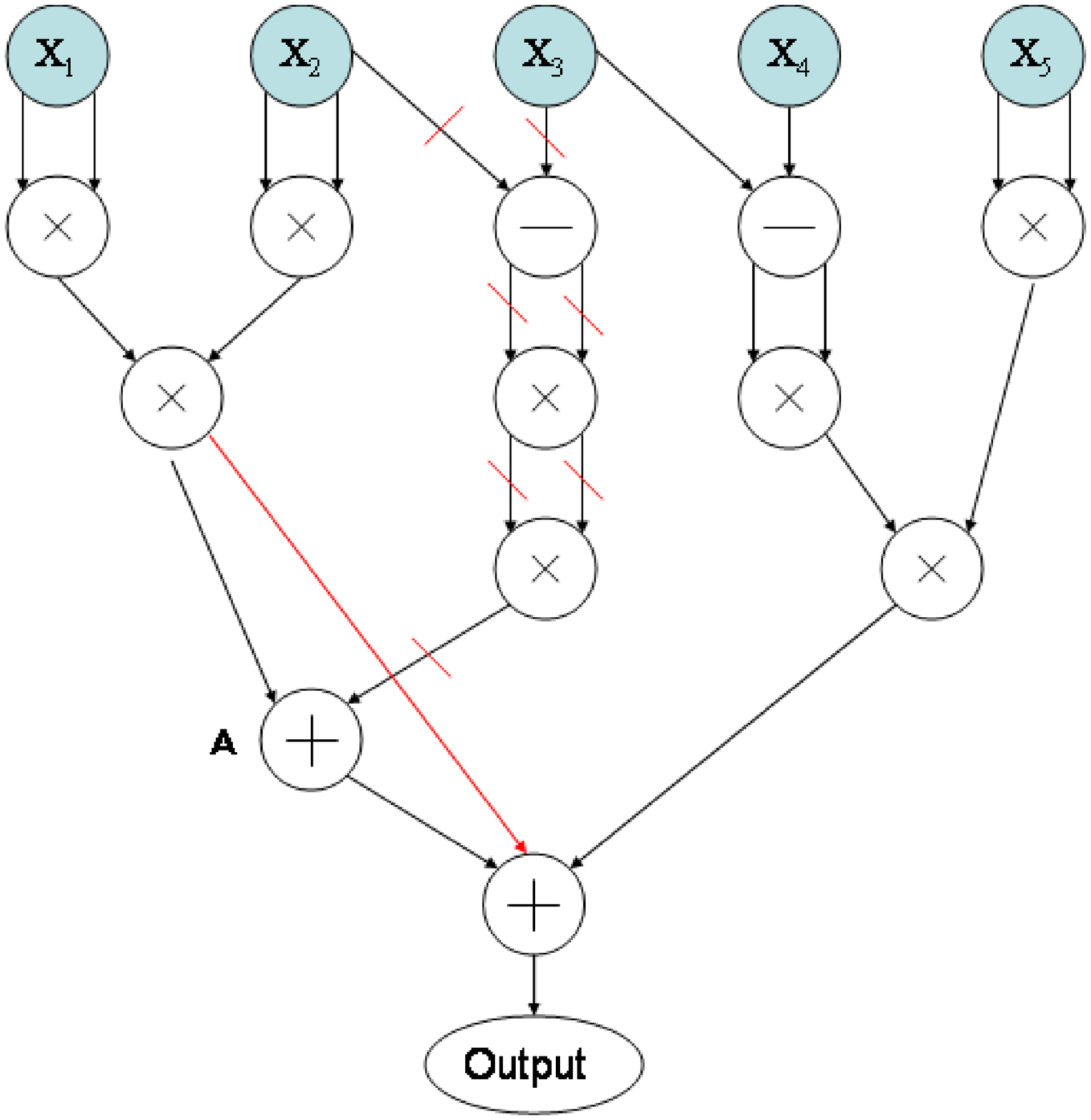}} \\
Figure 6. & Pruning an algorithm for $p(x)= x_1^2 x_2^2 + (x_2-x_3)^4+(x_3-x_4)^2x_5^2 $.  \\ \\
\end{tabular}
\vspace{0.2cm}

\noindent The result of pruning is an algorithm that
evaluates the dominant term $$ x_1^2x_2^2+(x_3-x_4)^2x_5^2.$$ 
One of two branches leading to the node $A$ is pruned due to the fact that
it computes a quantity of order $O(t^4)$ whereas the other branch produces
a quantity of order $O(t^2)$. 

The output of the original algorithm is given by  
\begin{eqnarray*}
&& 
 \left(  \left(   x_1^2(1+\delta_1)x_2^2(1+\delta_2) (1+\delta_3)+
(x_2-x_3)^4(1+\delta_4)^4(1+\delta_5)^2(1+\delta_6)\right)(1+\delta_7) \right. \\
&& \qquad \left. + (x_3-x_4)^2(1+\delta_8)^2(1+\delta_9)x_5^2(1+\delta_{10})(1+\delta_{11})\right)
(1+\delta_{12}).
\end{eqnarray*}
The output of the pruned algorithm is
\begin{eqnarray*}
 \left(   x_1^2(1+\delta_1)x_2^2(1+\delta_2) (1+\delta_3) 
 + (x_3-x_4)^2(1+\delta_8)^2(1+\delta_9)x_5^2(1+\delta_{10})(1+\delta_{11})\right)
(1+\delta_{12}).
\end{eqnarray*} \end{example}

%

Let us prove that this process will indeed produce an algorithm that accurately
evaluates the corresponding dominant term.


%
%

\begin{theorem} \label{p=>pdom}  Suppose a non-branching algorithm evaluates
a polynomial $p$ accurately on $\R^n$ by computing $p_{comp}(x,\delta)$. Suppose
$C$ is a standard change of variables~(\ref{change*}) associated with an
irreducible component of $V(p)$. Let $p_{dom_j,C}$ be one of the corresponding
dominant terms of $p$ and let $S_{\Lambda_j}$ satisfy~(\ref{exp_cond}). Then the 
pruned algorithm defined in Definition~\ref{def_prune} 
with output  $p_{dom_j,C, comp}(x,\delta)$ evaluates $p_{dom_j,C}$ accurately on $\R^n$. 
In other  words, being able to compute all such $p_{dom_j,C}$ for all components of the 
variety $V(p)$  and all standard changes of variables $C$ accurately is a condition 
necessary  to compute $p$ accurately.
\end{theorem}

\begin{proof} Directly from the definition of pruning it can be seen that the output of
each computational node is a homogeneous polynomial in $t$. This can be checked
inductively starting from computational nodes operating on two sources, using the pruning 
rules. Moreover, the pruning rules are equivalent to taking the lowest degree terms 
in $t$ (as well as setting some $\delta$s to zero). This can be checked inductively as 
well, once we rule out the situation when a $\pm$ node in the original algorithm leads 
to exact cancellation of lowest degree terms of the inputs, and at least one of the inputs
is not a source. Indeed, in that case one of the inputs contains a factor $(1+\delta)$
and that $\delta$ by acyclicity is not present in the other input. Therefore no
exact cancellation of lowest degree terms can occur.

Thus, the final output $p_{dom_j,C, comp}(x,\delta)$ of the pruned algorithm takes the 
lowest degree terms in $t$  of the final output of the original algorithm, so
 $p_{dom_j,C, comp}(x,\delta)$ is homogeneous in $t$ (of degree $d_j =\eta \lambda=\sum
\lambda_i n_i$). We write
\[\frac{p_{comp}(x,\delta) - p(x)}{p(x)} = \sum_{\alpha} \frac{p_{\alpha}(x)}{p(x)} \delta^{\alpha}, \qquad
\frac{p_{dom_j,C, comp}(x) - p_{dom_j,C}(x)}{p_{dom_j,C}(x)} = \sum_{\alpha} \frac{p_{\alpha, dom_j,C}(x)}{p_{dom_j,C}(x)} \delta^{\alpha}.\]
Note that, in eliminating nodes and redirecting the out-edges in the pruning process, we do the equivalent of setting some of the $\delta$s to $0$, and hence the set of monomials 
$\delta^\alpha$ present in the second sum is a subset of the set of monomials 
$\delta^\alpha$ present in the first sum.

Indeed, first of all, we can focus on the effect of deleting some nodes, since redirection does not affect $\delta$s at all, because only sources can be redirected. So, if $\delta^\alpha$ appears in the second sum, there is a path which yields the
corresponding multi-index, on which some term in $p_{\alpha, dom_j,C}(x)$ is computed. But 
since this path survived the pruning, there is a corresponding path in the original DAG 
which perhaps has a few more nodes that have been deleted in the pruning process
and a few source nodes that were redirected. In the 
computation, the effect of deleting a node was to set that $\delta$ equal to $0$ (and 
make some terms of higher degree disappear). So the surviving term was also present in 
the computation of $p_{comp}(x,\delta)$, with the same multi-index: just choose the $1$ 
in the $(1+\delta)$ each time when you hit a node that will be deleted (i.e., whose 
$\delta$  will be set to $0$).

Now note that $p_{\alpha, dom_j,C}(x)$ is the leading term of $p_{\alpha}(x)$, i.e., 
the term of smallest degree $d_j$ in $t$. This happens since each term of degree $d_j$ in $p_{\alpha}(x)$ must survive on the same path in the DAG, with the same choices of $1$ in $(1+\delta)$ each time we hit a deleted node.

We can now prove that $p_{dom_j, C, comp}(x)$ is accurate. To do that, it is enough to show that each
 $p_{\alpha, dom_j, C}(x)/p_{dom_j,C}(x)$ is bounded, provided that there is some constant $M$ such 
that $|p_{\alpha}(x)/p(x)| \leq M$ for all $x$.

Choose a point $x=(x_1, \ldots, x_n)$ not on the variety of $p_{dom_j,C}$ and consider the curve 
traced by the associated point $x(t)$ from~(\ref{input}) as $t$ tends to $0$. Since both $p_{\alpha, 
dom_j,C}(x(t))$ and $p_{dom_j,C}(x(t))$ are homogeneous of degree $d_j$ in $t$, we have 
\[
\left| \frac{p_{\alpha, dom_j,C}(x)}{p_{dom_j,C}(x)} \right|  = 
\left| \frac{p_{\alpha, dom_j,C}(x(t))}{p_{dom_j,C}(x(t))} \right| =  
\lim_{t \rightarrow 0} \left | \frac{p_{\alpha, dom_j,C}(x(t))}{p_{dom_j,C}(x(t))} \right | ~.
\]
Since both $p_{\alpha, dom_j,C}(x(t))$ \emph{and} $p_{dom_j,C}(x(t))$ are the dominant terms in 
$p_{\alpha}(x)$, respectively $p(x)$, along the curve $\{x(t): t\to 0\}$, we conclude that
\[ \left | \frac{p_{\alpha, dom_j,C}(x)}{p_{dom_j,C}(x)} \right|  =\lim_{t \rightarrow 0} \left | \frac{p_{\alpha, dom_j,C}(x(t))}{p_{dom_j,C}(x(t))} \right |= \lim_{t \rightarrow 0} \left | \frac{p_{\alpha}(x(t))}{p(x(t))} \right | \leq M~.
\]  
Invoking the density of the Zariski open set $\{ x: p_{dom_j,C}(x)\neq 0 \}$, we are done.
\end{proof}

\subsubsection{Sufficiency of evaluating dominant terms}  \label{sec_suff}

Our next goal is to prove a converse of a sort to Theorem~\ref{p=>pdom}. Strictly
speaking, the results that follow do not provide a true converse, since branching
is needed to construct an algorithm that evaluates a polynomial $p$ accurately
from algorithms that evaluate its dominant terms accurately.

For the rest of this section, we make two assumptions, viz., 
that our polynomial $p$ is homogeneous and irreducible. The latter assumption
effectively reduces the problem to that of accurate evaluation of a nonnegative
polynomial, due to the following lemma.

\begin{lemma} If a polynomial $p$ is irreducible and has an allowable variety $V(p)$,
then it is either a constant multiple of a linear form that defines an allowable
hyperplane or it does not change its sign in $\R^n$.
\end{lemma}

\begin{proof}
Suppose that $p$ changes its sign. Then the sets $\{x : p(x)>0\}$ and 
$\{ x : p(x)<0 \}$ are both open, hence the part of the variety $V(p)$
whose neighborhood contains points from both sets must have dimension $n-1$.
The only allowable sets of dimension $n-1$ are allowable hyperplanes. 
Therefore, the linear polynomial that defines an allowable hyperplane
must divide $p$. As $p$ is assumed to be irreducible, $p$ must be a constant
multiple of the linear polynomial. 

Thus, unless $p$ is a constant multiple of a linear factor of the type
$x_i$ or $x_i\pm x_j$, it must satisfy $p(x)\geq 0$ for all $x\in \R^n$
or $p(x)\leq 0$ for all $x\in \R^n$.
\end{proof}

From now on we therefore restrict ourselves to the nontrivial case when
a (homogeneous and irreducible) polynomial $p$ is nonnegative everywhere in $\R^n$.

\begin{theorem} \label{pdom=>p} Let $p$ be a homogeneous nonnegative polynomial 
whose variety $V(p)$ is allowable. Suppose that all dominant terms $p_{dom_j,C}$ 
for all components of the variety $V(p)$, all standard changes of variables $C$ 
and all subsets $\Lambda_j$ satisfying~(\ref{exp_cond}) are accurately evaluable. 
Then there exists a branching  algorithm that evaluates $p$ accurately over $\R^n$.  
\end{theorem}

\begin{proof} 
We first show how to evaluate $p$ accurately in a neighborhood of each
irreducible component of its variety $V(p)$. We next evaluate $p$ accurately
off these neighborhoods of $V(p)$. The final algorithm will involve branching
depending on which region the input belongs to, and the subsequent execution 
of the corresponding subroutine. We fix the relative accuracy $\eta$ that we
want to achieve. 

Consider a particular irreducible component $V_0$ of the variety $V(p)$. 
Using any standard change of variables $C$, say, a basic change of variables associated 
with $V_0$, we map $V_0$ to a set of the
form $\widetilde{x}_1=\cdots = \widetilde{x}_k=0$. Our goal is to create
an $\varepsilon$-neighborhood of $V_0$ where we can evaluate $p$ accurately.
It will be built up from semialgebraic $\varepsilon$-neighborhoods. We begin
with any set $S$ containing a neighborhood $V_0$, say, we let $S$ coincide with 
$[-1,1]^k \times \R^{n-k}$ after the change of variables $C$. 
We partition the cube $[-1,1]^k$ into sets $\widetilde{F}_{\Lambda_j}$ of 
type~(\ref{slices}) as follows:  We consider, as in Section~\ref{dom} above, 
the Newton polytope $P$ of $p$, form the intersection of the negative of its normal 
fan $-N(P)$ with the nonnegative quadrant $\R_+^k$ (in the new coordinate system),
and determine the sets $S_{\Lambda_j}$. If condition~(\ref{exp_cond}) fails for
some of the sets $S_{\Lambda_j}$, we transform them using a suitable standard 
change of variables as described in Remark~\ref{WLOG} so as to meet 
condition~(\ref{exp_cond}). For the rest of the argument, we can assume that all
sets $S_{\Lambda_j}$ satisfy~(\ref{exp_cond}).

To form conic neighborhoods of one-dimensional rays of $-N(P)\cap \R_+^k$ (which are normal 
to facets of the Newton polytope), we intersect $-N(P)\cap \R_+^k$ with, say, the 
hyperplane $\widetilde{x}_1+\cdots+\widetilde{x}_k=1$. Perform the Voronoi tesselation~(see, e.g.,~\cite{Z}) of the simplex $\widetilde{x}_1+\cdots+\widetilde{x}_k=1$, $\widetilde{x}_j\geq 0$, $j=1,\ldots, k$ relative to the intersection
points of $-N(P)\cap \R_+^k$ with the hyperplane $\widetilde{x}_1+\cdots+\widetilde{x}_k=1$. Connecting each Voronoi cell of the tesselation to the origin 
$\widetilde{x}_1=\cdots \widetilde{x}_k=0$ by straight rays, we obtain cones 
$\widetilde{S}_{\Lambda_j}$ and the corresponding sets $\widetilde{F}_{\Lambda_j}$
 of type~(\ref{slices}).  Note that the Voronoi cells and therefore the cones
$\widetilde{S}_{\Lambda_j}$ are determined by rational inequalities since the tesselation 
centers have rational coordinates. Hence the sets  $\widetilde{F}_{\Lambda_j}$ are
semialgebraic, and moreover are determined by polynomial inequalities with integer 
coefficients. Indeed, even though the sets $\widetilde{F}_{\Lambda_j}$ are defined using logarithms,
the resulting inequalities are among powers of absolute values of the variables and/or their 
sums and differences. For example, if a particular set  $\widetilde{F}_{\Lambda_j}$ is described
by the requirement that $(-\log|x_1|, -\log|x_2|)$ lie between two lines through thw origin, 
with slopes $1/2$ and $2/3$, respectively, then this translates into the condition
$|x_2|^4 \leq |x_1|^2 \leq |x_2|^3$. 
 
Consider a particular ``slice'' $\widetilde{F}_{\Lambda_j}$ and the dominant term
$p_{dom_j,C}$.  By Theorem~\ref{p=>pdom}, the dominant term $p_{dom_j,C}$ must be
accurately evaluable everywhere. Hence, in particular, its variety $V(p_{dom_j,C})$
must be allowable. Since the polynomial $p$ vanishes on $V_0$, the dominant term 
$p_{dom_j,C}$ must vanish on $V_0$ as well. So, there are two possibilities: 
$V(p_{dom_j,C})\cap \widetilde{F}_{\Lambda_j}$ either coincides with $V_0$ or
is strictly larger. 

In the first case we apply Lemma~\ref{true_dom} to show that we can evaluate $p$ 
accurately in $\widetilde{F}_{\Lambda_j}$ sufficiently  close to $V_0$, as follows:
Since the polynomial $p_{dom_j,C}$ is accurately evaluable everywhere, for any 
number $\eta_j>0$ there exists an algorithm with output $p_{dom_j,C,comp}$ such that
$|p_{dom_j,C,comp}-p_{dom_j,C}|/|p_{dom_j,C}|<\eta_j$ everywhere. Next, by 
Lemma~\ref{true_dom}, for any $\delta_j>0$ there exists $\varepsilon_j>0$ such 
that~(\ref{dom_inq}) holds. Choose $\eta_j$ and $\delta_j$ so that
$$ \eta_j (1+\delta_j)+\delta_j<\eta.$$ Then we have
\begin{eqnarray*} &&
\left| {p_{dom_j,C,comp} -p  \over p  } \right| \leq {|p_{dom_j,C,comp}-p_{dom_j,C}|
  + |p_{dom_j,C}-p| \over |p|}       \\ &&  ={|p_{dom_j,C,comp}-p_{dom_j,C}| \over 
|p_{dom_j,C}|}\cdot {|p_{dom_j,C}|\over|p|} + {|p_{dom_j,C}-p|\over |p|}   \leq \eta_j 
(1+\delta_j) +\delta_j <\eta  
\end{eqnarray*}
in the $\varepsilon_j$-neighborhood of $V_0$ within the set $\widetilde{F}_{\Lambda_j}$.
Therefore, $p$ can be evaluated by computing $p_{dom_j,C,comp}$ to accuracy $\eta$
in the $\varepsilon_j$-neighborhood of $V_0$ within $\widetilde{F}_{\Lambda_j}$.

In the second case $V(p_{dom_j,C})$ has an irreducible component, say, $W_0$, that is 
strictly larger than $V_0$ and intersects $\widetilde{F}_{\Lambda_j}$ nontrivially. Since
$V(p_{dom_j,C})$ is allowable, it follows from Definition~\ref{chng_var} that
there exists a standard change of variables $C_1$ associated with $V_0$ that maps $W_0$
to a set $\widetilde{x}_1=\cdots=\widetilde{x}_l=0$, $l<k$. Use that change of 
variables and consider the new Newton polytope and dominant terms. Since the polynomial
$p$ is positive in $S\setminus V_0$, there are terms in $p$ that do not contain variables
$\widetilde{x}_1$ through $\widetilde{x}_l$. Therefore each new $p_{dom_r,C_1}$
picks up some of those nonvanishing terms. Hence the (allowable) varieties 
$V(p_{dom_r,C_1})$ have irreducible components strictly smaller than $W_0$ 
(but still containing $V_0$). So, we can subdivide the set $\widetilde{F}_{\Lambda_j}$ 
further using the sets $\widetilde{F}_{\Lambda_r}$ coming from the change of variables 
$C_1$ and the resulting dominant terms will vanish on a set strictly smaller than $W_0$. 
In this fashion, we refine our subdivision repeatedly until we 
obtain a subdivision of the original set $S$ into semialgebraic pieces $(S_j)$
such that the associated dominant terms $p_{dom_j}$ vanish in $S_j$ only on $V_0$.
Applying Lemma~\ref{true_dom} in each such situation, we conclude that $p$ 
can be evaluated accurately sufficiently close to $V_0$ within each piece $S_j$. 

For each $V_0$, we therefore can find a collection $(S_j)$ of semialgebraic 
sets, all determined by polynomial inequalities with integer coefficients, 
and the corresponding numbers $\varepsilon_j$, so that the polynomial 
$p$ can be evaluated with accuracy $\eta$ in each $\varepsilon_j$-neighborhood of $V_0$
within the piece $S_j$. Note that we can assume that each $\varepsilon_j$ is a
reciprocal of an integer, so that testing whether a particular point $x$ is
within $\varepsilon_j$ of $V_0$ within $S_j$ can be done by branching based
on polynomial inequalities with integer coefficients.

The final algorithm will be organized as follows. Given an input $x$, determine by 
branching  whether $x$ is in $S_j$ and within the corresponding $\varepsilon_j$ of a 
component $V_0$.  If that is the case, evaluate $p(x)$ using the algorithm that is 
accurate in $S_j$ in that neighborhood of $V_0$. For $x$ not in any of the neighborhoods,
evaluate $p$ by Horner's rule. Since the polynomial $p$ is strictly positive
off the neighborhoods of the components of its variety, the reasoning of
Section~\ref{sec_PositivePolys} applies, showing that the Horner's rule
algorithm is accurate. If $x$ is on the boundary of a set $S_j$, any applicable
algorithm will do, since the inequalities we use are not strict. Thus the resulting 
algorithm for evaluating $p$ will have accuracy $\eta$ as required. 
\end{proof} 

\subsubsection{Obstacles to full induction}

The reasoning above suggests that there could be an inductive decision procedure that 
would allow us to determine whether or not a given polynomial is accurately evaluable 
by reducing the problem for the original polynomial $p$ to the same problem for
its dominant terms, then their dominant terms, and so forth, going all the way to 
monomials or other polynomials that are easy to analyze. However, this idea
would only work if the dominant terms were somehow ``simpler'' than the original
polynomial itself. This approach would require an induction variable that would 
decrease at each step. 

Two possible choices are the number of variables or the degree of the
polynomial under consideration. Sometimes, however, neither of the two 
goes down, and moreover, the dominant term may even coincide with the
polynomial itself. For example, if $$p(x)=A(x_{[3:n]}) x_1^2 +B(x_{[3:n]})
 x_1 x_2 +C(x_{[3:n]})x_2^2$$ where $A$, $B$, $C$ are nonnegative polynomials 
in  $x_3$ through $x_n$, then the only useful dominant term of $p$ in the neighborhood
of the set $x_1=x_2=0$ is the polynomial $p$ itself. Thus no progress whatsoever 
is made in this situation. 

Another possibility is induction on domains but we do not yet envision 
how to make this idea precise, since we do not know exactly when a given 
polynomial is accurately evaluable on a given domain. 

Further work to establish a full decision procedure is therefore highly desirable. 

\section{``Black-box'' arithmetic} \label{bb}

In this section we prove a necessary condition (for both the real and the complex cases) for a more general type of arithmetic, which allows for ``black-box'' polynomial operations. We describe the type of operations below.

\begin{definition} We call a black-box operation any type of operation that takes a number of inputs (real or complex) $x_1, \ldots, x_k$ and produces an output $q$ such that $q$ is a polynomial in $x_1, \ldots, x_k$. \end{definition}

\begin{example} $q(x_1, x_2, x_3) = x_1 + x_2 x_3$.
\end{example}

\begin{remark} Note that $+, -$, and $\cdot$ are all black-box operations.
\end{remark}

Consider a fixed set of multivariate polynomials $\{ q_j : j\in J\}$ with real or complex inputs (this set may be infinite). In our model under consideration, the arithmetic operations allowed are given by the black-box operations $q_1, \ldots, q_k$, and negation (which will be dealt with by way of dotted edges, as in Section~\ref{class}). With the exception of negation, which is exact, all the others yield a $rnd(op(a_1, \ldots, a_l)) = op(a_1, \ldots, a_l)(1+\delta)$, with $|\delta|< \epsilon$ ($\epsilon$ here is the machine precision). All arithmetic operations have unit cost. We consider the same arithmetical models as in Section~\ref{Sec_Axioms}, with this larger class of operations.

\subsection{Necessity: real and complex}

In order to see how the statement of the necessity Theorem~\ref{conj} changes, we need to introduce a different notion of allowability.

Recall that we denote by $\S$ the space of variables (which may be either $\R^n$ or $\C^n$). From now on we will denote the set $\{1, \ldots, n\}$ by $\k$.

\begin{definition} \label{all_subv}
Let $p(x_1, \ldots, x_n)$ be a multivariate polynomial over $\S$ with variety $V(p)$. Let $\k_Z \subseteq \k$, and let $\k_D, \k_S \subseteq \k \times \k$ . Modify $p$ as follows: impose conditions of the type $Z_i$ for each $i \in \k_Z$, and of type $D_{ij}$, respectively $S_{ij}$, on all pairs of variables in $\k_D$, respectively $\k_S$. Rewrite $p$ subject to those conditions (e.g. set $X_i = 0$ for all $i \in \k_Z$), and denote it by $\tilde{p}$, and denote by $\k_R$ the set of remaining independent variables (use the convention which eliminates the second variable in each pair in $\k_D$ or $\k_S$).

Choose a set $T \subseteq \k_R$, and let
$$ V_{T, \k_Z, \k_D, \k_S} (p) = \cap_{\alpha} V(q_{\alpha})~, $$
where the polynomials $q_{\alpha}$ are the coefficients of the expansion of $\tilde{p}$ in the variables $x_T$:
$$ \tilde{p} (x_1, \ldots, x_k) = \sum_{\alpha} q_{\alpha} x_{T}^{\alpha}~, $$
with $q_{\alpha}$ being polynomials in $x_{\k_R \setminus T}$ only.

Finally, let $\k_N$ be a subset of $\k_R \setminus T$. We negate each variable in $\k_N$, and let $V_{T, \k_Z, \k_D, \k_S, \k_N}(p)$ be the variety obtained from $V_{T, \k_Z, \k_D, \k_S}(p)$, with each variable in $\k_N$ negated. 
\end{definition}

\begin{remark} $V_{\emptyset, \emptyset, \emptyset, \emptyset, \emptyset}(p) = V(p)$. We also note that, if we have a black-box computing $p$, then the set of all polynomials $\tilde{p}$ that can be obtained from $p$ by permuting, repeating, and negating the variables (as in the definition above) is {\sl exactly\/} the set of all polynomials that can be evaluated with {\sl a single\/} rounding error, using that black box. \end{remark}

\begin{definition} For simplicity, we denote a set $(T, \k_Z, \k_D, \k_S, \k_N)$ by 
$\i$, and a set $(T, \k_Z, \k_D, \k_S)$ by $\i_{+}$.
\end{definition}

\vspace{.25cm}

\begin{example} Let $p(x,y,z) = x+y \cdot z$ (the fused multiply-add). We record below all possibilities for $\i = (T, \k_Z, \k_D, \k_S, \k_N)$, together with the obtained subvariety $V_{\i}(p)$.

Without loss of generality, assume that we have eliminated all redundant or complicated conditions, like $(x,y) \in \k_D$ and $(x, y) \in \k_S$ (which immediately leads to $x=y=0$, that is, $x, y, \in \k_Z$. We assume thus that all variables not present in $\k_Z$ cannot be deduced to be $0$ from conditions imposed by $\k_D$ or/and $\k_S$.

We obtain that all possibilities for $V_{\i}(p)$ are, up to a permutation of the variables,
\begin{enumerate} \item[{\mathversion{bold} $\diamond$}] $\{x= 0\}, ~\{x = 1\},~ \{x = -1\}$, 
\item[{\mathversion{bold} $\diamond$}] $\{x=0\} \cup \{x=1\}$, $\{x=0\} \cup \{x = -1\}$, 
\item[{\mathversion{bold} $\diamond$}] $\{x = 0\} \cup \{y = 0\}$, $\{x = 0\} \cup \{y = 1\}$, $\{x=0\} \cup \{y = -1\}$,
\item[{\mathversion{bold} $\diamond$}] $\{x = -y^2\}$, $\{x = y^2\}$, $\{x - y \cdot z = 0\}$, and $\{x+y \cdot z=0\}$. \end{enumerate} \end{example}

\begin{definition} \label{pmm} We define $q_{-2}(x_1, x_2) = x_1 x_2$, $q_{-1}(x_1, x_2) = x_1+x_2$, and $q_0(x_1, x_2) = x_1 - x_2$. \end{definition}

\begin{remark} The sets 
\begin{eqnarray} \label{q-unu}
& 1.&  Z_i =\{x~:~ x_i ~=~ 0\}~, \\
\label{q-doi}
& 2.&  S_{ij} = \{x~:~ x_i+x_j ~=~ 0\}~, \\
\label{q-trei}
& 3.&  D_{ij} = \{x~:~ x_i - x_j ~=~ 0\}~, 
\end{eqnarray}
and unions thereof, describe all non-trivial (neither $\emptyset$ nor $\S$) sets of type $V_{\i}$, for $q_{-2}, ~q_{-1}$, and $q_0$.
\end{remark}

\vspace{.25cm}

We will assume from now on that the black-box operations $q_{-2}, q_{-1}, q_0$ defined in~\ref{pmm}, and some arbitrary extra operations $q_j$, with $j \in J$ ($J$ may be infinite) are given and fixed. 

\vspace{.25cm}

\begin{definition} \label{q_allow} 
We call any set $V_{\i}(q_j)$ with $\i= (T, \k_Z, \k_D, \k_S, \k_N)$ as defined above and $q_j$ a black-box operation
{\sl basic $q$-allowable}.

We call any set $R$ {\sl irreducible $q$-allowable} if it is an irreducible component of a (finite) intersection of basic $q$-allowable sets, i.e., when $R$ is irreducible and 
\[ 
R \subseteq \cap_{l} ~Q_l~,
\]
where each $Q_l$ is a basic $q$-allowable set.

We call any set $Q$ {\sl $q$-allowable} if it is a (finite) union of irreducible $q$-allowable sets, i.e.\[
Q = \cup_{j} R_j~,
\]
where each $R_j$ is an irreducible $q$-allowable set. 

Any set $R$ which is not $q$-allowable we call $q$-unallowable. 
\end{definition}

\begin{remark} Note that the above definition of $q$-allowability is closed under taking union, intersection, and irreducible components. This parallels the definition of allowability for the classical arithmetic case -- in the classical case, every allowable set was already irreducible (being an intersection of hyperplanes). \end{remark}




Once again, we need to build the setup to state and prove our new necessity condition. To do this, we will modify the statements of the definitions and the statements and proofs of the lemmas from Section~\ref{class_nec}. Since most proofs just follow in the footsteps of those from Section~\ref{class_nec}, instead of repeating them, we will only point out the places where they differ and show how we modified them to work in the new context.

\begin{definition}
Given a polynomial $p$ with $q$-unallowable variety $V(p)$, consider 
all sets $W$ that are $q$-allowable (as in Definition~{\sl\ref{q_allow}}), and subtract from $V(p)$ those $W$ for which $W \subset V(p)$. 
We call the remaining subset of the variety {\em points in general position\/} and denote it by $\G(p)$.
\end{definition}

\begin{remark}
Since $V(p)$ is $q$-unallowable, $\G(p)$ is non-empty.
\end{remark}

\begin{definition}
Given $x\in \S$, define the set $\qallow(x)$ as the intersection of
all basic $q$-allowable sets going through $x$:
$$ \qallow(x) \eqbd \cap_{j \in J \cup \{-2, -1, 0\}}  \left( \cap_{I~:~ x\in V_{\i}(q_j)} ~~V_{\i}(q_j) \right) , $$
for all possible choices of $T, \k_Z, \k_D, \k_S, \k_N$.

The intersection in parentheses is $\S$  whenever $x \notin V_{\i}(q_j)$ for all possible $\i$.
\end{definition}

\begin{remark}
When $x \in \G(p)$, $\qallow(x) \not\subseteq \G(p)$.
\end{remark}

We can now state our necessity condition.

\begin{theorem} \label{gen_result} Given the black-box operations $q_{-2}, q_{-1}, q_0$, and $\{q_j: j \in J\}$, and the model of arithmetic described above, let $p$ be a polynomial defined over a domain $\mathcal{D} \subset \S$. Let $\mathcal{G}(p)$ be the set of points in general position on the variety $V(p)$. 
If there exists $x \in \mathcal{D} \cap \mathcal{G}(p)$ such that $\qallow(x) \cap \Int(\mathcal{D}) \neq \emptyset$, then $p$ is not accurately evaluable on $\mathcal{D}$. \end{theorem}

We proceed to the construction of the elements of the proof of Theorem~\ref{gen_result}. The algorithm will once again be represented by a DAG with input nodes, branching nodes, and output nodes. As in Section~\ref{class}, for simplicity in dealing with negation (since negation is exact), we will work with \emph{solid\/} edges, which convey a value unchanged, and \emph{dotted\/} edges, which indicate that negation of the conveyed quantity has occurred.

{From now on, unless specified, we will consider only non-branching algorithms}.


We will continue to use the definition of a Zariski set (Definition~\ref{zariski}) on a hypercube in $\delta$-space, and we work with the same definition of a non-trivial computational node (recalled below).

\begin{definition} \label{q-non-trivial} For a given $x\in \S$, we say that a computational node $N$ is {\em of $q$-non-trivial type\/} if its output is a nonconstant polynomial of $\delta$.
 \end{definition}



Recall the notation $H_{\epsilon}$ from Definition~\ref{eps-hypercube}.

The equivalent of Proposition~\ref{unu_unu} becomes the following.

\begin{proposition} \label{q-unu_unu} Given any algorithm, any $\epsilon>0$, and a point $x$ in $\mathcal{G}(p)$, there exists a  Zariski open set $\Delta \subset H_{\epsilon}$ such that no $q$-non-trivial computational node has a zero output on the input $x$ for all $\delta \in \Delta$.
\end{proposition}
\begin{proof}
The proof of Proposition~\ref{q-unu_unu} follows the same path as that of Proposition~\ref{unu_unu}. To each $q$-non-trivial node corresponds a Zariski open set in $\delta$-space; there is a finite number of them, and their intersection provides us with the Zariski open set we are looking for.
\end{proof}

We will now state and sketch the proof of the equivalent of Lemma~\ref{prima}.

\begin{lemma} \label{q-prima}
For a given algorithm, $x\in \mathcal{G}(p)$, and $\epsilon>0$, exactly one of the following holds:
\begin{enumerate} \item there exists a  Zariski open set $\Delta \in H_{\epsilon}$ such that the value $p_{comp}(x, \delta)$ computed by the algorithm is not zero when the algorithm 
is run with source input $x$ and errors $\delta \in \Delta$;
\item $p_{comp}(y, \delta)=0$ for all $y\in \qallow(x)$ and all $\delta \in H_{\epsilon}$.
\end{enumerate}
\end{lemma}

\noindent \textit{Proof of Lemma~\ref{q-prima}.} 
Give $x \in \G(p)$, choose $\delta$ from the Zariski open set $\Delta$ whose existence is given by Proposition~\ref{q-unu_unu}. Either the output node is of $q$-non-trivial type (in which case $p_{comp}(x, \cdot) \neq 0$ and we are done) or the output is a nonzero constant polynomial in $\delta$ (and again we are done) or the output is the zero polynomial in $\delta$. In this latter case, we trace back all the zeros again, as in the proof of Lemma~\ref{prima}, and get a set of paths of nodes that produced all $0$. 

Let us start from the last nodes on these paths and work our way up, level after level. The last node on such a path is either a source, or a node with all inputs from sources, or a node with at least one input which is not a source and not $0$ (and hence a polynomial in $\delta$). In the former two cases, we have traced back the zeros to basic $q$-allowable conditions. We will show that this is also true for the latter case. 



\begin{lemma} \label{case1} If the last zero occurs at a node which computes a black-box operation $q_j$ and which has some source inputs and some (nontrivial) polynomial inputs, then the sources lie in some $V_{\i}(q_j)$ (and this constraint causes the zero output of this node).
\end{lemma}

\noindent \textit{Proof of Lemma~\ref{case1}.} Label the node $\tilde{N}$ (assume it corresponds to the black-box $q_j$) with output $0$; then some inputs are sources, and some are polynomials of $\delta$ (since this was the last node on the path, it has no zero inputs). By the choice of $\delta$, it follows that the output has to be the $0$ polynomial in $\delta$. 

Some of the non-source inputs to $\tilde{N}$ might come from the same nodes; assume you have a total of $l$ distinct nodes which input to $\tilde{N}$ (nonconstant) polynomials of $\delta$, and let $I_1(\delta), I_2(\delta), \ldots$, $ I_l(\delta)$ denote these inputs. We will need the following lemma.

\begin{lemma} \label{alg}Since the DAG is acyclic, $I_1(\delta), I_2(\delta), \ldots, I_l(\delta)$ are  {\sl algebraically independent\/} polynomials in $\delta$. \end{lemma}

\noindent \textit{Proof of Lemma~\ref{alg}.} Suppose there is some polynomial dependence 
among $I_1(\delta), \ldots$, $I_l(\delta)$. Let $N_1, \ldots$, $N_l$ be the nodes
which have computed $I_1(\delta), \ldots, I_l(\delta)$, and
let $\delta_1, \delta_2, \ldots, \delta_l$ be the specific $\delta$s at
these nodes. Let $D_1, \ldots, D_l$ be the set of $\delta$s present
(non-trivially) in the output of each node, e.g., $\delta_i \in D_i$. At
least one $\delta_i$ is not present in $\cup_{i \neq j} D_j$; otherwise
we get a cycle. But then $I_i(\delta)$ is algebraically independent from
the other inputs, i.e., there is some dependence among the
inputs $I_j(\delta)$ with $j \neq i$. We use induction on the number of
remaining inputs, and exclude one input at a time, until we're left
with a contradiction. This proves Lemma~\ref{alg}. \qed

Replace each $I_i(\delta)$ (or $-I_i(\delta)$) in $Y$ by the same \emph{dummy\/} variable $z_i$ (respectively $-z_i$), for each $i \in \{1, \ldots, l\}$. The variables $z_i$ are algebraically independent by Lemma~\ref{alg}. 

Denote by $z$ the new vector of inputs to $\tilde{N}$ (both values \emph{and\/} variables). The value $q_j(z)=0$, regardless of the $z_i$s (since it was $0$ regardless of the $\delta$s and the $z_i$ are algebraically independent variables). It follows that the constraints which place $z$ on the variety of $q_j$ are twofold: they come from constraints of the type $D_{ij}$ and $S_{ij}$ which describe the places where we inputted the values $z_i$, and they come from imposing conditions on the \emph{other\/} inputs, which are \emph{sources\/.} Thus the constraints on $z$ are of the form $V_{T, \emptyset, \k_D, \k_S, \k_N}(q_j) \bdeq V_{\i}$. This concludes the proof of Lemma~\ref{case1}. \qed

Now that we have shown that the last marked vertices all provide basic $q$-allowable conditions, we proceed by induction: we look at a ``next'' marked vertex (here ``next'' means that all its marked ancestors have been examined already). It has some zero inputs, some source inputs, some of the inputs satisfy constraints of type $D_{ij}$ and $S_{ij}$, and some of the inputs are polynomial. From here on we proceed as in Lemma~\ref{case1}, and obtain a set of new constraints to be imposed on the sources, of the type $V_{\i}(q_j)$, which we will intersect with the rest of the constraints obtained so far. 

At the end of the examinations, we have found a set of basic $q$-allowable constraints which the sources must satisfy, i.e., a list of basic $q$-allowable sets with the property that the sources lie in their intersection; the fact that the sources satisfy these constraints is responsible for the zero output at the end of the computation. 

It is not hard to see that in this case, once again, it follows that for all $y$ in $\qallow(x)$ and any $\delta \in \Delta$, the output is $0$ (just as in Lemma~\ref{prima}). Thus we have proved Lemma~\ref{q-prima}. \qed

From Lemma~\ref{q-prima} we obtain the following corollary.

\begin{corollary} \label{q-partsial} For any algorithm, for any $\epsilon>0$, and any $x\in \mathcal{G}(p)$, exactly one of the following holds: the relative error of computation, 
$|p_{comp}-p|/|p|$, is either infinity at $x$ for all $\delta$ in a  Zariski open set or $1$ at all
points $y\in (\qallow(x)\setminus V(p))$ and all $\delta \in H_{\epsilon}$.
\end{corollary}

We will now consider algorithms with or without branches.

\begin{theorem} \label{q-finala} Given a (branching or non-branching) algorithm with output function $p_{comp}(\cdot)$, $x \in \mathcal{G}(p)$, and $\epsilon>0$, then one of the following is true:
\begin{enumerate} 
\item there exists a set $\Delta_1 \in H_{\epsilon}$ of positive measure such that $p_{comp}(x, \delta)$ is nonzero whenever the algorithm is run with errors $\delta \in \Delta_1$, or
\item there exists a set $\Delta_2 \in H_{\epsilon}$ of positive measure such that for every $\delta \in \Delta_2$, there exists a neighborhood $N_{\delta}(x)$ of $x$ such that 
for every $y \in N_{\delta}(x) \cap \left ( \qallow(x) \setminus V(p) \right )$, $p_{comp}(y, \delta) = 0$ when the algorithm is run with errors $\delta$.
\end{enumerate}
\end{theorem}

\begin{remark}
Just as before, this implies that, on a set of positive measure in $H_{\epsilon}$, the relative accuracy 
of any given algorithm is either $\infty$ or $1$.  
\end{remark}

\begin{proof} The proof is essentially the same as in Theorem~\ref{finala}; 
the only thing that needs to be examined is the existence in $\qallow(x) \setminus V(p)$ of an infinite sequence $\{y_n\}$ with $y_n \rightarrow x$. 

We will make use of the following basic result in the theory of algebraic varieties, which can for example be found as Theorem 1 in~\cite[Section 6.1]{Shafarevich}.

\begin{result} \label{irred_cont}
If $X$ and $Y$ are polynomial varieties such that $X \subseteq Y$ , then $dim(X) \leq dim(Y)$. If $Y$ is irreducible and $X \subseteq Y$ is a (closed) subvariety with $dim(X) = dim(Y)$, then $X = Y$. \end{result}


We write $\qallow(x)$ as a union of irreducible $q$-allowable components. By the way we defined $\mathcal{G}(p)$, it follows that none of these components is included in $V(p)$; by Result~\ref{irred_cont}, it follows that the intersection of any irreducible $q$-allowable component $P$ of $\qallow(x)$ with $V(p)$ has a smaller dimension than $P$. 

Choose the (unique) irreducible component $P$ that contains $x$; this component must have dimension at least $1$ (since if it contained only $x$, the set $\{x\}$ would be $q$-allowable, and hence we would have extracted it from $V(p)$, which is a contradiction with the fact that $x \in \mathcal{G}(p)$). Since $P \setminus V(p)$ has a smaller dimension than $P$, there \emph{must\/} be some infinite sequence $\{y_n\}$ in $P\setminus V(p)$, i.e. in $\qallow(x) \setminus V(p)$, such that $y_n \rightarrow x$.

The rest of the argument goes through just as in Theorem~\ref{finala}.
\end{proof}

Finally, as in Section~\ref{class_nec}, we have a corollary.

\begin{corollary} \label{q_in_sfirsit} Let $p$ be a polynomial over $\S$ with unallowable variety $V(p)$. Given any algorithm with output function $p_{comp}(\cdot)$, a point $x \in \mathcal{G}(p)$, $\epsilon>0$, and $\eta<1$, there exists a set $\Delta_x$ of positive measure {\sl arbitrarily close} to $x$ and a set $\Delta \in H_{\epsilon}$ of positive measure, such that $|p_{comp} - p|/|p|$ is strictly larger than $\eta$ when computed at a point $y \in \Delta_x$ with errors $\delta \in \Delta$.
\end{corollary}

The proof is based on the topology of $\S$, and is identical to the proof of Corollary~\ref{in_sfirsit}; we choose not to repeat it. 

\vspace{.3cm}

\noindent \textit{Proof of Theorem~\sl{\ref{gen_result}}.}  Follows immediately from Theorem~\ref{q-finala} and Corollary~\ref{q_in_sfirsit}. \qed

\subsection{Sufficiency: the complex case} \label{suf_c_bb}

In this section we obtain a sufficiency condition for the accurate evaluability of a complex polynomial, given a black-box arithmetic with operations $q_{-2}, q_{-1}, q_0$ and $\{q_j | j \in J\}$ ($J$ may be an infinite set). 

\vspace{.25cm}

Throughout this section, we assume our black-box operations include $q^{c}$, which consists of multiplication by a complex constant: $q^c(x) = c\cdot x$. Note that this operation is natural, and that most computers perform it with relative accuracy.

\vspace{.25cm}

We believe that the sufficiency condition we obtain here is sub-optimal in 
general, but it subsumes the sufficiency condition we found for the basic complex case with classical arithmetic $\{+, -, \cdot\}$. 

We assume that the black-box polynomials defining the operations $q_j$ with $j \in J$ are \emph{irreducible.}

\begin{lemma} \label{irred_af}The varieties $V_{\i}(q_j)$ are irreducible for any $j \in J$ and any $\i$ as in Definition~\ref{all_subv} if and only if all $q_j$, $j \in J$ are affine polynomials. \end{lemma}

\begin{proof} If $q_j$ is an affine polynomial then any $V_{\i}(q_j)$ is also affine, hence irreducible over $\C^n$. Conversely, if $q_j$ is not an affine polynomial, then by inputting a single value $x$ for all the variables, we obtain a one-variable polynomial of degree at least $2$, which is necessarily reducible over $\C^n$. \end{proof}

We state here the best sufficiency condition for the accurate evaluability of a polynomial we were able to find in the general case, and a necessary and sufficient condition for the all-affine black-box operations case.

\begin{theorem}[General case] \label{q-suff-c1} Given a polynomial $p~:~\C^n~\rightarrow~\C$ with $V(p)$ a finite union of intersections of hyperplanes $Z_{i}, ~S_{i j},~D_{ij}$, and varieties $V(q_j)$, for $j \in J$, then $p$ is accurately evaluable. \end{theorem}

\begin{theorem}[Affine case] \label{q-suff-c2} If all black-box operations $q_j$, $j \in J$ are affine, then a polynomial $p~:~\C^n~\rightarrow~\C$ is accurately evaluable iff $V(p)$ is a union of intersections of hyperplanes $Z_{i}, ~S_{i j},~D_{ij}$, and varieties $V_{\i}(q_j)$, for $j \in J$ and $\i$ as in Definition~\ref{all_subv}.
\end{theorem}

We will begin by proving Theorem~\ref{q-suff-c1}. We will once again make use of Theorem \ref{dimensions} and of Theorem~\ref{irred_cont}.

\begin{lemma} \label{key_q} If $V(p)$ is as in Theorem~\ref{q-suff-c1}, 
then $V(p)$ is a~\emph{simple\/} finite union of
hyperplanes  $Z_{i}, ~S_{i j},~D_{ij}$ and varieties $V(q_j)$ (\emph{with no intersections\/}). 
\end{lemma}

\begin{proof} Indeed, if that were not the case, then some irreducible $q$-allowable component $P$ of $V(p)$ would be an
intersection of two or more sets described in Theorem \ref{q-suff-c1}. If $P$ were contained in the intersection of two or more (distinct) hyperplanes, its dimension would be 
smaller than $n-2$, and we would
get a contradiction to Theorem~\ref{dimensions}. 

Suppose now that $P$ was contained in the intersection 
of a $V(q_j)$ with some
other variety or hyperplane. All such varieties, by Theorem~\ref{dimensions}, must have dimension $n-1$, and since all such varieties and hyperplanes are irreducible, by Result~\ref{irred_cont}, their intersection must have dimension strictly smaller than $n-2$. Contradiction; we have thus proved that the variety $V(p)$ is a \emph{simple\/} union of hyperplanes $Z_i,~D_{ij},~S_{ij}$, and varieties $V(q_j)$.
\end{proof}

\begin{corollary} \label{q-factors} If $p: \C^n \to \C$ is a 
polynomial whose variety $V(p)$ is $q$-allowable, then it is a product 
$p=c \prod_j p_j$, where each $p_j$ is a power of $x_i$, $(x_i-x_j)$, $(x_i +x_j)$, or $q_j$, and $c$ is a complex constant.
\end{corollary}

\begin{proof} By Lemma~\ref{key_q}, the variety $V(p)$ is
a union of basic $q$-allowable hyperplanes and varieties $V(q_j)$. 

Choose an irreducible $q$-allowable set in the union. If this set is a hyperplane, then by following the same argument as in Corollary~\ref{factors}, we obtain that $p$ factors into some $\widetilde{p}$ and some power of either $x_i$, $(x_i-x_j)$, or $(x_i+x_j)$.

Suppose now that the irreducible $q$-allowable set were a variety $V(q_j)$; since $p$ is $0$ whenever $q_j$ is $0$ and $q_j$ is irreducible, it follows that $q_j$ divides $p$. We factor then $p$ into the largest power of $q_j$ which divides it, and some other polynomial $\widetilde{p}$. 

In either of the two cases, we proceed by factoring $\widetilde{p}$ in the same fashion, until we encounter a polynomial $\widetilde{p}$ of degree $0$. That polynomial is the constant $c$.
\end{proof}

\noindent \textit{Proof of Theorem~\sl{\ref{q-suff-c1}}.}   
By Corollary~\ref{factors},  $ p=c\prod_j p_j $, with each $p_j$ a power of $x_k$, $(x_k \pm x_l)$, or $q_l$. 

Since each of the factors is accurately evaluable, the algorithm that forms their product evaluates $p$ accurately. Multiplication by $c$ (corresponding to the black-box $q^c$) is also accurate, hence $p$ is accurately evaluable. 
\qed

\vspace{.25cm}

The proof of Theorem~\ref{q-suff-c2} follows the path described above; we sketch it here.

\vspace{.25cm}

\noindent \textit{Proof of Theorem~{\sl\ref{q-suff-c2}}.}
The key fact in obtaining this condition is the irreducibility of all sets $V_{\i}(q_j)$, which is guaranteed by Lemma~\ref{irred_af}, together with the result of Lemma~\ref{key_q}. Once again, we can write the polynomial as a product of powers of $x_i$, $(x_i \pm x_j)$, or $V_{\i}(q_j)$, times a constant; this takes care of the sufficiency part, while the necessity follows from Theorem~\ref{gen_result}. \qed

\begin{remark}
Note that Theorem~\ref{q-suff-c2} is a more general necessary and sufficient condition than Theorem~\ref{sufficiency_c}, which only considered having $q_{-2}, q_{-1}$, and $q_0$ as operations, and restricted the polynomials to have integer coefficients (thus eliminating the need for $q^c$).
 \end{remark}

\section{Accurate Linear Algebra in Rounded Arithmetic}
\label{sec_LinearAlgebra}

Now we describe implications of our results to the question
of whether we can accurately do
numerical linear algebra on structured matrices.
By a {\em structured matrix\/} we mean a family of 
$n$-by-$n$ matrices $M$ whose entries $M_{ij}(x)$ 
are simple polynomial or rational 
functions of parameters $x$.
Typically there are only $O(n)$ parameters, 
and the polynomials $M_{ij}(x)$ are closely related
(for otherwise little can be said).
Typical examples include 
Cauchy matrices 
($M_{ij}(x,y) = 1/(x_i + y_{j})$), 
Vandermonde matrices 
($M_{ij}(x) = x_i^{j-1}$),
generalized Vandermonde matrices 
($M_{ij}(x) = x_i^{j-1+\lambda_j}$, where the $\lambda_j$ are 
a nondecreasing sequence of nonnegative integers),
Toeplitz matrices  
($M_{ij}(x) = x_{i-j}$),
totally positive matrices (where $M$ is expressed as a product of
simple nonnegative bidiagonal matrices arising from its
Neville factorization),
acyclic matrices,
suitably discretized elliptic partial differential
operations, 
and so on
\cite{demmelICM02,demmelkoevICIAM03,demmelkoev99,DGESVD,demmel99,mmatsvd,demmelkoevPolyVandermonde,demmelkoevTPGenVandermonde,koevAccTP04,koevAccEig}.

It has been recently shown that all the matrices on the
above list 
(except Toeplitz and non-totally-positive generalized Vandermonde matrices) 
admit accurate algorithms in rounded arithmetic for many
or all of the problems of numerical linear algebra:
\begin{itemize}
\item computing the determinant
\item computing all the minors
\item computing the inverse
\item computing the triangular factorization from Gaussian elimination,
with various kinds of pivoting
\item computing eigenvalues
\item computing singular values
\end{itemize}

We have gathered together these results in Table~\ref{table1}.

\newcommand{\ts}{\rm}
\newcommand{\IGNORE}[1]{}


{\small 
\begin{center}

\begin{table}[ht]
\caption{General Structured Matrices \label{table1} }

\smallskip

\begin{tabular}{|ll|c|c|c|c|c|c|}
\hline

\ts
&     &     &   &\ts Any   &\ts    &\ts    &\ts Sym   \\
\multicolumn{2}{|l|}{Type of matrix}  &\ts $\det A$&\ts $A^{-1}$&\ts  
minor &\ts  LDU & \ts  
SVD & 
\ts EVD 
\\

\hline
\multicolumn{2}{|l|}{Acyclic} 
 & $n$ & $n^2$ & $n$ & $\le n^2$ & $n^3$ & N/A 
\\
\multicolumn{2}{|l|}{(bidiagonal and other)} 
& \cite{DGESVD} & \cite{DGESVD} & \cite{DGESVD} & \cite{DGESVD}
& \cite{DGESVD} &   
\\

\hline
\multicolumn{2}{|l|}{Total Sign Compound } &
$n$ & $n^3$ & $n$ & $n^4$ & $n^4$ & $n^4$  
\\
(TSC)&&\cite{DGESVD} & \cite{DGESVD} & \cite{DGESVD}
&\cite{DGESVD} & \cite{DGESVD} & \cite{dopicomoleramoro03} 
\\
\hline
\multicolumn{2}{|l|}{Diagonally Scaled Totally}
& $n^3$ & $n^5$ & $n^3$ & $n^3$ & $n^3$ & $n^3$ 
\\
\multicolumn{2}{|l|}{Unimodular (DSTU)} &\cite{DGESVD}
&&\cite{DGESVD}&\cite{DGESVD} & \cite{DGESVD} &\cite{dopicomoleramoro03}
\\

\hline

 \multicolumn{2}{|l|}{Weakly diagonally}
 &\ts  $n^3$
               &\ts  $n^3$
               &\ts  {No}
               &\ts  \textcolor{black}{$ {n^3}$}
               &\ts  \textcolor{black}{$ {n^3}$}
               &\ts  \textcolor{black}{$ {n^3}$}
\\

\multicolumn{2}{|l|}{dominant
M-matrix}     &\ts  \cite{ocinneide96}
               &\ts  \cite{alfaxueye,alfaxueye2}
               &\ts  \cite{demmelkoev99}
               &\ts  \cite{mmatsvd}
               &\ts  \cite{mmatsvd}
&\cite{dopicomoleramoro03}
\\

\hline
\ts
&Cauchy  
&\ts \textcolor{black}{$n^2$}     
&\ts \textcolor{black}{$n^2$}
&\ts \textcolor{black}{$n^2$}
& $\le n^3$
&\ts \textcolor{black}{$n^3$}   
&\ts  \textcolor{black}{$n^3$}
\\
Displace-&&
&
&
&\ts \cite{demmel99}
&\ts \cite{demmel99}
&\cite{dopicomoleramoro03}
\\

\cline{2-8}

\ts
ment & Vandermonde &\ts  $n^2$
             &\ts  \textcolor{black}{{\rm No}}
             &\ts  \textcolor{black}{{\rm No}}
             &\ts  \textcolor{black}{{\rm No}}
& $n^3$ & $n^3$ 
\\
Rank One &
             &
             &\ts \cite{demmelkoev99}
             &\ts \cite{demmelkoev99}
             &\ts \cite{demmelkoev99}
             &\ts \cite{demmel99,demmelkoevPolyVandermonde}
&\cite{dopicomoleramoro03}
\\

\cline{2-8}

\ts
 & Polynomial  
&\ts  $n^2$ 
&\ts  \textbf{No}
&\ts  \textbf{No}
&\ts  \textbf{No}
& $*$ & $*$ 
\\
\ts & Vandermonde
&\ts \cite{higham90e}
&\ts Section \ref{vander} 
&\ts Section \ref{vander}
&\ts Section \ref{vander}
&\ts \cite{demmelkoevPolyVandermonde} & \cite{dopicomoleramoro03}
\\ 

\cline{2-8}



\hline

\end{tabular}
\end{table} 
\end{center}
} 


The proliferation of these accurate algorithms for some but not all
matrix structures motivates us to ask for which structures they exist.

To convert this to a question about polynomials, we begin by noting
that being able to compute the determinant accurately is a {\em necessary\/}
condition for most of the above computations. For example, if the diagonal
entries of a triangular factorization of $A$, or its eigenvalues, are
computable with small relative error, then so is their product, the
determinant. 

It is also true that being able to compute all the
minors accurately is a {\em sufficient\/} condition for many of the above
computations. For the inverse and triangular factorization, this
follows from Cramer's rule and Sylvester's theorem, resp., and for
the singular values an algorithm is described in~\cite{DGESVD,demmelkoev99}.

Thus, if the determinants $p_n(x) = {\rm det} M^{n \times n}(x)$
of a class of $n$-by-$n$ structured matrices $M$ do not satisfy the 
necessary conditions described in Theorem~\ref{gen_result} for {\em any\/} 
enumerable set of black-box operations (perhaps with other properties, like bounded degree), then we can conclude that 
accurate algorithms of the sort described in the above citations
are impossible.

In particular, to satisfy these necessary conditions would require that
the varieties $V(p_n)$ be allowable (or $q$-allowable).
For example, if $V$ is a Vandermonde matrix, then 
${\rm det}(V) = \prod_{i<j} (x_i-x_j)$ satisfies this condition, using
only subtraction and multiplication.

The following theorem states a condition which guarantees the
impossibility of an algorithm using {\em any\/} enumerable set of
black-box operations of bounded degree:

\begin{theorem}
\label{Thm_StructuredMatrixImpossibility}
Let $M(x)$ be an $n$-by-$n$ structured complex matrix with determinant
$p_n(x)$ as described above. Suppose $p_n(x)$ has an irreducible
factor $\hat{p}_n(x)$ whose degree goes to infinity as $n$ goes
to infinity. Then for any enumerable set of black-box arithmetic
operations of bounded degree, for sufficiently large $n$ it is impossible to 
accurately evaluate $p_n(x)$ over the complex numbers.
\end{theorem}

\begin{proof}
Let $q_1,...,q_m$ be any finite set of black-box operations.
To obtain a contradiction, suppose the complex variety
$V(p_n)$ satisfies the necessary conditions of Theorem~\ref{gen_result},
i.e., that $V(p_n)$ is allowable. 
This means that $V(p_n)$, which
includes the hypersurface $V(\hat{p}_n)$ as an irreducible component, 
can be written as the union of irreducible 
$q$-allowable sets (by Def.~\ref{q_allow}).
This means that $V(\hat{p}_n)$ must
itself be equal to an irreducible 
$q$-allowable set (a hypersurface), since representations
as unions of irreducible sets are unique.
The irreducible $q$-allowable sets of codimension 1 are
defined by single irreducible polynomials, which are
in turn derived by the process of setting variables
equal to one another, to one another's negation, or zero
(as described in Defs.~\ref{all_subv} and~\ref{q_allow}),
and so have bounded degree.
This contradicts the unboundedness of the degree of $V(\hat{p}_n)$.
\end{proof}

In the next theorem we apply this result to the set of 
complex Toeplitz matrices. We use the following notation.
Let $T$ be an $n$-by-$n$ Toeplitz matrix, with
$x_j$ on the $j$-th diagonal, so $x_0$ is on the main diagonal,
$x_{n-1}$ is in the top right corner, and $x_{1-n}$ is in the bottom left
corner. 

\begin{theorem}
\label{Thm_Toeplitz_C}
The determinant of a Toeplitz matrix $T$ is irreducible over any field.
\end{theorem}

\begin{corollary}
The determinants of the set of complex Toeplitz matrices 
cannot be evaluated accurately using any enumerable set of bounded-degree 
black-box operations.
\end{corollary}

\vspace{.25cm}

\noindent \textit{Proof of Theorem~\ref{Thm_Toeplitz_C}}.
We use induction on $n$.
We note that ${\rm det}~T$ depends on every variable $x_j$, because
${\rm det}~T$ includes the monomials $\pm x_j^{n-j}x_{j-n}^j$
for $j > 0$,
as well as $x_0^n$, and these monomials
contain the maximum powers of $x_j$ and $x_{j-n}$ appearing 
in the determinant.
Now $x_{n-1}$ appears exactly once in $T$, 
so ${\rm det}~T$ must be an affine function
of $x_{n-1}$, say ${\rm det}~T = x_{n-1} \cdot p_{1n} + p_{2n}$. 
By expanding ${\rm det}~T$
along the first row or column, 
we see that $p_{1n}$ is itself the determinant
of a Toeplitz matrix with diagonals $x_{1-n},...,x_{n-3}$, 
and $p_{2n}$ depends on $x_{1-n},...,x_{n-2}$ but not $x_{n-1}$.
If ${\rm det}~T = x_{n-1} \cdot p_{1n} + p_{2n}$ were reducible, 
its factorization would have to look like $x_{n-1} \cdot p_{1n} + p_{2n} = 
(x_{n-1} \cdot p_{3n} + p_{4n})p_{5n}$, 
where all the subscripted $p$ polynomials are independent of $x_{n-1}$,
implying either that
$p_{1n} = p_{3n} p_{5n}$ were reducible, a contradiction by our
induction hypothesis, or
$p_{3n}=1$ and so $p_{1n} | p_{2n}$. 
Now we can write 
$p_{2n} = x_{n-2}^2 q_{1n} + x_{n-2} q_{2n} + q_{3n}$
where $q_{1n} \neq 0$,
since ${\rm det}~T$ includes the monomial
$\pm x_{n-2}^2 x_{-2}^{n-2}$ and no higher powers of $x_{n-2}$.
Furthermore $q_{1n}$ is independent of $x_{n-1}$, $x_{n-2}$ and $x_{n-3}$,
and $q_{2n}$ and $q_{3n}$ are independent of $x_{n-1}$ and $x_{n-2}$.
Since $p_{1n}$ is independent of $x_{n-2}$,
the only way we could have $p_{1n} | p_{2n}$ is to have
$p_{1n} | q_{1n}$, 
$p_{1n} | q_{2n}$, and
$p_{1n} | q_{3n}$. But since $p_{1n}$ depends on $x_{n-3}$
and $q_{1n}$ is independent of $x_{n-3}$, this is a contradiction.
So the determinant of a Toeplitz matrix must be irreducible. \qed

\vspace{.25cm}

In the real case, irreducibility of $p_n$ 
is not enough to conclude that $p_n$ cannot be
evaluated accurately, because $V_{\R}(p_n)$ may
still be allowable (and even vanish). 
So we consider another necessary condition for
allowability: Since all black-boxes have a finite
number of arguments, their associated codimension-1
irreducible components
must have the property that
whether $x \in V_{\i}(q_j)$ depends on only 
a finite number of components of $x$.
Thus to prove that the hypersurface $V_{\R}(p_n)$ is not allowable,
it suffices to find at least one regular point $x^*$ in $V_{\R}(p_n)$
such that the tangent hyperplane at $x^*$ 
is not parallel to
sufficiently many coordinate directions, i.e., membership in
$V_{\R}(p_n)$ depends on more variables than any $V_{\i}(q_j)$.
This is easy to do for real Toeplitz matrices.

\begin{theorem}
\label{thm_Toeplitz_R}
Let $V$ be the variety of the determinant of real singular Toeplitz matrices.
Then $V$ has codimension 1, and at almost all regular points,
its tangent hyperplane is parallel to no coordinate directions.
\end{theorem}

\begin{corollary}
The determinants of the set of real Toeplitz matrices 
cannot be evaluated accurately using any enumerable set of bounded-degree 
black-box operations.
\end{corollary}

\vspace{.25cm}

\noindent \textit{Proof of Theorem~\ref{thm_Toeplitz_R}}.
Let $Toep(i,j)$ denote the Toeplitz matrix
with diagonal entries $x_i$ through $x_j$;
thus $Toep(i,j)$ has dimension $(j-i)/2+1$.
Let $U$ be the Zariski open set where ${\rm det}~Toep(i,j) \neq 0$
for all $1-n \leq i \leq j < n-1$ and $j-i$ even.
Then ${\rm det}~T$ is a nonconstant
affine function of $x_{n-1}$, and so for any choice
of $x_{1-n},...,x_{n-2}$ in $U$, ${\rm det}~T$ is zero
for a unique choice of $x_{n-1}$.
This shows that $V_{\R}({\rm det}~T)$
has real codimension 1. 

Furthermore, ${\rm det}~T$ has
highest order term in each $x_i$, $0<i \leq n-1$, equal to
$\pm Toep(1-n,2i-n-1) x_i^{n-i}$, i.e., with nonzero coefficient on $U$. It also has the 
highest order term in each $x_i$, $1-n \leq i <0$, equal to
$\pm Toep(n+2i+1,n-1) x_i^{n+i}$, i.e., with nonzero coefficient on $U$. Finally, the highest order term in $x_0$ is $x_0^n$, with coefficient 1.
Thus the gradient of ${\rm det}~T$ has all nonzero components on
a Zariski open set, and whether ${\rm det}~T = 0$ depends on all
variables.
\qed

\subsection{Vandermonde matrices and generalizations}
\label{vander}

In this section we will explain the entries filled ``\textbf{No}'' in Table~\ref{table1}. First we will show that polynomial Vandermonde matrices do not have algorithms for computing accurate inverses, by proving that certain minors needed in the expression cannot be accurately computed (this will also explain the ``\textbf{No}'' in the \emph{Any minor\/} column). Finally, we will show that the LDU factorization for polynomial Vandermonde matrices cannot be computed accurately. 

First we consider the class of generalized Vandermonde matrices $V$,
where $V_{ij} = P_{j-1}(x_i)$ is a polynomial function of $x_i$,
with $1 \leq i,j \leq n$.
This class includes the standard Vandermonde 
(where $P_{j-1}(x_i) = x_i^{j-1}$) and many others.

Consider a generalized Vandermonde matrix where
$P_{j-1}(x_i) = x_i^{j-1+ \lambda_{i}}$
with $0 \leq \lambda_1 \leq \lambda_2 \leq \cdots \leq \lambda_n$.
The tuple $\lambda = (\lambda_1 , \lambda_2 , ... , \lambda_n)$ is called a {\em partition.\/}
Any square submatrix of such a generalized Vandermonde matrix 
is also a generalized Vandermonde matrix.
A generalized Vandermonde matrix is known to have determinant of the form
$s_{\lambda}(x) \prod_{i<j} (x_i-x_j)$ where $s_{\lambda}(x)$
is a polynomial of degree $|\lambda| = \sum_i \lambda_i$,
and called a Schur function~\cite{macdonald}. 
In infinitely many variables (not our situation) the
Schur function is irreducible~\cite{farahat58}, but
in finitely many variables, the Schur function is sometimes irreducible
and sometimes not~\cite[Exer. 7.30]{EC2}.
But there are irreducible Schur functions of arbitrarily high degree.
Thus we conclude by Theorem~\ref{Thm_StructuredMatrixImpossibility}
that no enumerable set of black-box operations of bounded degree can compute all
Schur functions accurately when the $x_i$ are complex.

If we restrict the domain ${\cal D}$ to be nonnegative real numbers,
then the situation changes: The nonnegativity of the coefficients of
the Schur functions shows that they are positive in $\cal D$, and
indeed the generalized Vandermonde matrix is totally positive \cite{karlin}.
Combined with the homogeneity of the Schur function, 
Theorem~\ref{thm_positive_homo} implies that the Schur function,
and so determinants (and minors) of totally positive generalized
Vandermonde matrices can be evaluated accurately in classical arithmetic.
For accurate algorithms that are more efficient than the one in
Theorem~\ref{thm_positive_homo}, see~\cite{demmelkoevSchurJack}.





Now consider a polynomial Vandermonde matrix $V_P$ defined by a
family $\{P_k(x)\}_{k \in \mathbb{N}}$ of polynomials such that 
deg$(P_k) = k$, and $V_P(i,j) = P_{j-1}(x_i)$.

Note that any $V_P$ can be written as $V_P = V C$, with $V$ 
being a regular Vandermonde matrix, and $C$ being an upper 
triangular matrix of coefficients of the polynomials $P_k$, i.e., 
$$ P_{j-1}(x) = \sum_{i=1}^{j} C(i,j) x^{i-1}~, ~~\forall 1 \leq j \leq n~.
$$
Denote by $c_{i-1} := \tilde{D}(i,i)$, for all $1 \leq i \leq n$ 
the highest-order coefficients of the polynomials $P_0(x), \ldots, P_{n-1}(x)$.

To compute the inverse of the matrix $V_P$, we need to compute 
the minors that result from deleting a row and a column, 
i.e., (in MATLAB notation) 
det$(V_P([1\!:\!i\!-\!1, i\!+\!1\!:\!n], [1\!:\!j\!-\!1, j\!+\!1\!:\!n])$. 
We will focus our attention on the computation of the $(i, n-1)$ minors, 
i.e., the ones that result from deleting any of the rows and the $(n-1)$st 
column.

The resulting matrices look like
\[
M_{P,i} := V_P([1\!:\!i\!-\!1, i\!+\!1\!:\!n], [1\!:\!n\!-\!2, n]) 
= 
\left[ \begin{array}{ccccc}
c_0 & P_1(x_1) & \ldots & P_{n-3}(x_1)  & P_{n-1}(x_1) \\
\vdots & \vdots & \vdots & \vdots  & \vdots \\
c_0 & P_1(x_{i-1}) & \ldots & P_{n-3}(x_{i-1}) & P_{n-1}(x_{i-1}) \\
c_0 & P_1(x_{i+1}) & \ldots & P_{n-3}(x_{i+1}) & P_{n-1}(x_{i+1}) \\
\vdots & \vdots & \vdots & \vdots & \vdots \\
c_0 & P_1(x_n) & \ldots & P_{n-3}(x_n) & P_{n-1}(x_n) \end{array} \right ]~~.
 \]
Hence, we can manipulate the columns of det$(M_{P, i})$ 
by subtracting from them linear combinations of other columns, to obtain 
\begin{eqnarray*}
{\rm det}(M_{P,i}) &=& \left | \begin{array}{cccccc} 
c_0 & c_1 x_1 & c_2 x_1^2 & \ldots & c_{n-3} x_1^{n-3} & 
   c_{n-1} x_1^{n-1} + C(n -1, n) x_1^{n-2} \\
\vdots & \vdots & \vdots & \vdots & \vdots \\
c_0 & c_1 x_{i-1} & c_2 x_{i-1}^2 & \ldots & 
  c_{n-3} x_{i-1}^{n-3} & c_{n-1} x_{i-1}^{n-1} + C(n - 1, n) x_{i-1}^{n-2} \\
c_0 & c_1 x_{i+1} & c_2 x_{i+1}^2 & \ldots & 
   c_{n-3} x_{i+1}^{n-3} & c_{n-1} x_{i+1}^{n-1} + C(n - 1, n) x_{i+1}^{n-2} \\
\vdots & \vdots & \vdots & \vdots & \vdots \\
c_0 & c_1 x_n & c_2 x_n^2 & \ldots & c_{n-3} x_n^{n-3} & 
  c_{n-1} x_n^{n-1} + C(n-1,n) x_n^{n-2} \end{array} \right | 
\end{eqnarray*}

By expanding on the last column, and using the results from 
\cite{demmelkoev99}, we obtain that there are constants $E$ and $F$, 
specifically 
$E = C(n-1,n) \prod_{i=1}^{n-2} c_{i-1}$ and $F = E \frac{c_{n-1}}{C(n-1,n)}$, such that 
\[
{\rm det}(M_{P,i}) = 
\prod_{\mbox{\small{$\begin{array}{c} k<j \\ k, j \neq i \end{array}$}}} 
(x_j - x_k) ~\big [ E + F \cdot s_{[1]}(x_1, \ldots, x_{i-1}, x_{i+1}, \ldots, x_n) \big ]~,
\]
with $s_{[1]}$ being the Schur function corresponding to the 
partition $\lambda = (1, 0, \ldots, 0)$, i.e., 
\[
{\rm det}(M_{P,i}) = 
\!\!\prod_{\mbox{\small{$\begin{array}{c} k<j \\ k, j \neq i \end{array}$}}} 
(x_j - x_k) ~\big [ E~+~ F \cdot (x_1+\ldots +x_{i-1}+ x_{i+1}+ \ldots +x_n) \big ]~,
\]
and it is not hard to see that for any $n \geq 4$, the above polynomial in 
$x_1, \ldots, x_{i-1}, x_{i+1}, \ldots, x_n$ 
\emph{does not have an allowable variety,\/} and hence the inverse 
\emph{cannot be evaluated accurately} in classical arithmetic. 

Denote $\vec{x_{\neq i}} := (x_1, \ldots, x_{i-1}, x_{i+1}, \ldots, x_n)$.

Similarly, one can prove that the $(i, n-k)$ minor det $M_{P,i,k}$ can be obtained as
\[ \det(M_{P,i,k}) = \prod_{\begin{array}{c} k<j //k, k \neq i \end{array}} (x_j - x_k) ~\big[
A_1 +A_2 s_{[1]}(\vec{x_{\neq i}}) + \ldots + A_{k} s_{[1^k]}(\vec{x_{\neq i}}) \big ]~, \]
where $[1^l] = (1, 1, 1, \ldots, 0)$, the right side containing exactly $l$ ones and the rest $0$; $A_1, \ldots, A_k$ are constants which can be computed easily in terms of the entries of the matrix $C$.

Since for any $l$, $s_{[1^l]}$ is a homogeneous \emph{irreducible} function of degree $l$, the factor in the square brackets has degree $k$. Appropriate choices of $n$ and the matrix $C$ are likely to make this factor irreducible (for example, by making $|A_k| \gg |A_l|$ for all $l \neq k$). If this is the case, then by Theorem~\ref{Thm_StructuredMatrixImpossibility}, 
this family of matrices has inverses that 
cannot be evaluated accurately even with the addition of 
any enumerable set of bounded-degree black-boxes.



This explains why we have filled in with ``\textbf{No}'' 
the entries corresponding to columns ``$A^{-1}$'' and ``$Any$ $minor$'' 
in the Polynomial Vandermonde row of Table~\ref{table1}. 
Below we explain the ``\textbf{No}'' in the Polynomial Vandermonde row, in the column ``$LDU$''.

We can write $C = \tilde{D} \tilde{C}$, 
with $\tilde{D}$ being the diagonal matrix of highest-order coefficients, 
i.e., $\tilde{D}(i,i) = C(i,i)$ for all $1 \leq i \leq n$. 
We will assume that the matrices $C$ and $\tilde{D}$ 
are given to us exactly. 

If we let $V_P = L_P D_P U_P$ and $V = LDU$, it follows that 
\begin{eqnarray*}
L_P & = & L~;\\
D_P & = & D \tilde{D}~; \\
U_P & = & \tilde{D}^{-1} UC~.
\end{eqnarray*}

Since we cannot compute $L$ accurately in the general Vandermonde case, 
it follows that we cannot compute $L_P$ accurately in the polynomial 
Vandermonde case.

Finally, we explain the ``*'' entries in the polynomial Vandermonde row.
These depend on special properties of the polynomial.
In general, neither the SVD nor the symmetric eigenvalue decomposition (EVD)
are computable accurately, but if the polynomials are certain orthogonal
polynomials, then the accurate SVD is possible \cite{demmelkoevPolyVandermonde},
and an accurate symmetric EVD may also be possible \cite{dopicomoleramoro03}.


\section*{Acknowledgments}

We are grateful to Bernd Sturmfels, Gautam Bharali, Plamen Koev, William (Velvel) Kahan,
and Gregorio Malajovich for interesting discussions, to Jonathan Dorfman for his help 
with the graphs, and to the anonymous referee for helpful critique and suggestions. 


\bibliographystyle{plain}
\bibliography{linalg,biblio,sr,extra}

\end{document}